%% file: spinsl2paper.tex
\documentclass[twoside]{irmaems}
\usepackage{amsmath}
\usepackage{amssymb}
\usepackage{latexsym}

\usepackage{amsfonts,amscd,trace,calc,array}
\usepackage{tikz}
\usepackage{graphics,color}
\usetikzlibrary{snakes}
\usepackage[ps]{xy}\xyoption{all}
\usepackage{makeidx}

\newcommand{\C}{\mathbb{C}}
\newcommand{\N}{\mathbb{N}}
\def\cb{\C}
\newcommand{\xb}{\mathbf{x}}
\newcommand{\XX}{\mathfrak{X}}
\newcommand{\Z}{\mathbb{Z}}
\newcommand{\bfx}{\mathbf{x}}
\newcommand{\Id}{\mathbb{I}}
\newcommand{\id}{\mathbb{I}}
\newcommand{\bfy}{\mathbf{y}}
\newcommand{\bfz}{\mathbf{z}}

\newcommand{\Sym}{\mathrm{Sym}}
\newcommand{\End}{\mathrm{End}}
\newcommand{\Hom}{\mathrm{Hom}}
\newcommand{\tr}{\mathrm{tr}}
\def\gr{\mathsf{gr}}

\def\injects{\hookrightarrow}
\def\tensor{\otimes}
\def\isom{\cong}
\def\dmd{\diamondsuit}
\def\cross{\times}

\def\sss{\scriptscriptstyle}
\def\ssz{\scriptstyle}

\newcommand{\chxx}{\raisebox{2pt}{$\chi$}}

\newcommand{\fg}{\mathsf{F}_2}
\newcommand{\cgright}{\C[G]_{\sss\textup{right}}}
\newcommand{\cspan}{\C\,}
\def\mo{\!\!\!\mod}
\def\half{\frac12}
\def\like{\:\leftrightarrow\:}
\def\dtensor{\tensor\cdots\tensor}
\def\mtt{M_{2\times 2}}
\def\nth{$n$th }
\def\qimplies{\quad\implies\quad}
\newcommand{\sltc}{{\rm SL}(2,\C)}

\newcommand{\ch}[3]{\chxx^{\sss #1,#2,#3}}
\newcommand{\chx}[3]{\chxx_{\sss #1,#2,#3}}
\newcommand{\bs}[3]{{#1}^*_{#2}\otimes {#1}_{#3}}
\newcommand{\imx}[1]{\tmx{#1_{11} & #1_{12}}{#1_{21} & #1_{22}}}
\newcommand{\Adm}[3]{#3\in\lceil #1,#2\rfloor}

\newcommand{\inv}[1]{#1^{-1}}
\newcommand{\tmxt}[2]{\bigl[\begin{smallmatrix}#1\\#2\end{smallmatrix}\bigr]} 
\newcommand{\tmx}[2]{\begin{bmatrix}#1\\#2\end{bmatrix}}                      
\def\vprod#1{V^{\tensor #1}}
\def\pfrac#1#2{\left(\frac{#1}{#2}\right)}
\def\ptfrac#1#2{\left(\tfrac{#1}{#2}\right)}

\def\sixjt#1#2#3#4#5#6{\bigl[\begin{smallmatrix}#1 & #2 & #6\\#4 & #3 & #5\end{smallmatrix}\bigr]}

\input{tikz-defs}

\newtheorem{thm}{Theorem}[section]
\newtheorem{cor}[thm]{Corollary}
\newtheorem{prop}[thm]{Proposition}
\newtheorem{lem}[thm]{Lemma}

\newtheorem*{FKV}{Theorem \ref{fricke}}

\theoremstyle{definition}
\newtheorem{defn}[thm]{Definition}

\newtheorem{conv}[thm]{Convention}

\makeindex

\begin{document}

\title{Spin networks and $\mathrm{SL}(2,\mathbb{C})$-character varieties}

\author{Sean Lawton and Elisha Peterson}

\address{
Departamento de Matem\'{a}tica, Instituto Superior T\'ecnico, Lisbon, Portugal\\
email:\,\tt{slawton@math.ist.utl.pt}
\\[4pt]
Department of Mathematical Sciences, West Point, NY 10996\\
email:\,\tt{elisha.peterson@usma.edu}
}

\maketitle



\input{spinsl2-intro}

\ifnum1=1%
\section*{Acknowledgements}
We would like to thank Bill Goldman for introducing this problem to us and for many helpful suggestions, including generously sharing his {\it Mathematica} notebooks with us.  His correspondence with Nicolai Reshetikhin and Charles Frohman provided the foundation for the application of spin networks to this problem. Reshetikhin sketched proofs of both Theorems \ref{symmetry} and \ref{mult}. This work has benefitted from helpful conversations with Ben Howard, Tom Haines, and regular participation in the University of Maryland's Research Interaction Teams. The first author has received research support from John Millson, Richard Schwartz, and the University of Maryland's VIGRE grant. The second author has been supported by an NSF Graduate Fellowship. All diagrams in this chapter were generated by a suite of TikZ/PGF commands written by the second author.

We benefitted greatly from comments by Adam Sikora, Carlos Florentino, Athanase Papadopoulos, and the referee on early drafts of this chapter. Carlos Florentino provided many helpful corrections, and pointed out the significance of the symmetry in Theorem \ref{symmetry}. He also suggested the correspondence given in Proposition \ref{monic}. The referee gave valuable feedback regarding the organization and exposition.
\fi

\input{spinsl2-prelim-sean}
\input{spinsl2-prelim-elisha}

\input{spinsl2-decomposition}
\input{spinsl2-rank2}



\end{document}

%% file: tikz-defs.tex
\tikzstyle heightone=[scale=.7,xscale=.6,shift={(0,-.3)}]
\tikzstyle{every picture}=[semithick,baseline=0pt,heightone]

\tikzstyle heightones=[scale=.8,xscale=.35,shift={(0,.1)}]
\tikzstyle heightoneonehalf=[scale=.9,shift={(0,-.2)}]
\tikzstyle heighttwo=[scale=.8,shift={(0,-.4)}]
\tikzstyle heighttwos=[scale=.5,xscale=.6,shift={(0,-.1)}]
\tikzstyle heightthree=[scale=.6,shift={(0,-.9)}]
\tikzstyle heightthrees=[scale=.4,xscale=.7,shift={(0,-.2)}]

\tikzstyle trivalent=[very thick]
\tikzstyle dotedge=[pos=.5,shape=circle,fill=black,draw,minimum size=2pt,inner sep=0pt]

\tikzstyle matrix=[thick,circle,draw,fill=white,scale=.8,inner sep=1pt]
\tikzstyle smatrix=[thick,circle,draw,fill=white,scale=.6,inner sep=1pt]
\tikzstyle midmx=[matrix,pos=.5]
\tikzstyle midsmx=[smatrix,pos=.5]

\tikzstyle basiclabel=[draw=none,fill=none,shape=rectangle,inner sep=2pt,scale=.8]

\tikzstyle symlabel=[draw=none,fill=none,white,scale=.8]
\tikzstyle asymlabel=[draw=none,fill=none,black,scale=.8]
\tikzstyle symlabelleft=[basiclabel,left]
\tikzstyle symlabelright=[basiclabel,right]

\def\pup{\:\tikz[xscale=.7]{\tupp<(0,0)>{1}{1}}\:}
\def\px{\perm{0/1,1/0}}
\def\pid{\tikz[xscale=.7]{\tupp<(0,0)>{2}{1}}}
\def\pcc{\drawtl{0/1}{0/1}{}}
\def\binor{\px=\pid-\pcc}

\def\upcl{\tikz{
		\draw(0,0)..controls+(.1,.5)and+(-.1,-.5)..(0,1)node[dotedge](N){};
		\draw(N)--+(-.3,0);}}
\def\upcr{\tikz{
		\draw(0,0)..controls+(.1,.5)and+(-.1,-.5)..(0,1)node[dotedge](N){};
		\draw(N)--+(.3,0);}}
		
\def\cupp{\tikz{
	\draw(0,.7)..controls+(-.01,-.1)and+(.01,.1)..(0,.6)
		..controls+(-.04,-.4)and+(-.04,-.4)..(1,.6)..controls+(.01,.1)and+(-.01,-.1)..(1,.7);
}}
\def\cupdot{\tikz{
	\draw(0,.7)..controls+(-.01,-.1)and+(.01,.1)..(0,.6)
		..controls+(-.04,-.4)and+(-.04,-.4)..node[dotedge](N){}(1,.6)..controls+(.01,.1)and+(-.01,-.1)..(1,.7);
}}
\def\cupcu{\tikz{
	\draw(0,.7)..controls+(-.01,-.1)and+(.01,.1)..(0,.6)
		..controls+(-.04,-.4)and+(-.04,-.4)..node[dotedge](N){}(1,.6)..controls+(.01,.1)and+(-.01,-.1)..(1,.7);
		\draw(N)--+(0,.2);
}}
\def\cupcd{\tikz{
	\draw(0,.7)..controls+(-.01,-.1)and+(.01,.1)..(0,.6)
		..controls+(-.04,-.4)and+(-.04,-.4)..node[dotedge](N){}(1,.6)..controls+(.01,.1)and+(-.01,-.1)..(1,.7);
		\draw(N)--+(0,-.2);
}}
\def\cupcl{
	\tikz{\draw(0,.7)..controls+(-.01,-.1)and+(.01,.1)..(0,.5)
		..controls+(-.04,-.4)and+(-.04,-.4)..node[dotedge,pos=0](N){}(1,.5)..controls+(.01,.1)and+(-.01,-.1)..(1,.7);
		\draw(N)--+(-.3,0);
}}
		
\def\capp{\tikz{\draw(0,.3)..controls+(.01,.1)and+(-.01,-.1)..(0,.4)
		..controls+(.04,.4)and+(.04,.4)..(1,.4)..controls+(-.01,-.1)and+(.01,.1)..(1,.3);}}
\def\capdot{\tikz{\draw(0,.3)..controls+(.01,.1)and+(-.01,-.1)..(0,.4)
		..controls+(.04,.4)and+(.04,.4)..node[dotedge](N){}(1,.4)..controls+(-.01,-.1)and+(.01,.1)..(1,.3);}}
\def\capcu{\tikz{\draw(0,.3)..controls+(.01,.1)and+(-.01,-.1)..(0,.4)
		..controls+(.04,.4)and+(.04,.4)..node[dotedge](N){}(1,.4)..controls+(-.01,-.1)and+(.01,.1)..(1,.3);
		\draw(N)--+(0,.2);}}
\def\capcd{\tikz{\draw(0,.3)..controls+(.01,.1)and+(-.01,-.1)..(0,.4)
		..controls+(.04,.4)and+(.04,.4)..node[dotedge](N){}(1,.4)..controls+(-.01,-.1)and+(.01,.1)..(1,.3);
		\draw(N)--+(0,-.2);}}
		
\def\circb{\tikz{\draw(0,.5)..controls+(.05,.5)and+(.05,.5)..(1,.5)
		..controls+(-.05,-.5)and+(-.05,-.5)..(0,.5);}}
\def\circs{\tikz{\draw(0,.5)..controls+(.05,.5)and+(.05,.5)..(1,.5)
		..controls+(-.05,-.5)and+(-.05,-.5)..(0,.5);}}

\def\mbinor(#1,#2){
	\tikz[heighttwo]{
		\tperm<(0,0)>{0/1,1/0}
		\draw(0,1)..controls+(.04,.4)and+(-.04,-.4)..node[pos=.3,smatrix]{$#1$}(0,2);
		\draw(1,1)..controls+(.04,.4)and+(-.04,-.4)..node[pos=.3,smatrix]{$#2$}(1,2);
	}
	=
	\tikz[heighttwo]{
		\tperm<(0,0)>{0/0,1/1}
		\draw(0,1)..controls+(.04,.4)and+(-.04,-.4)..node[pos=.3,smatrix]{$#1$}(0,2);
		\draw(1,1)..controls+(.04,.4)and+(-.04,-.4)..node[pos=.3,smatrix]{$#2$}(1,2);
	}
	-
	\tikz[heighttwo]{
		\ttl<(0,0)>{0/1}{0/1}{}
		\draw(0,1)..controls+(.04,.4)and+(-.04,-.4)..node[pos=.3,smatrix]{$#1$}(0,2);
		\draw(1,1)..controls+(.04,.4)and+(-.04,-.4)..node[pos=.3,smatrix]{$#2$}(1,2);
	}
}
\def\mbinorrela(#1,#2){
	\tikz[heighttwo]{
		\draw(0,0)..controls+(.04,.4)and+(-.04,-.4)..node[midsmx]{$#1$}(0,2);
		\draw(1,0)..controls+(.04,.4)and+(-.04,-.4)..node[midsmx]{$#2$}(1,2);
	}
}
\def\mbinorrelb(#1){
	\tikz[heighttwo]{
		\ttl<(0,0)>{0/1}{0/1}{}
		\draw(0,1)..controls+(.04,.4)and+(-.04,-.4)..node[midsmx]{$#1$}(0,2);
		\draw(1,1)..controls+(.04,.4)and+(-.04,-.4)..(1,2);
	}
}

\def\capcux{\tikz[yscale=.7,xscale=.8,shift={(0,-.2)}]{
		\draw(0,1)..controls+(.05,.8)and+(.05,.8)..(1,1)node[dotedge](N){};
		\draw(N)--+(0,.3);
		\tperm<(0,0)>{0/1,1/0}
	}}
\def\kink{\tikz{
		\draw(0,0)..controls+(.02,.2)and+(-.02,-.2)..(0,.5)
			..controls+(.03,.3)and+(.03,.3)..(.5,.5)..controls+(-.03,-.3)and+(-.03,-.3)..(1,.5)
			..controls+(.02,.2)and+(-.02,-.2)..(1,1);}}
\def\kinkdot{\tikz{
		\draw(0,0)..controls+(.02,.2)and+(-.02,-.2)..(0,.5)
			 ..controls+(.03,.3)and+(.03,.3)..node[dotedge]{}(.5,.5)..controls+(-.03,-.3)and+(-.03,-.3)..node[dotedge]{}(1,.5)
			..controls+(.02,.2)and+(-.02,-.2)..(1,1);}}
\def\ccdot{\tikz{
	\draw(0,0)..controls+(.1,.5)and+(.1,.5)..node[dotedge]{}(1,0);
	\draw(0,1)..controls+(-.1,-.5)and+(-.1,-.5)..node[dotedge]{}(1,1);}}
		
\def\cupsxlr{
	\tikz[heightthree,scale=.8,shift={(0,-.2)}]{\tperm<(0,2)>{0/1,1/0,3/3,4/4}\tcupp<(2,2)>{2}}
	=\tikz[heightthree,scale=.8,shift={(0,-.2)},xscale=-1]{\tperm<(0,2)>{0/1,1/0,3/3,4/4}\tcupp<(2,2)>{2}}}
	
\def\upcll{
	\tikz{
		\draw(0,0)..controls+(.02,.2)and+(-.02,-.2)..node[dotedge](M){}(0,.5)
			..controls+(.02,.2)and+(-.02,-.2)..node[dotedge](N){}(0,1);
		\draw(M)--+(-.4,0)(N)--+(-.4,0);}
	=\tikz{\draw(0,0)..controls+(.05,.5)and+(-.05,-.5)..(0,1){};}
}

\def\triplel{
	\tikz[heighttwo,xscale=.8]{
		\draw
			(0,0)..controls+(.05,.5)and+(-.05,-.5)..(-1,1)..controls+(.05,.5)and+(-.05,-.5)..(0,2)
			(1,0)..controls+(.05,.5)and+(-.05,-.5)..(2,1)..controls+(.03,.3)and+(.03,.3)..node[dotedge](N){}(3,1)
				..controls+(-.06,-.6)and+(.06,.6)..(3,0)
			(4,0)..controls+(.06,.6)and+(-.06,-.6)..(4,1)..controls+(.1,1)and+(.1,1)..node[dotedge](M){}(1,1)
				..controls+(-.03,-.3)and+(-.03,-.3)..node[dotedge](L){}(0,1)..controls+(.05,.5)and+(-.05,-.5)..(1,2);
		\draw(L)--+(0,.3)(M)--+(0,.3)(N)--+(0,.3);}
	=-\tikz[heighttwo,xscale=.8]{
		\draw
			(0,0)..controls+(.08,.8)and+(-.08,-.8)..(2,2)
			(5,0)..controls+(.08,.8)and+(-.08,-.8)..(3,2)
			(1,0)..controls+(.1,1)and+(.1,1)..node[dotedge](N){}(4,0);
		\draw(N)--+(0,.4);}
}

\def\ims#1{\tikz{\draw(0,0)..controls+(.05,.5)and+(-.05,-.5)..node[midsmx]{$#1$}(0,1);}}

\def\cupml#1{\tikz{
		\draw(0,1)..controls+(-.03,-.3)and+(.03,.3)..node[midmx]{$#1$}(0,.3)
			..controls+(-.05,-.5)and+(-.05,-.5)..(1,.3)..controls+(.03,.3)and+(-.03,-.3)..(1,1);}}
\def\cupmr#1{\tikz{
		\draw(0,1)..controls+(-.03,-.3)and+(.03,.3)..(0,.3)
			..controls+(-.05,-.5)and+(-.05,-.5)..(1,.3)..controls+(.03,.3)and+(-.03,-.3)..node[midmx]{$#1$}(1,1);}}
\def\cupmm#1#2{\tikz{
		\draw(0,1)..controls+(-.03,-.3)and+(.03,.3)..node[midmx]{$#1$}(0,.3)
			..controls+(-.05,-.5)and+(-.05,-.5)..(1,.3)..controls+(.03,.3)and+(-.03,-.3)..node[midmx]{$#2$}(1,1);}}

\def\capmm#1#2{\tikz{
		\draw(0,0)..controls+(.03,.3)and+(-.03,-.3)..node[midmx]{$#1$}(0,.7)
			..controls+(.05,.5)and+(.05,.5)..(1,.7)..controls+(-.03,-.3)and+(.03,.3)..node[midmx]{$#2$}(1,0);}}

\def\circml#1{\tikz{
		\draw(0,.7)..controls+(-.03,-.3)and+(.03,.3)..node[midsmx]{$#1$}(0,.3)
			..controls+(-.05,-.5)and+(-.06,-.6)..(1,.5)..controls+(.06,.6)and+(.05,.5)..(0,.7);}}

\def\circmm#1#2{\tikz{
		\draw(0,.7)..controls+(-.03,-.3)and+(.03,.3)..node[midsmx]{$#1$}(0,.3)
			..controls+(-.05,-.5)and+(-.05,-.5)..(1,.3)..controls+(.03,.3)and+(-.03,-.3)..node[midsmx]{$#2$}(1,.7)
			..controls+(.05,.5)and+(.05,.5)..(0,.7);}}

\def\upcln#1{\tikz[xscale=.7]{
		\foreach\xa in{0,.5,2}{\draw(\xa,0)..controls+(.1,.5)and+(-.1,-.5)..node[dotedge](N){}(\xa,1);\draw(N)--+(-.4,0);}
		\foreach\xa/\xb in{1/.2,1.5/.2,1/.8,1.5/.8}{\node[basiclabel,scale=1.5]at(\xa,\xb){.};}
		\node[basiclabel,scale=1.3,right]at(2,.9){$#1$};}}

\def\cupcun#1{
	\tikz[scale=.8,shift={(0,.2)}]{
		\draw(-.5,1)..controls+(-.02,-.2)and+(.02,.2)..(-.5,.7)..controls+(-.05,-.5)and+(-.05,-.5)
			..node[dotedge](M){}(.5,.7)..controls+(.02,.2)and+(-.02,0)..(.5,1);
		\draw(-2,1)..controls+(0,-.1)and+(0,.1)..(-2,.7)..controls+(-.15,-1.5)and+(-.15,-1.5)
			..node[dotedge](N){}(2,.7)..controls+(0,.1)and+(0,-.1)..(2,1)node[symlabelright]{$#1$};
		\draw(M)--+(0,.3)(N)--+(0,.3);
		\foreach\xa/\xb in{-1.2/0,-.7/.2,1.1/0,.6/.2}{\node[basiclabel,scale=2]at(\xa,\xb){.};}
	}}
		
\def\circn#1{
	\tikz[heightthrees]{
		\ttl<(0,0)>{}{-5/-4}{}\tperm<(0,1)>{-5/-5,-4/-4,-1/-1}\ttl<(0,2)>{-5/-4}{}{}
		\draw(-1,2)..controls+(.2,1.5)and+(.2,1.5)..(-8,2)..controls+(-.05,-.5)and+(.05,.5)..(-8,1)
			..controls+(-.2,-1.5)and+(-.2,-1.5)..(-1,1);
		\filldraw[black](-4.5,1.3)rectangle+(4,.4);
		\foreach\xa/\xb in{-7/1.5,-6/1.5,-3.2/2.3,-2.2/2.4,-3.2/.7,-2.2/.6}{\node[basiclabel,scale=1.5]at(\xa,\xb){.};}
		\node[basiclabel,right]at(-.5,1.7){$#1$};}}

\def\mupn#1#2{\tikz[xscale=.8]{
		\foreach\xa in{0,.5,2}{\draw(\xa,0)..controls+(.1,.5)and+(-.1,-.5)..node[midsmx]{$#2$}(\xa,1);}
		\node[basiclabel,right]at(2,1){$#1$};
		\foreach\xa/\xb in{1/.2,1.5/.2,1/.8,1.5/.8}{\node[basiclabel,scale=2]at(\xa,\xb){.};}}}

\def\slantn#1{\tikz{
		\draw(0,0)..controls+(.04,.4)and+(-.04,-.4)..(.9,1);
		\draw(.6,0)..controls+(.04,.4)and+(-.04,-.4)..(1.5,1)node[basiclabel,right]{$#1$};
		\foreach\xa/\xb in{.65/.5,.85/.5}{\node[basiclabel,scale=1.5]at(\xa,\xb){.};}
}}
\def\kinkn#1{\tikz{
		\draw(.2,-.1)..controls+(.03,.3)and+(-.03,-.3)..(0,.5)
			..controls+(.04,.4)and+(.04,.4)..(1.2,.6)..controls+(-.01,-.1)and+(-.01,-.1)..(1.5,.6)
			..controls+(.03,.3)and+(-.03,-.3)..(1.3,1.1);
		\draw(.8,-.1)..controls+(.03,.3)and+(-.03,-.3)..(.6,.4)
			..controls+(.01,.1)and+(.01,.1)..(.9,.4)..controls+(-.04,-.4)and+(-.04,-.4)..(2.1,.5)
			..controls+(.03,.3)and+(-.03,-.3)..(1.9,1.1);
		\foreach\xa/\xb in{.3/.2,.5/.2,1.6/.8,1.8/.8,.75/.6,.75/.75,1.35/.4,1.35/.25}{
			\node[basiclabel,scale=1.5]at(\xa,\xb){.};}
		\node[basiclabel,right]at(1.9,1.1){$#1$};}}

\def\suppn#1{\tikz[xscale=.5]{
		\foreach\xa in{0,.5,2}{\draw(\xa,0)..controls+(.1,.5)and+(-.1,-.5)..(\xa,1);}
		\filldraw[black](-.3,.3)rectangle(2.3,.7);%
		\draw(1,.5)node[symlabel]{$#1$};
		\foreach\xa/\xb in{1/.1,1.5/.1,1/.9,1.5/.9}{\node[basiclabel,scale=2]at(\xa,\xb){.};}}}

\def\suppnbig#1{
	\tikz[heightthrees]{
		\tperm<(0,0)>{0/0,3/3}\tperm<(0,1)>{0/0,3/3}\tperm<(0,2)>{0/0,3/3}
		\filldraw[black](-.5,1.3)rectangle+(4,.4)node[symlabelright,scale=.8]{$#1$};
		\foreach\xa/\xb in{1/.5,2/.5,1/2.5,2/2.5}{\node[basiclabel,scale=1.5]at(\xa,\xb){.};}
}}
		
\def\cfupbig#1#2#3{\tikz[scale=.2,xscale=1.3,shift={(0,2)}]{
		\draw
			(-4,-5)--(-4,-4)..controls+(.5,3)and+(-.5,-3)..(-8,4)--(-8,5)
			(-1,-5)--(-1,-4)..controls+(.5,3)and+(-.5,-3)..(-5,4)--(-5,5)
			(1,-5)--(1,-4)..controls+(.5,3)and+(-.5,-3)..(5,4)--(5,5)
			(4,-5)--(4,-4)..controls+(.5,3)and+(-.5,-3)..(8,4)--(8,5)
			(-4,5)--(-4,4)..controls+(-.5,-5)and+(-.5,-5)..(4,4)--(4,5)
			(-1,5)--(-1,4)..controls+(-.3,-2)and+(-.3,-2)..(1,4)--(1,5);			
		\filldraw[black](-4.5,-4.5)rectangle+(9,.9)(-8.5,3.6)rectangle+(8,.9)(.5,3.6)rectangle+(8,.9);
		 \draw(-8.5,4.2)node[basiclabel,left]{$#1$}(8.5,4.2)node[basiclabel,right]{$#2$}(4.5,-4.2)node[basiclabel,right]{$#3$};
		\foreach\xa/\xb in{-5/0,-4/0,4/0,5/0,0/1,0/2}{\node[basiclabel,scale=1.5]at(\xa,\xb){.};}
}}

\def\cfdownbig#1#2#3{\tikz[scale=.2,xscale=1.3,shift={(0,2)},yscale=-1]{
		\draw
			(-4,-5)--(-4,-4)..controls+(.5,3)and+(-.5,-3)..(-8,4)--(-8,5)
			(-1,-5)--(-1,-4)..controls+(.5,3)and+(-.5,-3)..(-5,4)--(-5,5)
			(1,-5)--(1,-4)..controls+(.5,3)and+(-.5,-3)..(5,4)--(5,5)
			(4,-5)--(4,-4)..controls+(.5,3)and+(-.5,-3)..(8,4)--(8,5)
			(-4,5)--(-4,4)..controls+(-.5,-5)and+(-.5,-5)..(4,4)--(4,5)
			(-1,5)--(-1,4)..controls+(-.3,-2)and+(-.3,-2)..(1,4)--(1,5);			
		\filldraw[black](-4.5,-4.5)rectangle+(9,.9)(-8.5,3.6)rectangle+(8,.9)(.5,3.6)rectangle+(8,.9);
		 \draw(-8.5,4.2)node[basiclabel,left]{$#1$}(8.5,4.2)node[basiclabel,right]{$#2$}(4.5,-4.2)node[basiclabel,right]{$#3$};
		\foreach\xa/\xb in{-5/0,-4/0,4/0,5/0,0/1,0/2}{\node[basiclabel,scale=1.5]at(\xa,\xb){.};}
	}}

\def\cfdownpartial(#1,#2)(#3,#4){
	\tikz[heightthree,scale=.3,yscale=-.7,shift={(0,-8)}]{
		\draw
			(-5.5,-2)..controls+(.5,3)and+(-.5,-3)..(-8,4)--(-8,5)
			(-2.5,-2)..controls+(.5,3)and+(-.5,-3)..(-5,4)--(-5,5)
			(2.5,-2)..controls+(.5,3)and+(-.5,-3)..(5,4)--(5,5)
			(5.5,-2)..controls+(.5,3)and+(-.5,-3)..(8,4)--(8,5)
			(-4,5)--(-4,4)..controls+(-.5,-5)and+(-.5,-5)..(4,4)--(4,5)
			(-1,5)--(-1,4)..controls+(-.3,-2)and+(-.3,-2)..(1,4)--(1,5);			
		\filldraw[black](-8.5,3.6)rectangle+(8,.9)(.5,3.6)rectangle+(8,.9);
		\draw(-8.5,4.2)node[basiclabel,left]{$#1$}(8.5,4.2)node[basiclabel,right]{$#2$};
		\foreach\xa/\xb in{-5/0,-4/0,4/0,5/0,0/1,0/2}{\node[basiclabel,scale=1.5]at(\xa,\xb){.};}
		\draw[snake=brace,raise snake=-3pt]
			(-6,-3.5)--(-2,-3.5)node[basiclabel,pos=.5,above]{$#3$}
			(2,-3.5)--(6,-3.5)node[basiclabel,pos=.5,above]{$#4$};
	}}

\def\supp#1{\tikz[xscale=.7]{\tsupp<(0,0)>{#1}{1}}}
\def\asupp#1]{\tikz[xscale=.7{\tasupp<(0,0)>{#1}{1}}}

\def\perm#1{\tikz[xscale=.7]{\tperm<(0,0)>{#1}}}

\def\drawtl#1#2#3{\tikz[xscale=.7]{\ttl<(0,0)>{#1}{#2}{#3}}}

\def\tup(#1){\tikz{
		\draw[trivalent](0,0)..controls+(.1,.5)and+(-.1,-.5)..(0,1)node[symlabelright]{$#1$};
	}}

\def\tupc(#1){
	\tikz{
		\draw[trivalent](0,0)..controls+(.1,.5)and+(-.1,-.5)..node[dotedge](N){}(0,1)node[symlabelright]{$#1$};
		\draw(N)--+(-.4,0);
	}}

\def\tcup(#1){\tikz{
	\draw[trivalent](0,.7)..controls+(-.01,-.1)and+(.01,.1)..(0,.6)
		..controls+(-.04,-.4)and+(-.04,-.4)..(1,.6)..controls+(.01,.1)and+(-.01,-.1)..(1,.7)node[symlabelright]{$#1$};
	}}
	
\def\tcupc(#1){\tikz{
	 \draw[trivalent](0,.7)..controls+(-.01,-.1)and+(.01,.1)..(0,.6)..controls+(-.04,-.4)and+(-.04,-.4)..node[dotedge](N){}(1,.6)
		..controls+(.01,.1)and+(-.01,-.1)..(1,.7)node[symlabelright]{$#1$};
	\draw(N)--+(0,.2);
	}}

\def\tx(#1,#2){
	\tikz{
		\draw[trivalent]
			(1,0)..controls+(.1,.5)and+(-.1,-.5)..(0,1)node[symlabelleft]{$#1$}
			(0,0)..controls+(.1,.5)and+(-.1,-.5)..(1,1)node[symlabelright]{$#2$};
	}}

\def\txkink(#1,#2){
	\tikz[heighttwo]{
		\draw[trivalent]
			(.75,.4)..controls+(.1,.5)and+(-.1,-.5)..(1.25,1.6)node[symlabelleft]{$#1$}
			 (0,.4)..controls+(.1,1)and+(-.5,.7)..(1,1)..controls+(.5,-.7)and+(-.1,-1)..(2,1.6)node[symlabelright]{$#2$};
	}}
	
\def\tkink(#1){
	\tikz[heighttwo]{
		\draw[trivalent](-.5,0)..controls+(.05,.5)and+(-.05,-.5)..(-.8,1.3)..controls+(.1,.4)and+(.1,.6)..(0,1)
			 ..controls+(-.1,-.6)and+(-.1,-.4)..(.8,.7)..controls+(.05,.5)and+(-.05,-.5)..(.5,2)node[symlabelright]{$#1$};
	}}

\def\tfancyup(#1,#2){\tikz{
		 \draw[trivalent](0,0)..controls+(.05,.3)and+(-.05,-.3)..(.3,.5)..controls+(.05,.3)and+(-.05,-.3)..(0,1)node[symlabelleft]{$#1$};
		 \draw[trivalent](1,0)..controls+(.05,.3)and+(-.05,-.3)..(.7,.5)..controls+(.05,.3)and+(-.05,-.3)..(1,1)node[symlabelright]{$#2$};
	}}

\def\tbubble(#1,#2,#3,#4){
	\tikz[heighttwo]{
		\draw(0,1)node[coordinate](A){}(1,1)node[coordinate](B){}
			 (.5,.35)node[coordinate](R){}(.5,0)node[coordinate](Z){}
			 (.5,1.65)node[coordinate](RR){}(.5,2)node[coordinate](ZZ){};
		\foreach\xa/\xb/\xc/\xd/\xe in{R/A/.3//,R/B/.3/#2/right,Z/R/.2/#4/left,A/RR/.3/#1/left,B/RR/.3//,RR/ZZ/.2/#3/right}{
			\draw[trivalent](\xa)..controls+(0,\xc)and+(0,-\xc)..node[basiclabel,\xe,inner sep=3pt]{$\xd$}(\xb);}
	}}

\def\tcirc(#1){
	\tikz[heighttwo,xscale=1.5,shift={(0,1)}]{
		\draw[trivalent](0,0)..controls+(.1,1)and+(.1,1)..node[near end,right,basiclabel]{#1}(1,0)..controls+(-.1,-1)and+(-.1,-1)..(0,0);
	}}

\def\tbubbled(#1,#2,#3){
	\tikz[heighttwo]{
		\draw(0,1)node[coordinate](A){}(1,1)node[coordinate](B){}
			 (.5,.4)node[coordinate](R){}(.5,1.6)node[coordinate](RR){}(1.8,1)node[coordinate](Z){};
		\foreach\xa/\xb/\xc/\xd/\xe in{R/A/.3//,R/B/.3/#2/right,A/RR/.3/#1/left,B/RR/.3//}{
			\draw[trivalent](\xa)..controls+(0,\xc)and+(0,-\xc)..node[basiclabel,\xe,inner sep=3pt]{$\xd$}(\xb);}
	\draw[trivalent](RR)..controls+(0,.3)and+(0,1)..node[near end,right,basiclabel,inner sep=2pt]{#3}(Z)..controls+(0,-1)and+(0,-.3)..(R);
	}}
	
\def\tupy(#1,#2,#3){\tikz{
		 \draw(0,1)node[coordinate](A){}(1,1)node[coordinate](B){}(.5,.35)node[coordinate](R){}(.5,0)node[coordinate](Z){};
		\foreach\xa/\xb/\xc/\xd/\xe in{R/A/.3/$#1$/left,R/B/.3/$#2$/right,R/Z/-.2/$#3$/right}{
			\draw[trivalent](\xa)..controls+(.03,\xc)and+(-.03,-\xc)..(\xb)node[basiclabel,\xe]{\xd};}
	}}
	
\def\tdowny(#1,#2,#3){\tikz{
		 \draw(0,0)node[coordinate](A){}(1,0)node[coordinate](B){}(.5,.65)node[coordinate](R){}(.5,1)node[coordinate](Z){};
		\foreach\xa/\xb/\xc/\xd/\xe in{R/A/.3/$#1$/left,R/B/.3/$#2$/right,R/Z/-.2/$#3$/right}{
			\draw[trivalent](\xa)..controls+(-.03,-\xc)and+(.03,\xc)..(\xb)node[basiclabel,\xe]{\xd};}
	}}

\def\tupyx(#1,#2,#3){
	\tikz[heightoneonehalf]{
		 \draw(0,.8)node[coordinate](A){}(1,.8)node[coordinate](B){}(.5,.3)node[coordinate](R){}(.5,0)node[coordinate](Z){}
			(0,1.5)node[coordinate](AA){}(1,1.5)node[coordinate](BB){};
		\foreach\xa/\xb/\xc/\xd/\xe in{R/A/.3//left,R/B/.3//right,R/Z/-.2/#3/right,A/BB/.2/#2/right,B/AA/.2/#1/left}{
			\draw[trivalent](\xa)..controls+(.03,\xc)and+(-.03,-\xc)..(\xb)node[basiclabel,\xe]{$\xd$};}
	}}

\def\tdownyx(#1,#2,#3){
	\tikz[heightoneonehalf]{
		 \draw(0,.7)node[coordinate](A){}(1,.7)node[coordinate](B){}(.5,1.2)node[coordinate](R){}(.5,1.5)node[coordinate](Z){}
			(0,0)node[coordinate](AA){}(1,0)node[coordinate](BB){};
		\foreach\xa/\xb/\xc/\xd/\xe in{R/A/.3//left,R/B/.3//right,R/Z/-.2/#3/right,A/BB/.2/#2/right,B/AA/.2/#1/left}{
			\draw[trivalent](\xa)..controls+(-.03,-\xc)and+(.03,\xc)..(\xb)node[basiclabel,\xe]{$\xd$};}
	}}

\def\tupykink(#1,#2,#3){
	\tikz{
		\draw[trivalent]
			(.5,0)node[symlabelleft]{#1}..controls+(.03,.2)and+(-.03,-.2)..(.5,.3)
			 (.5,.3)..controls+(.02,.2)and+(-.02,-.2)..(.2,.6)..controls+(.02,.2)and+(-.02,-.2)..(.4,1)node[symlabelleft]{#3}
			 (.5,.3)..controls+(.02,.2)and+(-.1,-.3)..(.7,.7)..controls+(.1,.3)and+(-.05,1)..(1.5,0)node[symlabelright]{#2};
	}}

\def\tuptreel(#1,#2,#3,#4,#5){
	\tikz[heightoneonehalf]{
		\draw(0,1.5)node[coordinate](A){}(1,1.5)node[coordinate](B){}(2,1.5)node[coordinate](C){}
			 (.5,1)node[coordinate](R){}(1,.5)node[coordinate](S){}(1,0)node[coordinate](Z){};
		\foreach\xa/\xb/\xc/\xd/\xe in{R/A/.3/#1/left,R/B/.3/#2/right,S/C/.5/#3/right,S/Z/-.2/#4/right}{
			\draw[trivalent](\xa)..controls+(.03,\xc)and+(-.03,-\xc)..(\xb)node[basiclabel,\xe]{$\xd$};}
		\draw[trivalent](S)..controls+(.03,.3)and+(-.03,-.3)..node[basiclabel,left]{$#5\:$}(R);
	}}
\def\tuptreer(#1,#2,#3,#4,#5){
	\tikz[heightoneonehalf,xscale=-1]{
		\draw(0,1.5)node[coordinate](A){}(1,1.5)node[coordinate](B){}(2,1.5)node[coordinate](C){}
			 (.5,1)node[coordinate](R){}(1,.5)node[coordinate](S){}(1,0)node[coordinate](Z){};
		\foreach\xa/\xb/\xc/\xd/\xe in{R/A/.3/#3/right,R/B/.3/#2/right,S/C/.5/#1/left,S/Z/-.2/#4/right}{
			\draw[trivalent](\xa)..controls+(.03,\xc)and+(-.03,-\xc)..(\xb)node[basiclabel,\xe]{$\xd$};}
		\draw[trivalent](S)..controls+(.03,.3)and+(-.03,-.3)..node[basiclabel,right]{$\:#5$}(R);
	}}

\def\tdowntreel(#1,#2,#3,#4,#5){
	\tikz[yscale=-1,heightoneonehalf]{
		\draw(0,1.5)node[coordinate](A){}(1,1.5)node[coordinate](B){}(2,1.5)node[coordinate](C){}
			 (.5,1)node[coordinate](R){}(1,.5)node[coordinate](S){}(1,0)node[coordinate](Z){};
		\foreach\xa/\xb/\xc/\xd/\xe in{R/A/.3/#1/left,R/B/.3/#2/right,S/C/.5/#3/right,S/Z/-.2/#4/right}{
			\draw[trivalent](\xa)..controls+(.03,\xc)and+(-.03,-\xc)..(\xb)node[basiclabel,\xe]{$\xd$};}
		\draw[trivalent](S)..controls+(.03,.3)and+(-.03,-.3)..node[basiclabel,left]{$#5\:$}(R);
	}}
\def\tdowntreer(#1,#2,#3,#4,#5){
	\tikz[yscale=-1,heightoneonehalf,xscale=-1]{
		\draw(0,1.5)node[coordinate](A){}(1,1.5)node[coordinate](B){}(2,1.5)node[coordinate](C){}
			 (.5,1)node[coordinate](R){}(1,.5)node[coordinate](S){}(1,0)node[coordinate](Z){};
		\foreach\xa/\xb/\xc/\xd/\xe in{R/A/.3/#3/right,R/B/.3/#2/right,S/C/.5/#1/left,S/Z/-.2/#4/right}{
			\draw[trivalent](\xa)..controls+(.03,\xc)and+(-.03,-\xc)..(\xb)node[basiclabel,\xe]{$\xd$};}
		\draw[trivalent](S)..controls+(.03,.3)and+(-.03,-.3)..node[basiclabel,right]{$\:#5$}(R);
	}}

\def\tfused(#1,#2,#3,#4,#5){
	\tikz{
		\draw[trivalent]
			(0,-.1)node[symlabelleft]{$#4$}..controls+(0,.3)and+(0,-.3)..(.5,.4)
			(.5,.6)..controls+(0,.3)and+(0,-.3)..(0,1.1)node[symlabelleft]{$#1$}
			(.5,.4)..controls+(.02,.1)and+(-.02,-.1)..node[basiclabel,right]{$#5$}(.5,.6)
			(1,-.1)node[symlabelright]{$#3$}..controls+(0,.3)and+(0,-.3)..(.5,.4)
			(.5,.6)..controls+(0,.3)and+(0,-.3)..(1,1.1)node[symlabelright]{$#2$};
	}}

\def\thleft(#1,#2,#3,#4,#5){
	\tikz{
		 \draw(0,1)node[coordinate](A){}(1.7,1)node[coordinate](B){}(1.7,0)node[coordinate](C){}(0,0)node[coordinate](D){}
			(.3,.7)node[coordinate](R){}(1.4,.3)node[coordinate](S){};
		\foreach\xa/\xb/\xc/\xd/\xe in{R/A/.2/#1/left,S/B/.3/#2/right,S/C/-.2/#3/right,R/D/-.3/#4/left}{
			\draw[trivalent](\xa)..controls+(.02,\xc)and+(-.02,-\xc)..(\xb)node[basiclabel,\xe]{$\xd$};}
		\draw[trivalent](S)..controls+(.02,.2)and+(-.02,-.2)..node[basiclabel,above]{$#5$}(R);
	}}

\def\thright(#1,#2,#3,#4,#5){
	\tikz{
		 \draw(0,1)node[coordinate](A){}(1.7,1)node[coordinate](B){}(1.7,0)node[coordinate](C){}(0,0)node[coordinate](D){}
			(.3,.3)node[coordinate](R){}(1.4,.7)node[coordinate](S){};
		\foreach\xa/\xb/\xc/\xd/\xe in{R/A/.3/#1/left,S/B/.2/#2/right,S/C/-.3/#3/right,R/D/-.2/#4/left}{
			\draw[trivalent](\xa)..controls+(.02,\xc)and+(-.02,-\xc)..(\xb)node[basiclabel,\xe]{$\xd$};}
		\draw[trivalent](R)..controls+(.02,.2)and+(-.02,-.2)..node[basiclabel,above]{$#5$}(S);
	}}

\def\thtwisted(#1,#2,#3,#4,#5){
	\tikz[heighttwo,xscale=1.2]{
		\draw[trivalent](0,.5)..controls+(0,.3)and+(0,-.3)..(.5,.9)
									  (.5,1.1)..controls+(0,.3)and+(0,-.3)..(0,1.5);
		\draw[trivalent](.5,.9)..controls+(.02,.1)and+(-.02,-.1)..node[basiclabel,right]{$#5$}(.5,1.1);
		\draw[trivalent](1,.5)..controls+(0,.3)and+(0,-.3)..(.5,.9)
										(.5,1.1)..controls+(0,.3)and+(0,-.3)..(1,1.5);
		\draw[trivalent](0,0)node[symlabelright]{$#3$}..controls+(.02,.1)and+(-.02,-.1)..(0,.5)
										(1,1.5)..controls+(.02,.1)and+(-.02,-.1)..(1,2)node[symlabelleft]{$#1$}
										 (-.5,0)node[symlabelleft]{$#4$}..controls+(.1,.5)and+(-.1,-.5)..(-.5,1.5)..controls+(0,.5)and+(0,.5)..(0,1.5)
										 (1.5,2)node[symlabelright]{$#2$}..controls+(-.1,-.5)and+(.1,.5)..(1.5,.5)..controls+(0,-.5)and+(0,-.5)..(1,.5);
	}}
	
\def\thtwistedx(#1,#2,#3,#4,#5){
	\tikz[heighttwo,xscale=1.2]{
		\draw[trivalent](.2,.5)..controls+(0,.3)and+(0,-.3)..(.5,.9)
									  (.5,1.1)..controls+(0,.3)and+(0,-.3)..(0,1.5);
		\draw[trivalent](.5,.9)..controls+(.02,.1)and+(-.02,-.1)..node[basiclabel,right]{$#5$}(.5,1.1);
		\draw[trivalent](.7,.5)..controls+(0,.3)and+(0,-.3)..(.5,.9)
										(.5,1.1)..controls+(0,.3)and+(0,-.3)..(1,1.5);
		\draw[trivalent](0,0)node[symlabelright]{$#3$}..controls+(.02,.2)and+(-.03,-.3)..(.7,.5)
										(1,1.5)..controls+(.02,.1)and+(-.02,-.1)..(1,2)node[symlabelleft]{$#1$}
										 (-.5,0)node[symlabelleft]{$#4$}..controls+(.1,.5)and+(-.1,-.5)..(-.5,1.5)..controls+(0,.5)and+(0,.5)..(0,1.5)
										 (1.5,2)node[symlabelright]{$#2$}..controls+(-.1,-.5)and+(.1,.5)..(1.5,.5)..controls+(0,-.5)and+(0,-.4)..(.2,.5);
	}}

\def\thkink(#1,#2,#3,#4,#5){
	\tikz{
		 \draw(0,1)node[coordinate](A){}(1.5,1)node[coordinate](B){}(1.5,0)node[coordinate](C){}(0,0)node[coordinate](D){}
			(.3,.7)node[coordinate](R){}(1.2,.3)node[coordinate](S){};
		\foreach\xa/\xb/\xc/\xd/\xe in{R/A/.2/#1/left,S/B/.3/#2/right,S/C/-.2/#3/right,R/D/-.3/#4/left}{
			\draw[trivalent](\xa)..controls+(.02,\xc)and+(-.02,-\xc)..(\xb)node[basiclabel,\xe]{$\xd$};}
		\draw[trivalent](R)..controls+(.4,-1.3)and+(-.4,1.3)..node[basiclabel,above,inner sep=6pt]{$#5$}(S);
	}}
	
\def\tmup(#1)(#2){
	\tikz{
		\draw[trivalent](0,0)..controls+(.05,.2)and+(-.05,-.2)..(0,.5)node[smatrix]{$#2$}
			..controls+(.05,.2)and+(-.05,-.2)..(0,1)node[symlabelright]{$#1$};
	}}

\def\trela(#1,#2,#3)(#4,#5,#6){
	\tikz[heighttwo]{
		\draw(-2,2)node[coordinate](A){}(-.7,2)node[coordinate](AA){}
				 (.7,2)node[coordinate](B){}(2,2)node[coordinate](BB){}
				 (-.7,0)node[coordinate](C){}(.7,0)node[coordinate](CC){}
				 (-.7,.8)node[coordinate](R){}(.7,.8)node[coordinate](RR){};
		\draw[trivalent]
			(C)..controls+(.03,.5)and+(-.03,-.5)..node[basiclabel,left,very near start]{$#3$}(R)
			(R)..controls+(.03,.5)and+(-.03,-.5)..node[basiclabel,left,very near end]{$#1$}(A)
			(R)..controls+(.03,.5)and+(-.03,-.5)..node[basiclabel,right,very near end]{$#2$}(B)
			(CC)..controls+(.03,.5)and+(-.03,-.5)..node[basiclabel,right,very near start]{$#6$}(RR)
			(RR)..controls+(.03,.5)and+(-.03,-.5)..node[basiclabel,left,very near end]{$#4$}(AA)
			(RR)..controls+(.03,.5)and+(-.03,-.5)..node[basiclabel,right,very near end]{$#5$}(BB);
	}}
	
\def\trelb(#1,#2,#3)(#4,#5,#6)(#7){
	\tikz[heighttwo]{
		\draw(-2,2)node[coordinate](A){}(-.7,2)node[coordinate](AA){}
				 (.7,2)node[coordinate](B){}(2,2)node[coordinate](BB){}
				 (-.7,0)node[coordinate](C){}(.7,0)node[coordinate](CC){}
				 (-.7,.8)node[coordinate](R){}(.7,.8)node[coordinate](RR){};
		\draw(-.5,1.2)node[coordinate](I){}(.5,1.6)node[coordinate](J){};
		\draw[trivalent]
			(C)..controls+(.03,.5)and+(-.03,-.5)..node[basiclabel,left,very near start]{$#3$}(R)
			(R)..controls+(.03,.5)and+(-.03,-.5)..node[basiclabel,left,very near end]{$#1$}(A)
			(CC)..controls+(.03,.5)and+(-.03,-.5)..node[basiclabel,right,very near start]{$#6$}(RR)
			(RR)..controls+(.03,.5)and+(-.03,-.5)..node[basiclabel,right,very near end]{$#5$}(BB)
			(R)..controls+(.03,.3)and+(-.03,-.3)..node[basiclabel,right,very near end]{$#2$}(I)
			(RR)..controls+(.03,.3)and+(-.03,-.3)..node[basiclabel,right,very near end]{$#4$}(J)
			(I)..controls+(.03,.3)and+(-.03,-.3)..node[basiclabel,above]{$#7$}(J)
			(I)..controls+(.03,.3)and+(-.03,-.3)..node[basiclabel,left,very near end]{$#4$}(AA)
			(J)..controls+(.03,.3)and+(-.03,-.3)..node[basiclabel,right,very near end]{$#2$}(B);
	}}
	
\def\trelc(#1,#2,#3)(#4,#5,#6)(#7,#8){
	\tikz[heighttwo]{
		\draw(-2,2)node[coordinate](A){}(-.7,2)node[coordinate](AA){}
				 (.7,2)node[coordinate](B){}(2,2)node[coordinate](BB){}
				 (-.7,0)node[coordinate](C){}(.7,0)node[coordinate](CC){}
				 (-.7,.8)node[coordinate](R){}(.7,.8)node[coordinate](RR){};
		\draw(-1.4,1.4)node[coordinate](D){}(1.4,1.4)node[coordinate](DD){};
		\draw(-.5,1.2)node[coordinate](I){}(.5,1.6)node[coordinate](J){};
		\draw[trivalent]
			(CC)..controls+(.03,.5)and+(-.03,-.5)..node[basiclabel,right,very near start]{$#6$}(RR)
			(RR)..controls+(.03,.5)and+(-.03,-.5)..node[basiclabel,right,very near end]{$#5$}(BB)
			(C)..controls+(.03,.5)and+(-.03,-.5)..node[basiclabel,left,very near start]{$#3$}(R)
			(R)..controls+(.03,.3)and+(-.03,-.3)..node[basiclabel,below,near end]{$#8$}(D)
			(D)..controls+(.03,.3)and+(-.03,-.3)..node[basiclabel,left,very near end]{$#1$}(A)
			(D)..controls+(.03,.3)and+(-.03,-.3)..node[basiclabel,left,very near end]{$#4$}(AA)
			(R)..controls+(.03,.3)and+(-.4,.6)..node[basiclabel,left]{$#7$}(J)
			(RR)..controls+(.03,.3)and+(-.03,-.3)..node[basiclabel,left]{$#4$}(J)
			(J)..controls+(.03,.3)and+(-.03,-.3)..node[basiclabel,right,very near end]{$#2$}(B);	
	}}
	
\def\treld(#1,#2,#3)(#4,#5,#6)(#7,#8,#9){
	\tikz[heighttwo]{
		\draw(-2,2)node[coordinate](A){}(-.7,2)node[coordinate](AA){}
				 (.7,2)node[coordinate](B){}(2,2)node[coordinate](BB){}
				 (-.7,0)node[coordinate](C){}(.7,0)node[coordinate](CC){}
				 (-.7,.8)node[coordinate](R){}(.7,.8)node[coordinate](RR){};
		\draw(-1.4,1.4)node[coordinate](D){}(1.4,1.4)node[coordinate](DD){};
		\draw(-.5,1.2)node[coordinate](I){}(.5,1.6)node[coordinate](J){};
		\draw[trivalent]
			(C)..controls+(.03,.5)and+(-.03,-.5)..node[basiclabel,left,very near start]{$#3$}(R)
			(R)..controls+(.03,.3)and+(-.03,-.3)..node[basiclabel,below,near end]{$#8$}(D)
			(D)..controls+(.03,.3)and+(-.03,-.3)..node[basiclabel,left,very near end]{$#1$}(A)
			(D)..controls+(.03,.3)and+(-.03,-.3)..node[basiclabel,left,very near end]{$#4$}(AA)
			(CC)..controls+(.03,.5)and+(-.03,-.5)..node[basiclabel,right,very near start]{$#6$}(RR)
			(RR)..controls+(.03,.3)and+(-.03,-.3)..node[basiclabel,below,near end]{$#9$}(DD)
			(DD)..controls+(.03,.3)and+(-.03,-.3)..node[basiclabel,right,very near end]{$#2$}(B)
			(DD)..controls+(.03,.3)and+(-.03,-.3)..node[basiclabel,right,very near end]{$#5$}(BB)
			(R)..controls+(.03,.5)and+(.03,.5)..node[basiclabel,above]{$#7$}(RR);
	}}
	
\def\trele(#1,#2,#3)(#4,#5,#6)(#7,#8,#9){
	\tikz[heighttwo]{
		\draw(-2,2)node[coordinate](A){}(-.7,2)node[coordinate](AA){}
				 (.7,2)node[coordinate](B){}(2,2)node[coordinate](BB){}
				 (-.7,0)node[coordinate](C){}(.7,0)node[coordinate](CC){};
		\draw(-1.4,1.4)node[coordinate](D){}(1.4,1.4)node[coordinate](DD){};
		\draw(-.4,.4)node[coordinate](K){}(.4,.8)node[coordinate](L){};
		\draw[trivalent]
			(D)..controls+(.03,.3)and+(-.03,-.3)..node[basiclabel,left,very near end]{$#1$}(A)
			(D)..controls+(.03,.3)and+(-.03,-.3)..node[basiclabel,left,very near end]{$#4$}(AA)
			(DD)..controls+(.03,.3)and+(-.03,-.3)..node[basiclabel,right,very near end]{$#2$}(B)
			(DD)..controls+(.03,.3)and+(-.03,-.3)..node[basiclabel,right,very near end]{$#5$}(BB)
			(C)..controls+(.03,.2)and+(-.03,-.2)..node[basiclabel,left,very near start]{$#3$}(K)
			(K)..controls+(.03,.2)and+(-.03,-.2)..node[basiclabel,left]{$#8$}(D)
			(CC)..controls+(.03,.2)and+(-.03,-.2)..node[basiclabel,right,very near start]{$#6$}(L)
			(K)..controls+(.03,.2)and+(-.03,-.2)..node[basiclabel,above]{$#7$}(L)
			(L)..controls+(.03,.2)and+(-.03,-.2)..node[basiclabel,below,near end]{$#9$}(DD);
	}}
	
\def\trelf(#1,#2,#3)(#4,#5,#6)(#7,#8,#9){
	\tikz[heighttwo]{
		\draw(-2,2)node[coordinate](A){}(-.7,2)node[coordinate](AA){}
				 (.7,2)node[coordinate](B){}(2,2)node[coordinate](BB){}
				 (-.7,0)node[coordinate](C){}(.7,0)node[coordinate](CC){};
		\draw(-1.4,1.4)node[coordinate](D){}(1.4,1.4)node[coordinate](DD){};
		\draw(0,.5)node[coordinate](M){}(0,.8)node[coordinate](MM){};
		\draw[trivalent]
			(D)..controls+(.03,.3)and+(-.03,-.3)..node[basiclabel,left,very near end]{$#1$}(A)
			(D)..controls+(.03,.3)and+(-.03,-.3)..node[basiclabel,left,very near end]{$#4$}(AA)
			(DD)..controls+(.03,.3)and+(-.03,-.3)..node[basiclabel,right,very near end]{$#2$}(B)
			(DD)..controls+(.03,.3)and+(-.03,-.3)..node[basiclabel,right,very near end]{$#5$}(BB)
			(C)..controls+(.03,.3)and+(-.03,-.3)..node[basiclabel,left,very near start]{$#3$}(M)
			(CC)..controls+(.03,.3)and+(-.03,-.3)..node[basiclabel,right,very near start]{$#6$}(M)
			(M)..controls+(.03,.3)and+(-.03,-.3)..node[basiclabel,right]{$#9$}(MM)
			(MM)..controls+(.03,.3)and+(-.03,-.3)..node[basiclabel,above]{$#7$}(D)
			(MM)..controls+(.03,.3)and+(-.03,-.3)..node[basiclabel,above]{$#8$}(DD);
	}}
	
\def\trelba(#1,#2,#3)(#4,#5,#6){
	\tikz[heightthree,shift={(0,1.5)}]{
		\draw[trivalent]
			(-2,-1.5)..controls+(.03,.5)and+(-.03,-.5)..node[basiclabel,left,very near start]{$#1$}(-.7,-.3)
			(-.7,-1.5)..controls+(.03,.5)and+(-.03,-.5)..node[basiclabel,below,at start]{$#4$}(.7,-.3)
			(.7,-1.5)..controls+(.03,.5)and+(-.03,-.5)..node[basiclabel,below,at start]{$#2$}(-.7,-.3)
			(2,-1.5)..controls+(.03,.5)and+(-.03,-.5)..node[basiclabel,right,very near start]{$#5$}(.7,-.3)
			(-.7,-.3)..controls+(.03,.5)and+(-.03,-.5)..node[basiclabel,left]{$#3$}(-.7,.3)
			(.7,-.3)..controls+(.03,.5)and+(-.03,-.5)..node[basiclabel,right]{$#6$}(.7,.3)
			(-.7,.3)..controls+(.03,.5)and+(-.03,-.5)..node[basiclabel,left,very near end]{$#1$}(-2,1.5)
			(.7,.3)..controls+(.03,.5)and+(-.03,-.5)..node[basiclabel,above,at end]{$#4$}(-.7,1.5)
			(-.7,.3)..controls+(.03,.5)and+(-.03,-.5)..node[basiclabel,above,at end]{$#2$}(.7,1.5)
			(.7,.3)..controls+(.03,.5)and+(-.03,-.5)..node[basiclabel,right,very near end]{$#5$}(2,1.5);
	}}
\def\trelbb(#1,#2)(#3,#4)(#5,#6,#7,#8,#9){
	\tikz[heightthree,shift={(0,1.5)}]{
		\draw[trivalent]
			(-2,-1.5)..controls+(.03,.5)and+(-.03,-.5)..node[basiclabel,left,very near start]{$#1$}(-1.3,-.9)
			(-.7,-1.5)..controls+(.03,.5)and+(-.03,-.5)..node[basiclabel,below,at start]{$#3$}(-1.3,-.9)
			(.7,-1.5)..controls+(.03,.5)and+(-.03,-.5)..node[basiclabel,below,at start]{$#2$}(1.3,-.9)
			(2,-1.5)..controls+(.03,.5)and+(-.03,-.5)..node[basiclabel,right,very near start]{$#4$}(1.3,-.9)
			(-1.3,-.9)..controls+(.03,.5)and+(-.03,-.5)..node[basiclabel,left,near start]{$#6$}(0,-.2)
			(1.3,-.9)..controls+(.03,.5)and+(-.03,-.5)..node[basiclabel,right,near start]{$#8$}(0,-.2)
			(0,-.2)..controls+(.03,.5)and+(-.03,-.5)..node[basiclabel,right]{$#9$}(0,.2)
			(0,.2)..controls+(.03,.5)and+(-.03,-.5)..node[basiclabel,left,near end]{$#5$}(-1.3,.9)
			(0,.2)..controls+(.03,.5)and+(-.03,-.5)..node[basiclabel,right,near end]{$#7$}(1.3,.9)
			(-1.3,.9)..controls+(.03,.5)and+(-.03,-.5)..node[basiclabel,left,very near end]{$#1$}(-2,1.5)
			(-1.3,.9)..controls+(.03,.5)and+(-.03,-.5)..node[basiclabel,above,at end]{$#3$}(-.7,1.5)
			(1.3,.9)..controls+(.03,.5)and+(-.03,-.5)..node[basiclabel,above,at end]{$#2$}(.7,1.5)
			(1.3,.9)..controls+(.03,.5)and+(-.03,-.5)..node[basiclabel,right,very near end]{$#4$}(2,1.5);
	}}

\def\cfone(#1)(#2){
	\tikz[heighttwo]{
		\draw[trivalent](0,1)node[matrix](A){$#1$}..controls+(.1,1)and+(.1,1)..node[near end,right,basiclabel]{$#2$}(2,1)..controls+(-.1,-1)and+(-.1,-1)..(A);
	}}

\def\cftwo(#1,#2)(#3,#4,#5){
	\tikz[heighttwo]{
		\draw(0,1)node[smatrix](A){$#1$}(1.2,1)node[smatrix](B){$#2$};
		\draw(.6,2)node[coordinate](RR){}(.6,0)node[coordinate](R){}(2.4,1)node[coordinate](Z){};
		\foreach\xa/\xb/\xc/\xd in{R/A/.5/,R/B/.5/,A/RR/.5/#3,B/RR/.5/#4}{
			\draw[trivalent](\xa)..controls+(.05,\xc)and+(.05,-\xc)..node[sloped,basiclabel,above]{$\xd$}(\xb);}
		\draw[trivalent](RR)..controls+(.1,.3)and+(0,1.5)..node[near end,right,basiclabel]{$#5$}(Z)..controls+(0,-1.5)and+(0,-.3)..(R);
	}}

\def\cftwob(#1,#2)(#3,#4,#5){
	\tikz[heighttwo]{
		\draw(0,1.5)node[smatrix](A){$#1$}(1.2,1.5)node[smatrix](B){$#2$};		
		\draw[trivalent]
			 (.5,.5)..controls+(.02,.2)and+(-.02,-.2)..node[basiclabel,right]{$#5$}(.5,.8)..controls+(.02,.2)and+(-.02,-.5)..(A)
			(.5,.8)..controls+(.02,.2)and+(-.02,-.5)..(B);
		\draw[trivalent]
			 (A)..controls+(.05,.5)and+(0,1)..(-1.2,1)node[basiclabel,left]{$#3$}..controls+(0,-1)and+(-.14,-.7)..(.5,.5)
			 (B)..controls+(.05,.5)and+(0,1)..(2.4,1)node[basiclabel,right]{$#4$}..controls+(0,-1)and+(.14,-.7)..(.5,.5);
	}}

\def\cftwoc(#1,#2)(#3,#4,#5){
	\tikz[heightthree,scale=1.2]{
		\draw[trivalent](-.3,2)node[smatrix](A){$#1$}..controls+(.12,1.2)and+(.12,1.2)..node[near end,right,basiclabel]{$#3$}(2.3,2)
			..controls+(-.05,-.5)and+(.05,.5)..(2.3,1)..controls+(-.12,-1.2)and+(-.07,-.7)..(-.3,.5);
		 \draw[trivalent](.3,2)node[smatrix](B){$#2$}..controls+(.08,.8)and+(.08,.8)..(1.7,2)node[left,basiclabel]{$#4$}
			..controls+(-.05,-.5)and+(.05,.5)..(1.7,1)..controls+(-.08,-.8)and+(-.03,-.3)..(.3,.5);
		\draw[trivalent](-.3,.5)..controls+(.02,.2)and+(-.02,-.2)..(0,1)
			 ..controls+(.02,.2)and+(-.02,-.2)..node[basiclabel,left]{$#5$}(0,1.2)..controls+(.04,.4)and+(-.04,-.4)..(A)
			(.3,.5)..controls+(.02,.2)and+(-.02,-.2)..(0,1)(0,1.2)..controls+(.04,.4)and+(-.04,-.4)..(B);
	}}

\def\cfthree(#1,#2,#3)(#4,#5,#6,#7,#8,#9){
	\tikz[]{
	\draw(0,0)node[matrix](A){$#1$}(1,0)node[matrix](B){$#2$}(2,0)node[matrix](C){$#3$};
	\draw(.5,.7)node[coordinate](RR){}(1.2,1.2)node[coordinate](SS){}
			 (.5,-.7)node[coordinate](R){}(1.2,-1.2)node[coordinate](S){}
			 (3,0)node[coordinate](Z){};
	\foreach\xa/\xb/\xc/\xd in{R/A/.3/,R/B/.3/,S/R/.3/#8,S/C/.5/,A/RR/.3/#4,B/RR/.3/#5,RR/SS/.3/#7,C/SS/.5/#6}{
		\draw[trivalent](\xa)..controls+(0,\xc)and+(0,-\xc)..node[sloped,basiclabel,above]{\xd}(\xb);}
	\draw[trivalent](SS)..controls+(0,.7)and+(0,1.5)..node[near end,right,basiclabel]{$#9$}(Z)..controls+(0,-1.5)and+(0,-.7)..(S);
}}

\def\cfrelaa(#1,#2)(#3,#4){
	\tikz[heighttwo]{
		\draw[trivalent](-.3,1)node[smatrix](A){$#1$}..controls+(.12,1.2)and+(.12,1.2)..node[near end,right,basiclabel]{$#3$}(2.3,1)..controls+(-.12,-1.2)and+(-.12,-1.2)..(A);
		\draw[trivalent](.3,1)node[smatrix](B){$#2$}..controls+(.08,.8)and+(.08,.8)
			..(1.7,1)node[left,basiclabel]{$#4$}..controls+(-.08,-.8)and+(-.08,-.8)..(B);
	}}

\def\cfrelab(#1)(#2,#3,#4){
	\tikz[heightthree,scale=1.2]{
		\draw[trivalent](-.3,2.5)..controls+(.07,.7)and+(.12,1.2)..node[near end,right,basiclabel]{$#2$}(2.3,2)
			..controls+(-.05,-.5)and+(.05,.5)..(2.3,1)..controls+(-.12,-1.2)and+(-.07,-.7)..(-.3,.5);
		\draw[trivalent](.3,2.5)..controls+(.03,.3)and+(.08,.8)..(1.7,2)node[left,basiclabel]{$#3$}
			..controls+(-.05,-.5)and+(.05,.5)..(1.7,1)..controls+(-.08,-.8)and+(-.03,-.3)..(.3,.5);
		\draw[trivalent](-.3,.5)..controls+(.02,.2)and+(-.02,-.2)..(0,1)..controls+(.02,.2)and+(-.02,-.2)
			..node[smatrix](A){$#1$}(0,2)node[basiclabel,left]{$#4$}..controls+(.02,.2)and+(-.02,-.2)..(-.3,2.5)
			(.3,.5)..controls+(.02,.2)and+(-.02,-.2)..(0,1)(0,2)..controls+(.02,.2)and+(-.02,-.2)..(.3,2.5);
	}}
	
\def\cfrelac(#1)(#2,#3,#4){
	\tikz[heightthree,scale=1.2]{
		\draw[trivalent](2,.7)..controls+(-.05,-.5)and+(-.1,-1)..(0,1)..controls+(.05,.5)and+(-.05,-.5)
			..node[midmx]{$#1$}(0,2)..controls+(.1,1)and+(.05,.5)..node[left,basiclabel,near start]{$#4$}(2,2.3);
		\draw[trivalent]
			(2,.7)..controls+(.02,.2)and+(-.02,-.2)..(2,.8)..controls+(.02,.2)and+(-.02,-.2)
				 ..(1.5,1.5)node[basiclabel,left]{$#2$}..controls+(.03,.3)and+(-.03,-.3)..(2,2.2)..controls+(.03,.3)and+(-.03,-.3)..(2,2.3)
			 (2,.8)..controls+(.03,.3)and+(-.03,-.3)..(2.5,1.5)node[basiclabel,right]{$#3$}..controls+(.03,.3)and+(-.03,-.3)..(2,2.2);
	}}

\def\tupp<#1>#2#3{\foreach\x in{1,...,#2}%
	\draw[shift={#1}]([shift={(-1,0)}]\x,0)..controls+(.1,.5)and+(-.1,-.5)..([shift={(-1,0)}]\x,#3);}

\def\tmupp<#1>#2#3{\foreach\xa/\xb in{#3}{\draw[shift={#1}](\xa,0)..controls+(.1,.5)and+(-.1,-.5)..(\xa,#2)node[smatrix]{\xb}..controls+(.1,.5)and+(-.1,-.5)..([shift={(0,#2)}]\xa,#2);}}

\def\tbbasis<#1>#2{\foreach\xa/\xb in{#2}{\draw[shift={#1}](\xa,.5)node[rectangle]{$e_\xb$}..controls+(.1,.5)and+(-.1,-.5)..(\xa,1);}}
\def\ttbasis<#1>#2{\foreach\xa/\xb in{#2}{\draw[shift={#1}](\xa,0)..controls+(.1,.5)and+(-.1,-.5)..(\xa,.5)node[rectangle]{$e_\xb$};}}

\def\tperm<#1>#2{\foreach\xa/\xb in{#2}{\draw[shift={#1}](\xa,0)..controls+(.1,.5)and+(-.1,-.5)..(\xb,1);}}

\def\ttl<#1>#2#3#4{%
	\foreach\xa/\xb in{#2}{\draw[shift={#1}](\xa,0)..controls+(.1,.5)and+(.1,.5)..(\xb,0);}%
	\foreach\xa/\xb in{#3}{\draw[shift={#1}](\xa,1)..controls+(-.1,-.5)and+(-.1,-.5)..(\xb,1);}%
	\foreach\xa/\xb in{#4}{\draw[shift={#1}](\xa,0)..controls+(.1,.5)and+(-.1,-.5)..(\xb,1);}%
}
	
\def\tcupp<#1>#2{\foreach\x in{1,...,#2}{\draw[shift={#1}](-\x,0)..controls+(-.05,-\x)and+(-.05,-\x)..(\x,0);}}
\def\tcapn<#1>#2{\foreach\x in{1,...,#2}{\draw[shift={#1}](-\x,0)..controls+(.05,\x)and+(.05,\x)..(\x,0);}}

\def\tcloser<#1>#2#3{\tcupp<#1>{#2}\tcapn<([shift={#1}]0,#3)>{#2}\tupp<([shift={#1}]1,0)>{#2}{#3}}

\def\tsupp<#1>#2#3{\tupp<#1>{#2}{#3}\filldraw[black](-.5,.3)rectangle([shift={(-.5,-.3)}]#2,#3);%
			\draw([shift={(-.5,0)},scale=.5]#2,#3)node[symlabel]{$#2$};}
\def\tasupp<#1>#2#3{\tupp<#1>{#2}{#3}\filldraw[fill=white](-.5,.3)rectangle([shift={(-.5,-.3)}]#2,#3);
			\draw([shift={(-.5,0)},scale=.5]#2,#3)node[asymlabel]{$#2$};}

\def\dftn#1{\ifcase{#1}.7\or.78\or.85\or.9\or.95\or1\or1.04\or1.08\or1.11\fi}

\newcount\tempn%
\def\remainder#1#2{\tempn=#1\divide\tempn by#2\multiply\tempn by#2\advance#1 by-\tempn}%

\newcount\nnnn\newcount\perma\newcount\positn\newcount\mmmm%
\def\tpermb#1{
  \nnnn=#1
  \perma=\nnnn
  \positn=100000000
  \mmmm=0
  \loop%
    \perma=\nnnn\divide\perma by\positn\remainder{\perma}{10}%
    \ifnum\perma<9 \number\mmmm/\number\perma\fi
    \advance\mmmm by1%
    \ifnum\positn>1\divide\positn by10,%
  \repeat%
}%

\def\tbigexc{
	\tikz[heightthree,xscale=1.3,yscale=1.1]{
		\draw(0,0)..controls+(.03,.3)and+(-.03,-.3)..(0,1)..controls+(.03,.3)and+(-.03,-.3)..(1,2)
			 ..controls+(.03,.5)and+(.03,.5)..node[dotedge]{}(2,2)..controls+(-.03,-.4)and+(-.03,-.4)..node[dotedge]{}(3,2)
			..controls+(.03,.3)and+(-.03,.3)..node[dotedge](M){}(3,3);
		\draw(1,0)..controls+(.03,.3)and+(-.03,-.3)..(1,1)..controls+(.03,.3)and+(-.03,-.3)..(0,2)
			..controls+(.03,.3)and+(-.03,-.3)..(0,3);
		\draw(2,0)..controls+(.03,.3)and+(-.03,-.3)..(2,1)..controls+(.03,.4)and+(.03,.4)..node[dotedge]{}(3,1)
			..controls+(-.03,-.3)and+(.03,.3)..(4,0);
		\draw(3,0)..controls+(.03,.3)and+(-.03,-.3)..(4,1)..controls+(.03,.3)and+(-.03,-.3)..node[dotedge](N){}(4,3);
		\draw(M)--+(.4,0);
		\draw(N)--+(-.4,0);
	}}
\def\tbigexccut{
	\tikz[heightthree,xscale=1.3,yscale=1.1]{
		\draw(0,0)..controls+(.03,.3)and+(-.03,-.3)..(0,1)..controls+(.03,.3)and+(-.03,-.3)..(1,2)
			 ..controls+(.03,.5)and+(.03,.5)..node[dotedge]{}(2,2)..controls+(-.03,-.4)and+(-.03,-.4)..node[dotedge]{}(3,2)
			..controls+(.03,.3)and+(-.03,.3)..node[dotedge](M){}(3,3);
		\draw(1,0)..controls+(.03,.3)and+(-.03,-.3)..(1,1)..controls+(.03,.3)and+(-.03,-.3)..(0,2)
			..controls+(.03,.3)and+(-.03,-.3)..(0,3);
		\draw(2,0)..controls+(.03,.3)and+(-.03,-.3)..(2,1)..controls+(.03,.4)and+(.03,.4)..node[dotedge]{}(3,1)
			..controls+(-.03,-.3)and+(.03,.3)..(4,0);
		\draw(3,0)..controls+(.03,.3)and+(-.03,-.3)..(4,1)..controls+(.03,.3)and+(-.03,-.3)..(4,2)
			..controls+(.03,.3)and+(-.03,-.3)..node[dotedge](N){}(4,3);
		\draw(M)--+(.4,0);
		\draw(N)--+(-.4,0);
		\draw[red,dotted](-.5,0)--(4.5,0)(-.5,1)--(4.5,1)(-.5,2)--(4.5,2)(-.5,3)--(4.5,3);
	}}
	
\def\tcfsyma(#1,#2,#3)(#4,#5,#6){
	\tikz[scale=.5,xscale=.8]{
		\foreach\xa in{0,8,16}{
			\draw[shift={(\xa,0)}]
				(-3.5,1)..controls+(.05,.5)and+(.05,.5)..(-5.5,1)..controls+(-.05,-.5)and+(.05,.5)..(-5.5,-1)
					..controls+(-.05,-.5)and+(-.05,-.5)..(-3.5,-1)..controls+(.05,.5)and+(-.05,-.5)..(-3.5,1)
				(-1,1)..controls+(.2,1.5)and+(.2,1.5)..(-8,1)..controls+(-.05,-.5)and+(.05,.5)..(-8,-1)
					..controls+(-.2,-1.5)and+(-.2,-1.5)..(-1,-1)..controls+(.05,.5)and+(-.05,-.5)..(-1,1);
			\foreach\xb/\xc in{-7/-1.4,-6/-1.3,-3/-1.3,-2/-1.4}{
				\draw[shift={(\xa,0)}]node[basiclabel,scale=1.5]at(\xb,\xc){.};}
		}
		 \filldraw[black](-9.5,.3)rectangle+(4.2,.4)(-3.7,.3)rectangle+(6.4,.4)(4.3,.3)rectangle+(6.4,.4)(12.3,.3)rectangle+(4.2,.4);
		\draw
			(-8,-.5)node[smatrix,scale=.9]{$#1$}(-5.5,-.5)node[smatrix,scale=.9]{$#1$}
			(0,-.5)node[smatrix,scale=.9]{$#2$}(2.5,-.5)node[smatrix,scale=.9]{$#2$}
			(8,-.5)node[smatrix,scale=.9]{$#3$}(10.5,-.5)node[smatrix,scale=.9]{$#3$};
		\filldraw[white]
				(17,1)--(16,.5)--(17,0)--cycle
				(-9.5,.5)--(-8.5,1)--(-10,1)--(-10,0)--(-8.5,0)--cycle;
		\draw[snake=brace,raise snake=-3pt]
			(-8,2.6)--(-1,2.6)node[basiclabel,pos=.5,above]{$#4$}
			(0,2.6)--(7,2.6)node[basiclabel,pos=.5,above]{$#5$}
			(8,2.6)--(15,2.6)node[basiclabel,pos=.5,above]{$#6$};
}}
\def\tcfsymb(#1,#2)(#3,#4,#5){
	\tikz[scale=.5,xscale=.8]{
		\foreach\xa in{0,8,16}{
			\draw[shift={(\xa,0)}]
				(-3.5,1)..controls+(.05,.5)and+(.05,.5)..(-5.5,1)..controls+(-.05,-.5)and+(.05,.5)..(-5.5,-1)
					..controls+(-.05,-.5)and+(-.05,-.5)..(-3.5,-1)..controls+(.05,.5)and+(-.05,-.5)..(-3.5,1)
				(-1,1)..controls+(.2,1.5)and+(.2,1.5)..(-8,1)..controls+(-.05,-.5)and+(.05,.5)..(-8,-1)
					..controls+(-.2,-1.5)and+(-.2,-1.5)..(-1,-1)..controls+(.05,.5)and+(-.05,-.5)..(-1,1);
			\foreach\xb/\xc in{-7/-1.4,-6/-1.3,-3/-1.3,-2/-1.4}{
				\draw[shift={(\xa,0)}]node[basiclabel,scale=1.5]at(\xb,\xc){.};}
		}
		 \filldraw[black](-9.5,.3)rectangle+(4.2,.4)(-3.7,.3)rectangle+(6.4,.4)(4.3,.3)rectangle+(6.4,.4)(12.3,.3)rectangle+(4.2,.4);
		\draw
			(10.5,-.5)node[smatrix,scale=.9]{$#2$}(8,-.5)node[smatrix,scale=.9]{$#2$}
			(7,-.5)node[smatrix,scale=.9]{$#2$}(4.5,-.5)node[smatrix,scale=.9]{$#2$}
			(2.5,-.5)node[smatrix,scale=.9]{$#1$}(0,-.5)node[smatrix,scale=.9]{$#1$}
			(-1,-.5)node[smatrix,scale=.9]{$#1$}(-3.5,-.5)node[smatrix,scale=.9]{$#1$};
		\filldraw[white]
				(17,1)--(16,.5)--(17,0)--cycle
				(-9.5,.5)--(-8.5,1)--(-10,1)--(-10,0)--(-8.5,0)--cycle;
		\draw[snake=brace,raise snake=-3pt]
			(-8,2.6)--(-1,2.6)node[basiclabel,pos=.5,above]{$#3$}
			(0,2.6)--(7,2.6)node[basiclabel,pos=.5,above]{$#4$}
			(8,2.6)--(15,2.6)node[basiclabel,pos=.5,above]{$#5$};
}}

\def\srinva#1#2{
	\tikz[heighttwo]{
		\foreach\xa in{0,.5,2}{\draw(\xa,0)..controls+(.1,.5)and+(-.1,-.5)..(\xa,1);}
		\filldraw[black](-.3,.3)rectangle+(2.6,.3)node[symlabelright]{$#1$};
		\foreach\xa in{0,.5,2}{\draw(\xa,1)..controls+(.1,.5)and+(-.1,-.5)..node[midsmx]{$#2$}(\xa,2);}
		\foreach\xa/\xb in{1/1,1.5/1,1/1,1.5/1}{\node[basiclabel,scale=2]at(\xa,\xb){.};}
	}}
\def\srinvb#1#2{
	\tikz[heighttwo]{
		\foreach\xa in{0,.5,2}{\draw(\xa,1)..controls+(.1,.5)and+(-.1,-.5)..(\xa,2);}
		\filldraw[black](-.3,1.4)rectangle+(2.6,.3)node[symlabelright]{$#1$};
		\foreach\xa in{0,.5,2}{\draw(\xa,0)..controls+(.1,.5)and+(-.1,-.5)..node[midsmx]{$#2$}(\xa,1);}
		\foreach\xa/\xb in{1/1,1.5/1}{\node[basiclabel,scale=2]at(\xa,\xb){.};}
	}}
	
\def\srstack#1#2{
	\tikz[heighttwo]{
		\foreach\xa in{-1.5,0,.5,2}{\draw(\xa,0)..controls+(.1,.5)and+(-.1,-.5)..(\xa,2);}
		\filldraw[black]
			(-1.8,.3)rectangle+(4.1,.3)node[symlabelright]{$#1$}
			(.3,1.4)rectangle+(1.9,.3)node[symlabelright]{$#2$};
		\foreach\xa/\xb in{1/1,1.5/1,-1/1,-.5/1}{\node[basiclabel,scale=2]at(\xa,\xb){.};}}
	=
	\tikz[heighttwo]{
		\foreach\xa in{-1.5,0,.5,2}{\draw(\xa,0)..controls+(.1,.5)and+(-.1,-.5)..(\xa,2);}
		\filldraw[black](-1.8,.3)rectangle+(4.1,.3)node[symlabelright]{$#1$};
		\foreach\xa/\xb in{1/1,1.5/1,-1/1,-.5/1}{\node[basiclabel,scale=2]at(\xa,\xb){.};}}
}

\def\srcapping#1{
	\tikz[heightoneonehalf]{
		\draw(-1,0)..controls+(.05,.5)and+(-.05,-.5)..(-1,.8)..controls+(.05,.5)and+(.05,.5)..(0,.8)
			..controls+(-.05,-.5)and+(.05,.5)..(0,0);
		\foreach\xa in{.5,2}{\draw(\xa,0)..controls+(.05,.5)and+(-.05,-.5)..(\xa,1.5);}
		\filldraw[black](-1.3,.3)rectangle+(3.6,.3)node[symlabelright]{$#1$};
		\foreach\xa/\xb in{1/1,1.5/1}{\node[basiclabel,scale=2]at(\xa,\xb){.};}}
	=0=
	\tikz[heightoneonehalf]{
		\draw(-1,1.5)..controls+(-.05,-.5)and+(.05,.5)..(-1,.8)..controls+(-.05,-.5)and+(-.05,-.5)..(0,.7)
			..controls+(.05,.5)and+(-.05,-.5)..(0,1.5);
		\foreach\xa in{.5,2}{\draw(\xa,0)..controls+(.05,.5)and+(-.05,-.5)..(\xa,1.5);}
		\filldraw[black](-1.3,.9)rectangle+(3.6,.3)node[symlabelright]{$#1$};
		\foreach\xa/\xb in{1/.5,1.5/.5}{\node[basiclabel,scale=2]at(\xa,\xb){.};}}
}

\def\srextrema#1{
	\tikz[heightoneonehalf,scale=.8]{
		\draw(-.5,1.5)..controls+(-.02,-.2)and+(.02,.2)..(-.5,.7)..controls+(-.05,-.5)and+(-.05,-.5)
			..node[dotedge](M){}(.5,.7)..controls+(.02,.2)and+(-.02,0)..(.5,1.5);
		\draw(-2,1.5)..controls+(0,-.1)and+(0,.1)..(-2,.7)..controls+(-.15,-1.5)and+(-.15,-1.5)
			..node[dotedge](N){}(2,.7)..controls+(0,.1)and+(0,-.1)..(2,1.5);
		\draw(M)--+(0,.3)(N)--+(0,.3);
		\foreach\xa/\xb in{-1.2/0,-.7/.2,1.1/0,.6/.2}{\node[basiclabel,scale=2]at(\xa,\xb){.};}
		\filldraw[black](-.3,.8)rectangle+(-1.9,.3)node[symlabelleft]{$#1$};
	}
	=\tikz[heightoneonehalf,scale=.8,xscale=-1]{
		\draw(-.5,1.5)..controls+(-.02,-.2)and+(.02,.2)..(-.5,.7)..controls+(-.05,-.5)and+(-.05,-.5)
			..node[dotedge](M){}(.5,.7)..controls+(.02,.2)and+(-.02,0)..(.5,1.5);
		\draw(-2,1.5)..controls+(0,-.1)and+(0,.1)..(-2,.7)..controls+(-.15,-1.5)and+(-.15,-1.5)
			..node[dotedge](N){}(2,.7)..controls+(0,.1)and+(0,-.1)..(2,1.5);
		\draw(M)--+(0,.3)(N)--+(0,.3);
		\foreach\xa/\xb in{-1.2/0,-.7/.2,1.1/0,.6/.2}{\node[basiclabel,scale=2]at(\xa,\xb){.};}
		\filldraw[black](-.3,.8)rectangle+(-1.9,.3)node[symlabelright]{$#1$};
	}
}

\def\srrecaa#1{\tikz[heighttwos]{
	\tperm<(0,0)>{0/0,1/1,2/2,5/5}
	\tperm<(0,1)>{0/0,1/1,2/2,5/5}
	\filldraw[black](-.5,1.3)rectangle+(6,.4)node[symlabelright,scale=.8]{$#1$};
	\foreach\xa/\xb in{3/.5,4/.5}{\node[basiclabel,scale=1.5]at(\xa,\xb){.};}
}}
\def\srrecab#1{\tikz[heighttwos]{
	\tperm<(0,0)>{0/0,1/1,2/2,5/5}
	\tperm<(0,1)>{0/0,1/1,2/2,5/5}
	\filldraw[black](.6,1.3)rectangle+(5,.4)node[symlabelright,scale=.8]{$#1-1$};
	\foreach\xa/\xb in{3/.5,4/.5}{\node[basiclabel,scale=1.5]at(\xa,\xb){.};}
}}
\def\srrecac#1{\tikz[heighttwos]{
	\ttl<(0,0)>{0/1}{0/1}{2/2,5/5}
	\tperm<(0,1)>{0/0,1/1,2/2,5/5}
	\filldraw[black](.6,1.3)rectangle+(5,.4)node[symlabelright,scale=.8]{$#1-1$};
	\foreach\xa/\xb in{3/.5,4/.5}{\node[basiclabel,scale=1.5]at(\xa,\xb){.};}
}}
\def\srrecad#1{\tikz[heighttwos]{
	\ttl<(0,0)>{1/2}{0/1}{0/2,3/3,6/6}
	\tperm<(0,1)>{0/0,1/1,2/2,3/3,6/6}
	\filldraw[black](.6,1.3)rectangle+(6,.4)node[symlabelright,scale=.8]{$#1-1$};
	\foreach\xa/\xb in{4/.5,5/.5}{\node[basiclabel,scale=1.5]at(\xa,\xb){.};}
}}
\def\srrecae#1{\tikz[heighttwos]{
	\ttl<(0,0)>{4/5}{0/1}{0/2,3/5,6/6,9/9}
	\tperm<(0,1)>{0/0,1/1,2/2,5/5,6/6,9/9}
	\filldraw[black](.6,1.3)rectangle+(9,.4)node[symlabelright,scale=.8]{$#1-1$};
	\foreach\xa/\xb in{2/.5,3/.5,7/.5,8/.5}{\node[basiclabel,scale=1.5]at(\xa,\xb){.};}
}}
\def\srrecaf#1{\tikz[heighttwos]{
	\ttl<(0,0)>{4/5}{0/1}{0/2,3/5}
	\tperm<(0,1)>{0/0,1/1,2/2,5/5}
	\filldraw[black](.6,1.3)rectangle+(5,.4)node[symlabelright,scale=.8]{$#1-1$};
	\foreach\xa/\xb in{2/.5,3/.5}{\node[basiclabel,scale=1.5]at(\xa,\xb){.};}
}}

\def\trarrowa{\tikz[heightones]{
	\tperm<(0,0)>{0/0,1/1,2/2,5/5}
	\foreach\xa/\xb in{3/.5,4/.5}{\node[basiclabel,scale=1.5]at(\xa,\xb){.};}}}
\def\trarrowb{\tikz[heightones]{
	\tperm<(0,0)>{0/1,1/0,2/2,5/5}
	\foreach\xa/\xb in{3/.5,4/.5}{\node[basiclabel,scale=1.5]at(\xa,\xb){.};}}}
\def\trarrowc{\tikz[heightones]{
	\tperm<(0,0)>{0/1,1/2,2/0,5/5}
	\foreach\xa/\xb in{3/.5,4/.5}{\node[basiclabel,scale=1.5]at(\xa,\xb){.};}}}
\def\trarrowd{\tikz[heightones]{
	\tperm<(0,0)>{0/1,3/4,4/0,5/5,8/8}
	\foreach\xa/\xb in{1/.2,2/.2,2/.8,3/.8,6/.5,7/.5}{\node[basiclabel,scale=1.5]at(\xa,\xb){.};}}}
\def\trarrowe{\tikz[heightones]{
	\tperm<(0,0)>{0/1,1/2,4/5,5/0}
	\foreach\xa/\xb in{2/.2,3/.2,3/.8,4/.8}{\node[basiclabel,scale=1.5]at(\xa,\xb){.};}}}

\def\srrecb#1#2{
	\tikz[heightthrees]{
		\tperm<(0,0)>{0/0,1/1,2/2,5/5}
		\tperm<(0,1)>{0/0,1/1,2/2,5/5}
		\tperm<(0,2)>{0/0,1/1,2/2,5/5}
		\filldraw[black](-.5,1.3)rectangle+(6,.4)node[symlabelright,scale=.8]{$#1$};
		\foreach\xa/\xb in{3/.5,4/.5,3/2.5,4/2.5}{\node[basiclabel,scale=1.5]at(\xa,\xb){.};}
	}
	=\tikz[heightthrees]{
		\tperm<(0,0)>{0/0,1/1,4/4,6/6,9/9}
		\tperm<(0,1)>{0/0,1/2,4/5,6/6,9/9}
		\tperm<(0,2)>{0/0,2/2,5/5,6/6,9/9}
		\filldraw[black]
			(4.5,.3)rectangle+(-5,.4)node[symlabelleft,scale=.8]{$#2$}
			(5.5,.3)rectangle+(4,.4)node[symlabelright,scale=.8]{$#1-#2$}
			(1.5,2.3)rectangle+(8,.4)node[symlabelright,scale=.8]{$#1-1$};
		\foreach\xa in{2.5,3.5,7,8}{\node[basiclabel,scale=1.5]at(\xa,1.5){.};}
	}
	+(-1)^{#2}\left(\frac{#1-#2}{#1}\right)
	\tikz[heightthrees]{
		\tperm<(0,0)>{0/0,3/3,4/4,6/6,7/7,10/10}
		\ttl<(0,1)>{4/6}{0/2}{0/3,3/6,7/7,10/10}
		\tperm<(0,2)>{0/0,2/2,3/3,6/6,7/7,10/10}
		\filldraw[black]
			(4.5,.3)rectangle+(-5,.4)node[symlabelleft,scale=.8]{$#2$}
			(5.5,.3)rectangle+(5,.4)node[symlabelright,scale=.8]{$#1-#2$}
			(1.5,2.3)rectangle+(9,.4)node[symlabelright,scale=.8]{$#1-1$};
		\foreach\xa/\xb in{1/1,2/1,4/2,5/2,8/1.5,9/1.5}{\node[basiclabel,scale=1.5]at(\xa,\xb){.};}
	}
}

\def\srrecca#1{
	\tikz[heightthrees]{
		\tperm<(0,0)>{0/0,1/1,2/2,5/5}
		\tperm<(0,1)>{0/0,1/1,2/2,5/5}
		\tperm<(0,2)>{0/0,1/1,2/2,5/5}
		\filldraw[black](-.5,1.3)rectangle+(6,.4)node[symlabelright,scale=.8]{$#1$};
		\foreach\xa/\xb in{3/.5,4/.5,3/2.5,4/2.5}{\node[basiclabel,scale=1.5]at(\xa,\xb){.};}
}}
\def\srreccb#1{
	\tikz[heightthrees]{
		\tperm<(0,0)>{-.2/-.2,1/1,2/2,5/5}
		\tperm<(0,1)>{-.2/-.2,1/1,2/2,5/5}
		\tperm<(0,2)>{-.2/-.2,1/1,2/2,5/5}
		\filldraw[black](.6,1.3)rectangle+(5,.4)node[symlabelright,scale=.8]{$#1-1$};
		\foreach\xa/\xb in{3/.5,4/.5,3/2.5,4/2.5}{\node[basiclabel,scale=1.5]at(\xa,\xb){.};}
}}
\def\srreccc#1{
	\tikz[heightthrees]{
		\tperm<(0,0)>{-1/-1,1/1,2/2,5/5}
		\ttl<(0,1)>{-1/1}{-1/1}{2/2,5/5}
		\tperm<(0,2)>{-1/-1,1/1,2/2,5/5}
		\filldraw[black]
			(.5,2.3)rectangle+(5,.4)node[symlabelright,scale=.8]{$#1-1$}
			(.5,.3)rectangle+(5,.4)node[symlabelright,scale=.8]{$#1-1$};
		\foreach\xa in{3,4}{\node[basiclabel,scale=1.5]at(\xa,1.5){.};}
}}

\def\srloopa#1{
	\tikz[heightthrees]{
		\ttl<(0,0)>{}{-3/-1}{0/0,3/3}
		\tperm<(0,1)>{-3/-3,-1/-1,0/0,3/3}
		\ttl<(0,2)>{-3/-1}{}{0/0,3/3}
		\filldraw[black](-1.5,1.3)rectangle+(5,.4)node[symlabelright,scale=.8]{$#1$};
		\foreach\xa/\xb in{1/.5,2/.5,1/2.5,2/2.5}{\node[basiclabel,scale=1.5]at(\xa,\xb){.};}
}}
\def\srloopb#1{
	\tikz[heightthrees]{
		\ttl<(0,0)>{}{-3/-1}{0/0,3/3}
		\tperm<(0,1)>{-3/-3,-1/-1,0/0,3/3}
		\ttl<(0,2)>{-3/-1}{}{0/0,3/3}
		\filldraw[black](-.5,1.3)rectangle+(4,.4)node[symlabelright,scale=.8]{$#1$};
		\foreach\xa/\xb in{1/.5,2/.5,1/2.5,2/2.5}{\node[basiclabel,scale=1.5]at(\xa,\xb){.};}
}}
\def\srloopc#1#2{
	\tikz[heightthrees]{
		\ttl<(0,0)>{}{-5/-4}{0/0,3/3}
		\tperm<(0,1)>{-5/-5,-4/-4,-1/-1,0/0,3/3}
		\ttl<(0,2)>{-5/-4}{}{0/0,3/3}
		\draw(-1,2)..controls+(.2,1.5)and+(.2,1.5)..(-8,2)..controls+(-.05,-.5)and+(.05,.5)..(-8,1)
			..controls+(-.2,-1.5)and+(-.2,-1.5)..(-1,1);
		\filldraw[black](-4.5,1.3)rectangle+(8,.4)node[symlabelright,scale=.8]{$#1$};
		\foreach\xa/\xb in{-7/1.5,-6/1.5,-3.2/2.3,-2.2/2.4,-3.2/.7,-2.2/.6,1/.5,2/.5,1/2.5,2/2.5}{\node[basiclabel,scale=1.5]at(\xa,\xb){.};}
		\node[basiclabel,left]at(-7.5,1.5){$#2\Bigl\{$};
}}

\def\cfrelba(#1,#2,#3)(#4,#5,#6,#7){
	\tikz[heighttwo,yscale=.7,shift={(0,1)}]{
	\draw
		(0,-1.3)node[coordinate](A){}(1,-1.3)node[coordinate](B){}(2,-1.3)node[coordinate](C){}
		(0,1.3)node[coordinate](AA){}(1,1.3)node[coordinate](BB){}(2,1.3)node[coordinate](CC){}
		(.5,.7)node[coordinate](ABU){}(1.5,.7)node[coordinate](BCU){}
		(.5,-.7)node[coordinate](ABD){}(1.5,-.7)node[coordinate](BCD){}
		(1,.2)node[coordinate](ABCU){}(1,-.2)node[coordinate](ABCD){}
		(1,0)node[coordinate](Z){};
	\foreach\xa/\xb/\xc/\xd/\xe in
		 {A/Z/.5//left,Z/AA/.5/$#7$/left,B/BCD/.3//left,C/BCD/.3//right,BCD/BCU/.5/$#6$/right,BCU/BB/.3//left,BCU/CC/.3//right}
		{\draw[trivalent](\xa)..controls+(.03,\xc)and+(.03,-\xc)..node[basiclabel,\xe,inner sep=3pt,scale=.8]{\xd}(\xb);}
	\draw[trivalent]
		(AA)--node[midsmx]{$#1$}(0,2.3)..controls+(.03,.3)and+(.03,.3)..(-1,2.3)
			..controls+(-.1,-1)and+(.1,1)..(-1,-1.3)..controls+(-.03,-.3)and+(-.03,-.3)..(A)
		(BB)--node[midsmx]{$#2$}(1,2.3)..controls+(.06,.6)and+(.06,.6)..node[very near end,left,scale=.8]{$#4$}(-1.5,2.3)
			..controls+(-.1,-1)and+(.1,1)..(-1.5,-1.3)..controls+(-.06,-.6)and+(-.06,-.6)..(B)
		(CC)--node[midsmx]{$#3$}(2,2.3)..controls+(.05,.5)and+(.05,.5)..node[very near end,right,scale=.8]{$#5$}(3,2.3)
			..controls+(-.1,-1)and+(.1,1)..(3,-1.3)..controls+(-.05,-.5)and+(-.05,-.5)..(C);
	}}
	
\def\cfrelbb(#1,#2,#3)(#4,#5,#6,#7,#8,#9){
	\tikz[heighttwo,yscale=.7,shift={(0,1)}]{
	\draw
		(0,-1.3)node[coordinate](A){}(1,-1.3)node[coordinate](B){}(2,-1.3)node[coordinate](C){}
		(0,1.3)node[coordinate](AA){}(1,1.3)node[coordinate](BB){}(2,1.3)node[coordinate](CC){}
		(.5,.7)node[coordinate](ABU){}(1.5,.7)node[coordinate](BCU){}
		(.5,-.7)node[coordinate](ABD){}(1.5,-.7)node[coordinate](BCD){}
		(1,.2)node[coordinate](ABCU){}(1,-.2)node[coordinate](ABCD){}
		(1,0)node[coordinate](Z){};
	\foreach\xa/\xb/\xc/\xd/\xe in
		{A/ABCD/.5//left,B/BCD/.3//left,C/BCD/.3//right,BCD/ABCD/.3/$#7$/right,ABCD/ABCU/.2/$#8$/left,
			ABCU/AA/.5/$#9$/left,ABCU/BCU/.3/$#6$/right,BCU/BB/.3//left,BCU/CC/.3//right}
		{\draw[trivalent](\xa)..controls+(.03,\xc)and+(.03,-\xc)..node[basiclabel,\xe,inner sep=5pt,scale=.7]{\xd}(\xb);}
	\draw[trivalent]
		(AA)--node[midsmx]{$#1$}(0,2.3)..controls+(.03,.3)and+(.03,.3)..(-1,2.3)
			..controls+(-.1,-1)and+(.1,1)..(-1,-1.3)..controls+(-.03,-.3)and+(-.03,-.3)..(A)
		(BB)--node[midsmx]{$#2$}(1,2.3)..controls+(.06,.6)and+(.06,.6)..node[very near end,left,scale=.8]{$#4$}(-1.5,2.3)
			..controls+(-.1,-1)and+(.1,1)..(-1.5,-1.3)..controls+(-.06,-.6)and+(-.06,-.6)..(B)
		(CC)--node[midsmx]{$#3$}(2,2.3)..controls+(.05,.5)and+(.05,.5)..node[very near end,right,scale=.8]{$#5$}(3,2.3)
			..controls+(-.1,-1)and+(.1,1)..(3,-1.3)..controls+(-.05,-.5)and+(-.05,-.5)..(C);
	}}
	
\def\cfrelbc(#1,#2,#3)(#4,#5,#6,#7,#8,#9){
	\tikz[heighttwo,yscale=.7,shift={(0,1)}]{
	\draw
		(0,-1.3)node[coordinate](A){}(1,-1.3)node[coordinate](B){}(2,-1.3)node[coordinate](C){}
		(0,1.3)node[coordinate](AA){}(1,1.3)node[coordinate](BB){}(2,1.3)node[coordinate](CC){}
		(.5,.7)node[coordinate](ABU){}(1.5,.7)node[coordinate](BCU){}
		(.5,-.7)node[coordinate](ABD){}(1.5,-.7)node[coordinate](BCD){}
		(1,.2)node[coordinate](ABCU){}(1,-.2)node[coordinate](ABCD){}
		(1,0)node[coordinate](Z){};
	\foreach\xa/\xb/\xc/\xd/\xe in
		{A/ABD/.3//left,B/ABD/.3//left,ABD/ABCD/.3/$#7$/left,C/ABCD/.5//right,ABCD/ABCU/.2/$#8$/right,
			ABCU/ABU/.2/$#6$/left,ABU/AA/.3/$#9$/left,ABU/BB/.3//left,ABCU/CC/.5//right}
		{\draw[trivalent](\xa)..controls+(.03,\xc)and+(.03,-\xc)..node[basiclabel,\xe,inner sep=5pt,scale=.7]{\xd}(\xb);}
	\draw[trivalent]
		(AA)--node[midsmx]{$#1$}(0,2.3)..controls+(.03,.3)and+(.03,.3)..(-1,2.3)
			..controls+(-.1,-1)and+(.1,1)..(-1,-1.3)..controls+(-.03,-.3)and+(-.03,-.3)..(A)
		(BB)--node[midsmx]{$#2$}(1,2.3)..controls+(.06,.6)and+(.06,.6)..node[very near end,left,scale=.8]{$#4$}(-1.5,2.3)
			..controls+(-.1,-1)and+(.1,1)..(-1.5,-1.3)..controls+(-.06,-.6)and+(-.06,-.6)..(B)
		(CC)--node[midsmx]{$#3$}(2,2.3)..controls+(.05,.5)and+(.05,.5)..node[very near end,right,scale=.8]{$#5$}(3,2.3)
			..controls+(-.1,-1)and+(.1,1)..(3,-1.3)..controls+(-.05,-.5)and+(-.05,-.5)..(C);
	}}
	
\def\cfrelbd(#1,#2)(#3,#4,#5){
	\tikz[heighttwo,yscale=.7,shift={(0,1)}]{
	\draw
		(0,-1.3)node[coordinate](A){}(2,-1.3)node[coordinate](B){}
		(0,1.3)node[coordinate](AA){}(2,1.3)node[coordinate](BB){}
		(.5,.7)node[coordinate](ABU){}(1.5,.7)node[coordinate](BCU){}
		(.5,-.7)node[coordinate](ABD){}(1.5,-.7)node[coordinate](BCD){}
		(1,.2)node[coordinate](ABCU){}(1,-.2)node[coordinate](ABCD){}
		(1,0)node[coordinate](Z){};
	\foreach\xa/\xb/\xc/\xd/\xe in
		{A/ABCD/.5//left,C/ABCD/.5//right,ABCD/ABCU/.2/$#5$/right,ABCU/AA/.5//left,ABCU/CC/.5//right}
		{\draw[trivalent](\xa)..controls+(.03,\xc)and+(.03,-\xc)..node[basiclabel,\xe,inner sep=1pt,scale=.7]{\xd}(\xb);}
	\draw[trivalent]
		(AA)--node[midsmx]{$#1$}(0,2.3)..controls+(.03,.3)and+(.03,.3)..node[very near end,left,scale=.8]{$#3$}(-1,2.3)
			..controls+(-.1,-1)and+(.1,1)..(-1,-1.3)..controls+(-.03,-.3)and+(-.03,-.3)..(A)
		(BB)--node[midsmx]{$#2$}(2,2.3)..controls+(.05,.5)and+(.05,.5)..node[very near end,right,scale=.8]{$#4$}(3,2.3)
			..controls+(-.1,-1)and+(.1,1)..(3,-1.3)..controls+(-.05,-.5)and+(-.05,-.5)..(B);
	}}
	
\def\cfrelca(#1,#2)(#3,#4,#5)(#6,#7,#8){
	\tikz[heighttwo,yscale=.7,shift={(0,1)}]{
	\draw
		 (0,-1.3)node[coordinate](A){}(1,-1.3)node[coordinate](B){}(2,-1.3)node[coordinate](C){}(3,-1.3)node[coordinate](D){}
		 (0,1.3)node[coordinate](AA){}(1,1.3)node[coordinate](BB){}(2,1.3)node[coordinate](CC){}(3,1.3)node[coordinate](DD){}
		(.5,.7)node[coordinate](ABU){}(.5,-.7)node[coordinate](ABD){}
		(1,.4)node[coordinate](ACU){}(1,-.4)node[coordinate](ACD){}
		(1.5,.2)node[coordinate](ABCDU){}(1,-.2)node[coordinate](ABCDD){}
		(2,.4)node[coordinate](BDU){}(2,-.4)node[coordinate](BDD){}
		(2.5,.7)node[coordinate](CDU){}(2.5,-.7)node[coordinate](CDD){};
	\foreach\xa/\xb/\xc/\xd/\xe in
		{A/ACD/.5//left,C/ACD/.5//right,B/BDD/.5//left,D/BDD/.5//right,ACD/ACU/.3/$c$/left,BDD/BDU/.3/$c'$/right,
			ACU/AA/.5//left,ACU/CC/.5//right,BDU/BB/.5//left,BDU/DD/.5//right}
		{\draw[trivalent](\xa)..controls+(.03,\xc)and+(.03,-\xc)..node[basiclabel,\xe,inner sep=2pt]{\xd}(\xb);}
	\draw[trivalent]
		(AA)--node[midsmx]{$#1$}(0,2.3)..controls+(.03,.3)and+(.03,.3)..(-1,2.3)
			 ..controls+(-.1,-1)and+(.1,1)..node[basiclabel,pos=.4,right]{$a$}(-1,-1.3)..controls+(-.03,-.3)and+(-.03,-.3)..(A)
		(BB)--node[midsmx]{$#1$}(1,2.3)..controls+(.06,.6)and+(.06,.6)..(-1.5,2.3)
			..controls+(-.1,-1)and+(.1,1)..node[basiclabel,at start,left]{$a'$}(-1.5,-1.3)..controls+(-.06,-.6)and+(-.06,-.6)..(B)
		(CC)--node[midsmx]{$#2$}(2,2.3)..controls+(.06,.6)and+(.06,.6)..(4.5,2.3)
			..controls+(-.1,-1)and+(.1,1)..node[basiclabel,at start,right]{$b$}(4.5,-1.3)..controls+(-.06,-.6)and+(-.06,-.6)..(C)
		(DD)--node[midsmx]{$#2$}(3,2.3)..controls+(.03,.3)and+(.03,.3)..(4,2.3)
			 ..controls+(-.1,-1)and+(.1,1)..node[basiclabel,pos=.4,left]{$b'$}(4,-1.3)..controls+(-.03,-.3)and+(-.03,-.3)..(D);
	}}

\def\cfrelcb(#1,#2)(#3,#4)(#5,#6)(#7,#8,#9){
	\tikz[heighttwo,yscale=.7,shift={(0,1)}]{
	\draw
		 (0,-1.3)node[coordinate](A){}(1,-1.3)node[coordinate](B){}(2,-1.3)node[coordinate](C){}(3,-1.3)node[coordinate](D){}
		 (0,1.3)node[coordinate](AA){}(1,1.3)node[coordinate](BB){}(2,1.3)node[coordinate](CC){}(3,1.3)node[coordinate](DD){}
		(.5,.7)node[coordinate](ABU){}(.5,-.7)node[coordinate](ABD){}
		(1,.4)node[coordinate](ACU){}(1,-.4)node[coordinate](ACD){}
		(1.5,.2)node[coordinate](ABCDU){}(1.5,-.2)node[coordinate](ABCDD){}
		(2,.4)node[coordinate](BDU){}(2,-.4)node[coordinate](BDD){}
		(2.5,.7)node[coordinate](CDU){}(2.5,-.7)node[coordinate](CDD){};
	\foreach\xa/\xb/\xc/\xd/\xe in
		{A/ABD/.5//left,B/ABD/.5//right,C/CDD/.5//left,D/CDD/.5//right,
			ABD/ABCDD/.3/$#7_2$/left,CDD/ABCDD/.3/$#8_2$/right,ABCDD/ABCDU/.3/$#9$/right,ABU/ABCDU/-.3/$#7_1$/left,
			CDU/ABCDU/-.3/$#8_1$/right,ABU/AA/.5//left,ABU/BB/.5//right,CDU/CC/.5//left,CDU/DD/.5//right}
		{\draw[trivalent](\xa)..controls+(.03,\xc)and+(.03,-\xc)..node[basiclabel,\xe,near start,scale=.7]{\xd}(\xb);}
	\draw[trivalent]
		(AA)--node[midsmx]{$#1$}(0,2.3)..controls+(.03,.3)and+(.03,.3)..(-1,2.3)
			 ..controls+(-.1,-1)and+(.1,1)..node[basiclabel,pos=.4,right]{$a$}(-1,-1.3)..controls+(-.03,-.3)and+(-.03,-.3)..(A)
		(BB)--node[midsmx]{$#1$}(1,2.3)..controls+(.06,.6)and+(.06,.6)..(-1.5,2.3)
			..controls+(-.1,-1)and+(.1,1)..node[basiclabel,at start,left]{$a'$}(-1.5,-1.3)..controls+(-.06,-.6)and+(-.06,-.6)..(B)
		(CC)--node[midsmx]{$#2$}(2,2.3)..controls+(.06,.6)and+(.06,.6)..(4.5,2.3)
			..controls+(-.1,-1)and+(.1,1)..node[basiclabel,at start,right]{$b$}(4.5,-1.3)..controls+(-.06,-.6)and+(-.06,-.6)..(C)
		(DD)--node[midsmx]{$#2$}(3,2.3)..controls+(.03,.3)and+(.03,.3)..(4,2.3)
			 ..controls+(-.1,-1)and+(.1,1)..node[basiclabel,pos=.4,left]{$b'$}(4,-1.3)..controls+(-.03,-.3)and+(-.03,-.3)..(D);
	}}



%% file: spinsl2-intro.tex
\section{Introduction}

The purpose of this chapter is to demonstrate the utility of a graphical calculus in the algebraic study of $\sltc$-representations of the fundamental group of an oriented surface of Euler characteristic $-1$.

Let $\mathsf{F}_2$ be a rank 2 free group, the fundamental group of both the three-holed sphere and the one-holed torus. The set $\mathcal{R}=\Hom(\mathsf{F}_2,\sltc)$ of representations inherits the structure of an algebraic set from $\sltc$. The subset of representations that are \emph{completely reducible},\index{completely reducible representation} denoted by $\mathcal{R}^{ss}$, have closed orbits under conjugation. Consequently, the orbit space $\mathcal{R}^{ss}/\sltc=\mathcal{R}/\!/\sltc$ is an algebraic set referred to as the \emph{character variety} \index{character variety}. The character variety encodes both Teichm\"uller Space and moduli of geometric structures \cite{Gol2}.

Graphs known as \emph{spin networks} permit a concise description of a natural additive basis for the coordinate ring of the
character variety
	$$\C[\mathcal{R}/\!/\sltc]=\C[\mathcal{R}]^{\sltc}.$$
We will refer to the basis elements as \emph{central functions}\index{central function}. The central functions are indexed by Clebsch-Gordan injections
	$$V_c\injects V_a\tensor V_b,$$
where $V_c=\Sym^c(\C^2)$ denotes an irreducible representation of $\sltc$. Our main results use the spin network calculus to
describe a strong symmetry within the central function basis, a graphical means of computing the product of two central
functions, and an algorithm for computing central functions. This provides a concrete description of the regular functions on the $\sltc$-character variety of $\fg$ and a new proof of a classical result of Fricke, Klein, and Vogt.

We are motivated by a greater understanding of the invariant ring, and the subsequent knowledge of various geometric objects of interest encoded within the character variety. Consequently, the main results in this chapter concern the structure of the central function basis. The results and methods of this chapter may also provide new insight into gauge theoretic questions. However, we are most interested in a methodology and point of view that allows for generalizations to other Lie groups and other surface groups.

\medskip\noindent\emph{History of Central Functions and Spin Networks.}

The first reference to the central function basis in the literature appears in \cite{Ba}, where Baez used spin networks to
describe a basis of quantum mechanical ``state vectors.'' He considered the basis abstractly, showing that the space of square integrable functions on a related space of connections modulo gauge transformations is spanned by a set of labelled graphs. He also demonstrated that the basis is orthonormal with respect to the $L^2$ inner product. His basis, when restricted to $\mathrm{SU}(2)$, is precisely the one under consideration here.

More recently, Florentino, Mour\~ao, and Nunes use a like basis to produce distributions related to geometric quantization of moduli spaces of flat connections on a surface \cite{FMN}. Adam Sikora has also used spin networks to study the character variety for $\mathrm{SL}(n,\C)$, although without using the central function basis \cite{Sik}. The construction of arbitrary rank $\sltc$ central functions is described in \cite{Pet}, while much of the diagrammatic theory required for the $\mathrm{SL}(n,\C)$ case is covered in \cite{Cvi,Cvi2,Pet,Sik}.

The history of the diagrammatic calculus in this chapter is hard to trace, due to the historical difficulty in publishing papers making extensive use of figures. While it is likely that many works on diagrammatic notation have been lost over the years, the specific notation used in this chapter is due to Roger Penrose. In a 1981 letter to Predrag Cvitanovi\'c, a physicist who also used diagrams extensively, Penrose recalls developing the notation in the early 1950s while ``trying to cope with Hodge's lectures on differential geometry'' \cite{Pen2}.

Diagrammatic notations have also played an important role in modern physics. Feynman diagrams are probably the most famous example, but spin networks have also been used for many years, as a graphical description of quantum angular momentum \cite{Pen1}. The use of diagrams in physics is probably best summarized in \cite{Ste}. Cvitanovi\'c also has a thorough description of such notations, which he calls \emph{birdtracks}\index{birdtracks} in \cite{Cvi,Cvi2}. In his work, birdtracks play a starring role in a new classification of semi-simple Lie algebras. Using primitive invariants, which have unique diagrammatic depictions, the exceptional Lie algebras arise in a single series in a construction that he calls the ``Magic Triangle.''

\medskip\noindent\emph{The remainder of this chapter is organized as follows.}

Section \ref{prelim} gives some basic definitions and results from invariant theory, as well as a short history of
$\mathrm{SL}(2,\C)$ invariant theory. It also covers necessary material from representation theory.

In Section \ref{spincalc}, we introduce spin networks, which are special types of graphs that may be identified with functions between tensor powers of $\C^2$. We give a full treatment of the {\it spin network calculus}, a powerful means for working with regular functions on $\mathcal{R}/\!/\sltc$.

Section \ref{decomp} begins by constructing an additive basis for $\C[\mathcal{R}/\!/\sltc]$. This basis, denoted by $\{\ch
abc\}$, is indexed by triples of nonnegative integers $(a,b,c)$ satisfying the \emph{admissibility condition}:
  $$\tfrac12(-a+b+c), \:\tfrac12(a-b+c), \:\tfrac12(a+b-c)\in\N.$$
The functions $\ch abc\in\C[\mathcal{R}/\!/\sltc]$ are central in%
  $$\End(V_c)\injects\End(V_a)\tensor\End(V_b),$$
and are referred to as \emph{central functions}. The construction of the central function basis uses the decomposition
  $$\C[\sltc]\isom\sum_{n\ge0}V_n^*\tensor V_n.$$
We include a constructive proof of this decomposition, since it is hard to find in the literature. The section concludes by
examining the $\sltc$-central functions of a rank one free group.

Section \ref{rank2} contains the main results of this chapter, which concern the case of a rank two free group. In this case,
central functions may be written as polynomials in three trace variables, a consequence of a theorem due to Fricke, Klein, and Vogt \cite{FK,Vogt}. The results we prove are summarized below.
\begin{itemize}
  \item Theorem \ref{symmetry} describes a symmetry property of the central function basis: permuting the indices of a central function is equivalent to permuting the variables of its polynomial representation.
  \item Corollary \ref{recurrence} states that, with an appropriate definition of rank, any central function may be written in terms of at most four central functions of lower rank:
    \begin{multline*}
      \ch abc=x\cdot\ch{a-1}b{c-1}-\tfrac{(a+b-c)^2}{4a(a-1)}\ch{a-2}bc%
        -\tfrac{(-a+b+c)^2}{4c(c-1)}\ch ab{c-2}\\-\tfrac{(a+b+c)^2(a-b+c-2)^2}{16a(a-1)c(c-1)}\ch{a-2}b{c-2}.
    \end{multline*}
    Together with Theorem \ref{symmetry}, this result gives an algorithm for computing central functions explicitly.
  \item Proposition \ref{monic} states that central functions are monic, and gives the leading term
    of the central function $\ch abc$.
  \item Proposition \ref{grading} describes a $\Z_2\times\Z_2$ grading on the central function basis.
  \item Theorem \ref{mult} gives the coefficients in the expression of the product of two central functions as a sum of central functions, and therefore a precise description of the ring structure of $\C[\mathcal{R}]^{\sltc}$ in terms of central functions.
\end{itemize}
Finally, as another consequence of the recurrence relation and Theorem \ref{symmetry}, we provide a new constructive proof of the following classical theorem \cite{FK,Vogt}:
\begin{FKV}[Fricke-Klein-Vogt Theorem]
	Let $G=\sltc$ act on $G\times G$ by simultaneous conjugation.  Then
		$$\mathbb{C}[G\times G]^G \isom \mathbb{C}[t_x,t_y,t_z],$$
	the complex polynomial ring in three indeterminates. In particular, every regular function $f:\sltc\times\sltc\to\C$ 		 satisfying%
		$$f(\bfx_1,\bfx_2)=f(g\bfx_1g^{-1},g\bfx_2g^{-1})\qquad\text{for all }g\in\sltc,$$%
	can be written uniquely as a polynomial in the three trace variables $x=\tr(\bfx_1)$, $y=\tr(\bfx_2)$, and
	$z=\tr(\bfx_1\inv{\bfx_2})$.
\end{FKV}

%% file: spinsl2-prelim-sean.tex
\section{Preliminaries} \label{prelim}

\subsection{Algebraic Structure of the Character Variety $\mathcal{R}/\!/\sltc$}%

The group $G=\mathrm{SL}(2,\cb)$ has the structure of an irreducible algebraic set, since it is the zero set of the irreducible
polynomial $\det(\xb)-1$.  Since the product of two varieties is again a variety, the \emph{representation variety} \index{representation variety}
$\mathcal{R}=\Hom({\sf F}_2,G)\isom G\times G$ of a rank $2$ free group ${\sf F}_2$ is an irreducible algebraic set as well. The
\emph{coordinate ring} of $\mathcal{R}$ is
  $$\C[\mathcal{R}]=\frac{\C[x^k_{ij} : 1\leq i,j,k \leq2]}{(\det(\xb_1)-1,\det(\xb_2)-1)}.$$
Stated otherwise, it is the free commutative polynomial ring in $8$ indeterminates over $\C$ subject to the ideal generated by
the two polynomials $\det(\xb_k)-1$, where $\xb_k=(x^k_{ij})$ are called {\it generic matrices} \index{generic matrices}.

There is an action of $G$ on $\mathcal{R}$ by simultaneous conjugation. Given $(\bfx_1,\bfx_2)\in G\times G$, then
	$$g\cdot(\bfx_1,\bfx_2)=(g\bfx_1g^{-1},g\bfx_2g^{-1}).$$
This is a \emph{polynomial action}, since $\mathcal{R}\times G\to\mathcal{R}$ is a regular mapping.
  \begin{defn}
    The \emph{ring of invariants} \index{ring of invariants} $\C[\mathcal{R}]^G$ consists of elements of the coordinate ring 			
    $\C[\mathcal{R}]$ which are invariant under the action of simultaneous conjugation:
      $$\C[\mathcal{R}]^G=\{f\in \C[\mathcal{R}] : g\cdot f=f\}.$$
  \end{defn}
Recall that an algebraic group is {\it linearly reductive}\index{reductive} if its finite dimensional rational representations are decomposable as direct sums of irreducible representations. Since $G=\mathrm{SL(2,\C)}$ is linearly reductive, the ring of invariants $\C[\mathcal{R}]^G=\{f\in \C[\mathcal{R}] : g\cdot f=f\}$ is finitely generated \cite{Dol}. This implies that the space of maximal ideals of $\C[R]^G$ is also an irreducible algebraic set, permitting the following definition:
  \begin{defn}
    The $G$-{\it character variety}\index{character variety} of ${\sf F}_2$ is the space of maximal ideals
      $$\XX=\mathrm{Spec}_{max}(\cb[\mathcal{R}]^G)=\mathcal{R}/\!/G.$$
  \end{defn}
The character variety $\XX$ is identified with conjugacy classes of {\it completely reducible} representations in $\mathcal{R}$
\cite{Art, Pro2}.  Procesi \cite{Pro} has shown that $\cb[\mathcal{R}]^G$ is generated by traces of products of matrix variables
of word length less than or equal to three \cite{Pro}.  Hence $\C[\XX]$ is generated, although not minimally, by
$$\{\tr(\bfx_1),\tr(\bfx_2),\tr(\bfx_1\bfx_2),\tr(\bfx_1\bfx_2^2),\tr(\bfx_2\bfx_1^2)\}.$$

\subsection{History of $\mathrm{SL}(2,\C)$ Invariant Theory}\index{invariant theory}

The invariant theory of $\mathrm{SL}(2,\cb)$ has a long history. Two pioneering papers on the subject were authored by Vogt in
1889 \cite{Vogt}, and by Fricke and Klein in 1896 \cite{FK}. Both investigated the invariants of pairs of unimodular $2\times2$
matrices with respect to simultaneous conjugation.  They showed this ring of invariants to be the free commutative polynomial
ring in three indeterminants, given by the trace of each generic matrix and the trace of their product. This chapter concludes with
a reproof of this classical result using the spin network calculus.

In 1972, Horowitz investigated the algebraic structure of this ring, saying that Fricke's approach was principally analytic, and
partially incomplete \cite{Hor}.  In 1980, Magnus made clear the priority of Vogt's approach \cite{Vogt} and worked out the
defining polynomial relations for an arbitrary number of matrices under simultaneous conjugation \cite{Mag}.  In 1983, Culler and
Shalen defined the character variety and showed that it is in fact an algebraic set \cite{CS}; the set is the image under a
``trace'' map. Gonz\'{a}lez-Acu\~{n}a and Montesinos-Amilibia showed in 1993 that the relations of Magnus in fact determine the
algebraic set that Culler and Shalen had defined \cite{GM}.  In 2001, Sikora, using results of Procesi \cite{Pro}, showed that
the character variety of $\mathrm{SL}(n,\C)$ can be realized as spaces of graphs subject to topologically motivated relations
\cite{Sik}. These graphs correspond to the spin networks discussed in this chapter when $n=2$.

Closely related is the ring of invariants of \emph{arbitrary} generic $2\times 2$ matrices under simultaneous conjugation. The
works of Procesi (1976) and Razmyslov (1974) generalized the work above to the case of $n\times n$ matrices \cite{Pro,Raz}, and
showed that the invariant ring is generated by traces of words in generic matrices. Methods from geometric invariant theory (see
Dolgachev \cite{Dol}) show that the character variety is the variety whose coordinate ring is the ring of invariants. Restricting
to unimodular matrices gives like results for the unimodular ring of invariants. From this point of view, the character variety
begins as an algebraic set and so is obviously closed.  However, the defining relations and minimal generators are not at all
obvious.

A central question in invariant theory is a description of the generators and relations of an invariant ring. Indeed, a theorem
that characterizes the generators of an invariant ring is called a {\it first fundamental theorem}, and a theorem giving the
relations is called a {\it second fundamental theorem}. In \cite{Pro, Raz} both Procesi and Razmyslov gave the two fundamental
theorems, although they offered only sufficient generators and an implicit description of the relations.

It is much more difficult to determine \emph{minimal} generators and \emph{explicit} relations. In this more general context,
which bears strongly on the unimodular case, minimal generators and defining relations for the invariants of an arbitrary number
of generic $2\times 2$ matrices were found only recently by Drensky in 2003 \cite{Dre}.


\subsection{Representation Theory of $\sltc$}\label{reptheory}
The coordinate ring $\C[G]$ decomposes into a direct sum
of tensor products of the finite-dimensional irreducible representations of $G$. We will use this
decomposition, given explicitly by Theorem \ref{pwdecomp}, to understand $\C[\XX]$. To this end, we
review the representation theory of $G$ (see \cite{BtD,Dol,FH}).

The symmetric powers of the standard representation of $G$ are all irreducible representations and
moreover they comprise a complete list. Let $V_0=\cb =V_0^*$ be the trivial representation of $G$.
Denote the standard basis for $\C^2$ by $e_1=\tmxt10$ and $e_2=\tmxt01$, and the dual basis by
$e_1^*=e_1^T$ and $e_2^*=e_2^T$. Then the standard representation and its dual are
$$V=V_1 = \cb e_1\oplus \cb e_2\quad\text{and}\quad V^*=V_1^*=\cb e_1^*\oplus \cb e_2^*,$$ respectively.
Denote the symmetric powers of these representations by
$$V_n=\mathrm{Sym}^n(V) \textrm{ and } V_n^*=\mathrm{Sym}^n(V^*).$$%
Since $V_n$ admits an invariant non-degenerate bilinear form, $V_n\isom (V_n)^*$.

Moreover, $V_n^*$ is naturally isomorphic to $(V_n)^*$, so elements in $V_n$ pair with elements in
$V_n^*$. Denote the projection of $v_1\otimes v_2 \otimes \cdots \otimes v_n \in V^{\otimes n}$ to
$V_n$ by $v_1\circ v_2 \circ \cdots \circ v_n$. There exist bases for $V_n$ and $V_n^*$, given by
the elements%
\begin{align*}
    {\sf n}_{n-k}&=e_1^{n-k}e_2^k=\underbrace{e_1\circ e_1\circ \cdots \circ e_1}_{n-k}\circ
        \underbrace{e_2\circ e_2\circ \cdots \circ e_2}_{k}\quad\text{and}\\
    {\sf n}_{n-k}^*&=(e_1^*)^{n-k}(e_2^*)^k=\underbrace{e_1^*\circ e_1^*\circ \cdots \circ e_1^*}_{n-k}\circ
        \underbrace{e_2^*\circ e_2^*\circ \cdots \circ e_2^*}_{k}\:,
\end{align*}
respectively, where $0\le k\le n$. In these terms, this pairing is given by
$${\sf n}^*_{n-k}(v_1\circ v_2 \circ \cdots \circ v_n)=
\frac{1}{n!} \sum_{\sigma \in \Sigma_n}({\sf n}_{n-k})^*(v_{\sigma(1)}\otimes v_{\sigma(2)} \otimes
\cdots \otimes v_{\sigma(n)}),$$ where $\Sigma_n$ is the symmetric group on $n$ elements. In
particular,
$${\sf n}^*_{n-k}({\sf n}_{n-l})=\frac{(n-k)!k!}{n!}\delta_{kl}=\raisebox{2pt}{$\delta_{kl}$}\!\Big/\!\raisebox{-2pt}{$\tbinom{n}{k}$}.$$

Let $g=\tmx{g_{11} & g_{12}}{g_{21} & g_{22}}\in G$. Then the $G$-action on $V_n$ is given by
\begin{align*}
  g\cdot {\sf n}_{n-k}&=(g_{11}e_1+g_{21}e_2)^{n-k}(g_{12}e_1+g_{22}e_2)^k\\%
    &=\sum_{\substack{0\le j\le n-k\\0\le i \le k}}\tbinom{n-k}{j}\tbinom{k}{i}
        \left(g_{11}^{n-k-j}g_{12}^{k-i}g_{21}^{j}g_{22}^{i}\right){\sf n}_{n-(i+j)}.
\end{align*}
For the dual, $G$ acts on $V_n^*$ in the usual way:
\begin{displaymath}
(g \cdot {\sf n}^*_{n-k})(v)= {\sf n}^*_{n-k}(g^{-1}\cdot v ) \textrm{ for } v\in V_n.
\end{displaymath}

The tensor product $V_a\tensor V_b,$ where $a,b \in \N $, is also a representation of $G$ and
decomposes into irreducible representations as follows:
\begin{prop}[Clebsch-Gordan formula]\label{cgform}\index{Clebsch-Gordan formula}
  $$V_a\otimes V_b\isom\bigoplus^{\mathrm{min}(a,b)}_{j=0} V_{a+b-2j}.$$
\end{prop}

Finally, we give several versions of Schur's Lemma, which will be used frequently.
\begin{prop}[Schur's Lemma]\label{schurs}\index{Schur's Lemma}
Let $G$ be a group, $V$ and $W$ representations of $G$, and $f\in\Hom_G(V,W)$ with $f\neq0$.
\begin{enumerate}
\item If $V$ is irreducible, then $f$ is injective.
\item If $W$ is irreducible, then $f$ is surjective.
\item If $V=W$ is irreducible, then $f$ is a homothety.
\item Suppose $V,W$ are irreducible:
\subitem if $V\isom W$, then $\dim_\C\mathrm{Hom}_G(V,W)=1$;%
\subitem if $V\not\isom W$, then $\dim_\C\mathrm{Hom}_G(V,W)=0$.
\end{enumerate}
\end{prop}
See \cite{BtD} or \cite{CSM} for proof of Propositions \ref{cgform} and \ref{schurs}.


%% file: spinsl2-prelim-elisha.tex

\section{The Spin Network Calculus} \label{spincalc}

This section provides a self-contained introduction to spin networks and the spin network calculus. Our treatment employs a
nonstandard definition of spin networks which is more natural when working with traces. This definition leads to different
versions of the usual spin network relations in the literature \cite{CFS,Cvi,Cvi2,Kau,Pen1,Ste}.

\subsection{Spin Networks and Representation Theory}
At its heart, a spin network is a graph that is identified with a specific function between tensor
powers of $V=\C^2$, the standard $\sltc$ representation.

In order for this function to be well-defined, the edges incident to each vertex of the spin
network must have a cyclic ordering. This ordering is often called a \emph{ciliation},\index{ciliation} since it is
represented on paper by a small mark drawn between two of the edges. The edges adjacent to a
ciliated vertex are ordered by proceeding in a clockwise fashion from this mark. For example, in
the degree 2 case, there are two possible ciliations: $\upcl_2^1$ and
$\upcr^{\!\!2}_{\!\!1}.$

\begin{defn}
A \emph{spin network} \index{spin network} $\S$ is a graph with vertex set $\S_i\sqcup\S_o\sqcup\S_v$ consisting of
degree 1 `inputs' $\S_i$, degree 1 `outputs' $\S_o$ and degree 2 `ciliated vertices' $\S_v$. If
there are $k_i=|\S_i|$ inputs and $k_o=|\S_o|$ outputs, then $\S$ is identified with a function
$f_\S:\vprod{k_i}\to\vprod{k_o}$. If the spin network is \emph{closed}, meaning $k_i=0=k_o$, it is
identified with a complex scalar $f_\S\in\C$.
\end{defn}

Spin networks are drawn in \emph{general position} inside an oriented rectangle with inputs at the
bottom and outputs at the top. This convention allows us to equate the composition of functions
$f_{\S'}\circ f_\S$ with the concatenation of diagrams $\S'\circ\S$ formed by placing $\S'$ on top
of $\S$.

For example, the following spin network has two ciliated vertices and represents a function from
$\vprod5\to\vprod3$:
	$$\underbrace{\overbrace{\tbigexccut}^{\text{3 outputs}}}_{\text{5 inputs}}
		=\left(\pup\:\capdot\:\:\upcr\:\upcl\right)\circ\left(\px\:\:\ccdot\:\pup\right)
		\circ\left(\perm{1/1,2/2,3/3,4/5,5/4}\right).$$
Note that the marks on the local extrema do not indicate vertices of the graph, but are indicators of how to decompose the graph.

Since spin networks are just graphs with ciliations, it does not matter how the graph is
represented inside the square. Strands may be moved about freely and ciliations may ``slide'' along
the strands. As long as the endpoints remain fixed, the underlying spin network does not change.

Let $v,w\in V$ and let $\{e_1,e_2\}$ be the standard basis for $\C^2$. The function $f_\S$ of a
spin network $\S$ is computed by decomposing $\S$ into the four \emph{spin network component maps}:
\index{spin network!component maps}
\begin{itemize}
\item the \emph{identity} $\pup:V\to V$,\quad$v\mapsto v$;
\item the \emph{cap} $\capdot:V\tensor V\to\C$,\quad$v\tensor w\mapsto v^Tw$ (inner product);
\item the \emph{cup} $\cupdot:\C\to V\tensor V$,\quad$1\mapsto e_1\tensor e_1+e_2\tensor e_2$;
\item the \emph{cap vertex} $\capcu:V\tensor V\to \C$,\quad$v\tensor w\mapsto\det[v\:w]$.
\end{itemize}
For example, since $\capcd$ and $\capcux$ are the same ciliated graph,
  $$\capcd(v\tensor w)=\capcux(v\tensor w)=\capcu\circ\px(v\tensor w)=\capcu(w\tensor
    v)=\det[w\:v].$$

The definition given here differs from the literature \cite{CFS,Kau,Pen1}. In particular, we omit the
$i=\sqrt{-1}$ factor in the definition of $\capcu$ to gain an advantage in trace calculations.
Also, the maps $\capdot$ and $\cupdot$ are included in order to simplify the proof that
$f_\S$ is well-defined.

\begin{thm}
The spin network function $f_\S$ is well-defined.
\begin{proof}
We need to show that every decomposition of $\S$ into the component maps gives the same function.

If $\S$ has $n$ ciliated vertices, then any decomposition of $\S$ into component maps has $n$
occurrences of $\capcu$. The remainder of the diagram consists of loops or arcs without any
vertices. Two corresponding arcs in different decompositions will differ only by the insertion or
deletion of a number of `kinks' of the form \kinkdot. Finally, since%
$$\kinkdot(v)=\capdot\:\pup\circ\pup\:\cupdot(v)=\pup(v)$$%
for all $v\in V$, these kinks do not change the resulting function. For alternate proofs, see
\cite{CFS,Kau}.
\end{proof}
\end{thm}

This theorem allows us to freely interpret a spin network $\S$ as a function. The computation of
$f_\S$ will be easier once the functions for a few simple spin networks are known.

\begin{prop} As spin network functions,
  \begin{enumerate}
    \item the swap $\px:V\tensor V\to V\tensor V$ takes $v\tensor w\mapsto w\tensor v$;
    \item the vertex on a straight line $\upcl:V\to V$ takes $v\mapsto\tmxt{0&-1}{1&0}v$;
    \item the vertex on a cup $\cupcu:\C\to V\tensor V$ takes $1\mapsto e_1\tensor e_2-e_2\tensor e_1$;
    \item with opposite ciliations, $\capcd=-\capcu$, $\upcr=-\upcl$, and $\cupcd=-\cupcu$.
  \end{enumerate}
\begin{proof}%
First (1) is the statement that crossings change only the order of the outputs. %
Statement (2) follows from, for $v=\tmxt{v^1}{v^2}$:
\begin{align*}
  \upcl(v)&=\left(\capcu\:\pup\right)\circ\left(\pup\:\cupdot\right)(v)
              =\left(\capcu\:\pup\right)(v\tensor e_1\tensor e_1+v\tensor e_2\tensor e_2)\\
          &=\det[v\: e_1] e_1+\det[v\: e_2] e_2
              =-v^2 e_1+v^1 e_2=\tmxt{0&-1}{1&0}v.
\end{align*}
Statement (3) is computed similarly, using the decomposition
$$\cupcu=\left(\pup\:\capcu\:\pup\right)\circ\left(\cupdot\:\:\cupdot\right).$$
Finally, (4) follows from the observation $\capcd=\capcux=-\capcu,$ which has
already been demonstrated.
\end{proof}
\end{prop}

Given these facts, the function of the earlier example can be computed. The reader may check that the function given by
	$$\tbigexc$$
takes $e_1\tensor e_2\tensor e_2\tensor e_1\tensor e_2$ to
$-e_2\tensor e_2\tensor e_2$.
 
The maps $\capdot$ and $\cupdot$ are unnecessary for trace computations, and so we make the following assumption:
\begin{conv}
For the remainder of this chapter, the set of ciliated vertices will \emph{coincide exactly with the
set of local extrema}. The ciliations are usually omitted, with the understanding that
\begin{align*}
                  \cupp&=\cupcu:1\mapsto e_1\tensor e_2-e_2\tensor e_1\\%
  \text{and}\qquad\capp&=\capcu:v\tensor w\mapsto\det[v\:w].%
\end{align*}
\end{conv}
Under this assumption, each straightened kink $\kink\like\pup$ introduces a sign, and more
generally
$$\kinkn{n}=(-1)^n\slantn{n}.$$
Thus, any diagram manipulation in which kinks are straightened must be done carefully.


Spin networks exhibit considerable symmetry, which can be exploited for calculations. For example:
\begin{prop}\label{reflprop}
Let $\S$ be a spin network with function $f_\S:\vprod{k_i}\to\vprod{k_o}$. Denote its images under
reflection \index{spin network!reflection} through vertical and horizontal lines by $\overleftrightarrow\S$ and $\S^\updownarrow$,
respectively. Then
$$f_{\overleftrightarrow\S}=(-1)^{|\S_v|} \overleftrightarrow{f_\S}:\vprod{k_i}\to\vprod{k_o},$$ %
where $|\S_v|$ is the number of local extrema in the diagram and $\overleftrightarrow f$ indicates that the ordering of inputs
and outputs is reversed. Also, $f_{\S^\updownarrow}=(f_\S)^*$ where
$$(f_\S)^*(v_1\dtensor v_{k_i})=\sum_{e_b\in\mathcal{B}(\vprod{k_i})}\left(f_\S(e_b)\cdot(v_1\dtensor v_{k_o})\right)e_b,$$%
where $\cdot$ indicates the dot product with respect to the standard basis for $\vprod{k_o}$ and $\mathcal{B}(\vprod{k_i})$ is
the basis for $\vprod{k_i}$. That is, $(f_\S)^*$ and $f_\S$ are dual with respect to the standard inner product on $V$.
\begin{proof}
The first statement is an extension of the fact that reflecting $\capcd$ through a vertical
line gives $\capcu=-\capcd.$

For the second statement, consider $\S=\cupcu$. If $v_i=\tmxt{v_i^1}{v_i^2}$, then%
\begin{align*}
  (f_\S)^*(v_1\tensor v_2)&=\cupp(1)\cdot(v_1\tensor v_2)=(e_1\tensor e_2-e_2\tensor e_1)\cdot(v_1\tensor v_2)\\ %
                        &=v_1^1v_2^2-v_1^2v_2^1=\det[v_1\: v_2]=\capcu(v_1\tensor v_2).
\end{align*}
This computation, together with the corresponding one for $\S=\capcu$, are sufficient to prove
the second claim (see \cite{Pet} for details).
\end{proof}
\end{prop}

The next theorem, which follows from Proposition \ref{reflprop}, describes how to apply these symmetries to relations among spin
networks:
\begin{thm}[Spin Network Reflection Theorem]\label{reflector}\index{spin network!reflection} 
A relation $$\sum_m\alpha_m\S^m=0$$ among some collection of spin networks $\{\S^m\}$ is equivalent
to the same relation for the vertically reflected spin networks $\{{\S^\updownarrow}^m\}$ and (up
to sign) for the horizontally reflected spin networks $\{\overleftrightarrow\S^m\}$, that is
$$\sum_m\alpha_m{\S^\updownarrow}^m=0\qquad\text{and}\qquad\sum_m\alpha_m(-1)^{|\S^m_v|}\overleftrightarrow\S^m=0.$$
\end{thm}


\subsection{Basic Diagram Manipulations} 

In this section, we describe the \emph{spin network calculus}\index{spin network!calculus}, which governs diagram manipulations.

\begin{prop}
Any spin network can be expressed as a sum of diagrams with no crossings or loops. In particular,
\begin{equation}\index{fundamental binor identity}\index{binor identity}
  \binor;\qquad\circs\:\S=\tr(I)\S=2\S.
\end{equation}%
\end{prop}

The proof is given in \cite{Pet}. The first of these relations is called the \emph{Fundamental Binor Identity}, and represents a
fundamental type of structure in mathematics; it is the core concept in defining both the \emph{Kauffman Bracket Skein Module} 
\index{Kauffman Bracket Skein Module} in knot theory \cite{BFK} and the \emph{Poisson bracket} \index{Poisson bracket} on the set of loops on a surface, which Goldman describes in \cite{Gol1}. It can also be identified with the \emph{characteristic equation} \index{characteristic equation} for $2\times 2$ matrices \cite{Pet,Sik}.

Since $2\times 2$ matrices act on $V$, the definition of spin networks may be extended to allow matrices to act on diagrams: $\ims{\bfx}$ is the action $v\mapsto \bfx\cdot v$. The corresponding action on the tensor product $\vprod n$ is represented by%
$$\mupn{n}{\bfx}(v_1\tensor\cdots\tensor v_n)=\bfx v_1\tensor\cdots\bfx v_n.$$

The matrices $\bfx\in\sltc$ of interest in this chapter satisfy the following special property:
\begin{prop}\index{spin network!equivariance}
The spin network component maps $\pup,\cupp=\cupcu,$ and $\capp=\capcu$, and therefore all spin networks, are
equivariant under the natural action of $\sltc$ on $V$ described above.
\begin{proof}
The case for the identity $\pup$ is clear, while%
\begin{multline*}
  \capmm{\bfx}{\bfx}(v\tensor w)=\det[\bfx v\;\bfx w]=\det(\bfx\cdot[v\;w])\\
    =\det(\bfx)\cdot\det[v\;w]=1\cdot\det[v\;w]=\capp(v\tensor w)
\end{multline*}
shows that $\capcu\circ\bfx=\capcu=\bfx\circ\capcu$.

The proof for $\cupcu$ follows by reflecting this relation.
\end{proof}
\end{prop}

This means that matrices in such a diagram can ``slide across'' a vertex (local extremum) by simply
inverting the matrix, so that
$$\text{if}\quad\ims{\bar\bfx}=\bfx^{-1}\in\sltc,\quad\text{then}\quad\cupml{\bfx}=\cupmr{\bar\bfx}.$$
For a general matrix $\bfx\in\mtt$, the determinant is introduced in such relations since
$\cupmm{\bfx}{\bfx}=\det\left(\ims{\bfx}\right)\cupp$. If $\bfx$ is invertible, this implies
$$\cupml{\bfx}=\det\left(\ims{\bfx}\right)\cupmr{\bar\bfx}.$$

A \emph{closed} spin network with one or more matrices is called a \emph{trace diagram},\index{trace diagram} and may be
identified with a map $G\times\cdots\times G\to\C$. One of the primary motivations for this chapter
is the study of invariance properties of such maps. The simplest cases are given by:%
\begin{prop} For $\bfx\in\mtt$ and $\Id=\tmxt{1&0}{0&1}$,
\begin{equation}
  \circb = 2 = \tr(\Id); \qquad \circml{\bfx}=\tr(\bfx); \qquad \circmm{\bfx}{\bfx}=\det(\bfx)\cdot\tr(\Id).
\end{equation}
\end{prop}


\subsection{Symmetrizers and Irreducible Representations} %
Another important $\sltc$-equivariant map is the symmetrizer, defined by:
\begin{defn}
The \emph{symmetrizer} \index{symmetrizer} $\suppn{n}:\vprod n\to\vprod n$ is the map taking%
\begin{equation} \label{symeq}
v_1\tensor v_2\dtensor v_n\mapsto\frac1{n!}\sum_{\sigma\in\Sigma_n}
v_{\sigma(1)}\tensor v_{\sigma(2)}\dtensor v_{\sigma(n)},%
\end{equation}
where $v_i\in V$ and $\Sigma_n$ is the group of permutations on $n$ elements.
\end{defn}
For example,
\begin{equation*}
  \supp{2}=\tfrac12\left(\pid+\px\right)=\pid-\tfrac12\left(\pcc\right);
\end{equation*}
\begin{align*}
  \supp{3}&=\tfrac16
  	\left(\perm{1/1,2/2,3/3}+\perm{1/2,2/1,3/3}+\perm{1/3,2/2,3/1}+\perm{1/1,2/3,3/2}+\perm{1/2,2/3,3/1}+\perm{1/3,2/1,3/2}\right)\\
  &=\drawtl{}{}{1/1,2/2,3/3}-\tfrac23\left(\drawtl{1/2}{1/2}{3/3}+\drawtl{2/3}{2/3}{1/1}\right)
  	+\tfrac13\left(\drawtl{1/2}{2/3}{3/1}+\drawtl{2/3}{1/2}{1/3}\right)%
\end{align*}
Note that the crossings are removed by applying the Fundamental Binor Identity.

The defining equation \eqref{symeq} of $\suppn{n}$ should look familiar: its image is a subspace of
$\vprod n$ isomorphic to the \nth symmetric power $\Sym^nV$, and thus it can be thought of as
either the projection $\pi:\vprod n\to\Sym^nV$ or as the inclusion $i:\Sym^nV\to\vprod n$ (see
\cite{FH}, page 473).

What does this mean for us? If a diagram from $\vprod{k_i}$ to $\vprod{k_o}$ has symmetrizers at its top and bottom, it can be
thought of as a map between $V_{k_i}$ and $V_{k_o}$. We freely interpret such spin networks as maps between tensor powers of
these irreducible $\sltc$-representations.

\begin{prop}[Basic Symmetrizer Properties]
\begin{align}
\label{invariance} \text{{Invariance: }}       	&\quad \srinva{n}{A}=\srinvb{n}{A};\index{symmetrizer!invariance}\\%
\label{stacking} \text{{stacking relation: }}  	&\quad \srstack{n}{k};\index{symmetrizer!stacking}\\%
\label{capping} \text{{capping/cupping: }} 		 	&\quad \srcapping{n};\\%
\label{rotation} \text{{symmetrizer sliding: }}	&\quad \srextrema{n};%
\end{align}
\begin{proof}
The first relation \eqref{invariance} is evident if one expands the symmetrizer in terms of permutations, since permutations are
$\sltc$-equivariant.

The \emph{stacking relation} is the statement that symmetrizing the last $k$ elements of a
symmetric tensor has no effect, since they are already symmetric.

For the \emph{capping} and \emph{cupping relations}, notice that
$$\capp\circ\supp{2}(v\tensor w)=\capp(\tfrac12(v\tensor w+w\tensor v))=\tfrac12(\det[v\:w]+\det[w\:v])=0.$$
This implies the general case because, by the stacking relation, one may insert $\!\supp{2}\!$
between $\capp$ and $\suppn{n}$. The other case is similar.

There are a number of ways to demonstrate \eqref{rotation}. It follows by reflection (Proposition \ref{reflprop}) or as a special
case of $\sltc$-equivariance, since $\cupcd=\cupcl=\cupml{\bfx}$ for $\ims{\bfx}=g=\tmxt{0&1}{-1&0}\in\sltc$. More
directly, expand the symmetrizer into a sum of permutations. Since each permutation is a product of transpositions, then
\eqref{rotation} follows from the simple relation $\cupsxlr$. See \cite{Pet} for more details.
\end{proof}
\end{prop}

We now move on to some more involved relations among symmetrizers. Although it is easy to write
down an arbitrary $\suppn{n}$ in terms of permutations, it is usually rather difficult to write it
down in terms of diagrams without crossings (the Temperley-Lieb algebra).\index{Temperley-Lieb algebra}
The next two propositions
describe how to do exactly this. As such, they are a fundamental step in the proof of Theorem \ref{xmult},
which permits a fast computation of rank two central functions. 
\begin{prop}
The symmetrizer $\suppn{n}$ satisfies:
\begin{multline}\label{rec1}\index{symmetrizer!recurrence}
\srrecaa{n}=\srrecab{n}-\pfrac{n-1}{n}\srrecac{n}+\pfrac{n-2}{n}\srrecad{n}+\cdots\\%
                    + (-1)^i\pfrac{n-i}{n}\srrecae{n}+\cdots+(-1)^{n-1}\pfrac{1}{n}\srrecaf{n}.%
\end{multline}
\begin{proof}
If $\Sigma_n$ is the group of permutations on the set $N_n=\{1,2,\ldots,n\}$, then
$$|\Sigma_n|=|N_n|\:|\Sigma_{n-1}|.$$
Interpret $|\Sigma_n|$ as the number of ways to arrange $n$ people in a line. To do this, one may first select someone to be at
the front of the line ($|N_n|$ choices), and then rearrange the remaining $n-1$ people ($|\Sigma_{n-1}|$ choices).

In diagram form, the selection of someone to head the line corresponds to one of the diagrams
$$\trarrowa, \trarrowb, \trarrowc, \ldots, \trarrowd, \ldots, \trarrowe.$$
The arrangement of the remaining people corresponds to $\srrecab{n}$. Thus, the diagrammatic form of the above
interpretation is:
$$\srrecaa{n}=\tfrac1n\srrecab{n}\circ\left(\trarrowa+\trarrowb+\trarrowc+\cdots+\trarrowd+\cdots+\trarrowe\right).$$

Now, use the binor identity to remove crossings. Most of the resulting terms disappear, since any
term whose cups are not in the `first position' on top will vanish due to the \emph{capping
relation}. In particular:
$$\srrecab{n}\!\!\!\circ\trarrowd=\srrecab{n}\!\!\!-\:\srrecac{n}\!\!\!+\:\srrecad{n}\!\!\!+\cdots+(-1)^i\:\srrecae{n},$$ where
$i$ is the number of `kinks' \slantn{} in \trarrowd or 1 plus the number of kinks in
\srrecae{n}. Finally, group the number of terms on the righthand side with the same number of kinks
together: there will be $n-i-1$ terms with $i$ kinks.
\end{proof}
\end{prop}


\begin{prop} $\suppn{n}$ also satisfies the recurrence relations:
\begin{align} \label{rec2}\index{symmetrizer!recurrence}
&\srrecb{n}{i};\\ %
\label{rec3}
&\srrecca{n}=\srreccb{n}-\pfrac{n-1}{n}\srreccc{n}. %
\end{align}
\begin{proof}
Compose relation \eqref{rec1} with $\suppnbig{i}\tensor\suppnbig{n-i}$. This has no effect on the lefthand
side, by the \emph{stacking relation}. On the righthand side, all but one of the terms with a cap
on the bottom vanish, due to the \emph{capping relation}, since they will cap off either the
$\suppnbig{i}$ or the $\suppnbig{n-i}$. The one term which remains `caps between' these two symmetrizers.
The coefficient is $(-1)^i\pfrac{n-i}{n}$ since in recurrence \eqref{rec1}, $i$ is equal to one
more than the number of kinks \slantn{} in \srrecae{n}.

Relation \eqref{rec3} is a special case of \eqref{rec2} for $i=1$.
\end{proof}
\end{prop}

The next relations follow directly from these recurrences:
\begin{prop}[Looping Relations]\index{spin network!looping relation}
\begin{align} \label{rec4}
\srloopa{n} &= \pfrac{n+1}{n} \suppnbig{n-1}. \\%
\text{When $k$ strands of $\suppn{n}$ are}&\text{ closed off:}\hspace{2in}\nonumber \\%
\label{rec5}%
\srloopc{n}{k} &= \pfrac{n+1}{n-k+1} \suppnbig{n-k}. \\%
\label{deltaval}%
\circn{n} &= n+1.
\end{align}

\begin{proof}
Close off the left strand in \eqref{rec3} above. Then, $\srrecca{n},\srreccb{n},$ and $\srreccc{n}$ become $\srloopa{n}$, $\srloopb{n-1}=2\:\suppnbig{n-1}$ and $\suppnbig{n-1}$, respectively. Now collect terms
to get \eqref{rec4}, and proceed to \eqref{rec5} by applying the first relation $k$ times. Finally, \eqref{deltaval} is a special
case of \eqref{rec5} with $k=n$.
\end{proof}
\end{prop}


\subsection{Symmetrizers and Trivalent Spin Networks}%
Recall the Clebsch-Gordan decomposition (Proposition \ref{cgform}):\index{Clebsch-Gordan formula}
$$V_a \tensor V_b \isom \bigoplus_{\Adm abc} V_c, \qquad \lceil a,b \rfloor = \{a+b,a+b-2,\ldots,|a-b|\}.$$%
The requirement $\Adm abc$ is equivalent to the following symmetric condition:

\begin{defn} \label{admissible}\index{admissible triple}
A triple $(a,b,c)$ of nonnegative integers is \emph{admissible}, and we write $\Adm abc$, if
\begin{equation}\label{admcond}
  \tfrac12(-a+b+c),\quad\tfrac12(a-b+c),\quad\tfrac12(a+b-c)\in\N.
\end{equation}
\end{defn}

Two maps arise from the Clebsch-Gordon decomposition: an injection $\iota^{a,b}_c:V_c\to V_a\tensor
V_b$ and a projection $(\iota^*)_{a,b}^c:V_a\tensor V_b\to V_c$. Both have simple diagrammatic
depictions \cite{CFS}: %
$$\iota^{a,b}_c=\!\!\!\!\cfupbig{a}{b}{c} \hspace{-10pt}: V_c\to V_a\tensor V_b; \qquad
  (\iota^*)_{a,b}^c=\!\!\!\!\cfdownbig{a}{b}{c} \hspace{-5pt}: V_a\tensor V_b\to V_c.$$
The admissibility condition \eqref{admcond} is the requirement that there is a nonnegative number
of strands connecting each pair of symmetrizers. These ``strand numbers'' appear frequently in
diagram manipulations, and will be referenced by the Greek letters $\alpha,\beta,\gamma$:

\begin{conv}\label{alphabeta}
Given an admissible triple $(a,b,c)$, denote by $\alpha$, $\beta$, and $\gamma$ the total number of
strands connecting $V_b$ to $V_c$, $V_a$ to $V_c$, and $V_a$ to $V_b$, respectively. Also, denote
by $\delta$ the total number of strands in the diagram. Then:%
$$\alpha=\tfrac12(-a+b+c),\quad\beta=\tfrac12(a-b+c),\quad\gamma=\tfrac12(a+b-c); \quad\delta=\tfrac12(a+b+c).$$
Note that $(a,b,c)$ is admissible if and only if $\alpha,\beta,\gamma\in\N$.
\end{conv}

\begin{conv}
Because the maps $\iota^{a,b}_c$ and $(\iota^*)_{a,b}^c$ will be so important for the remainder of this chapter, we introduce a notation which
simplifies their depiction. Let $n$ lines with a symmetrizer be represented by one {\bf thick} line labelled $n$, so that
$\tup(n)\equiv\suppn{n}$.
\end{conv}
\begin{defn}
A \emph{trivalent spin network}\index{trivalent spin network}\index{spin network!trivalent} $\S$ is a graph drawn on the plane with vertices of degree $\leq 3$
and edges labelled by positive integers such that:
\begin{itemize}
\item 2-vertices are ciliated and coincide with local extrema;
\item 3-vertices are drawn `up' \tupy(,,) or `down' \tdowny(,,);
\item any two edges meeting at a 2-vertex have the same label;
\item the three labels adjacent to any vertex form an admissible triple.
\end{itemize}
If there are $m$ input edges with labels $l_i$ for $i=1,\ldots,m$ and $n$ output edges with labels
$l'_i$ for $i=1,\ldots,n$, the network is identified with a map between tensor products of
irreducible $\sltc$ representations,%
  $$f_{\S}:V_{l_1}\dtensor V_{l_m}\to V_{l'_1}\dtensor V_{l'_n}.$$
This map is computed by identifying $\S$ with a regular spin network using the following identifications:
\begin{gather*}
  \tup(n)\equiv\suppnbig{n}\qquad
  \tupc(n)\equiv\upcln{n}\qquad
  \tcup(n)=\tcupc(n)\equiv\cupcun{n}\\
  \tupy(a,b,c)\equiv\cfupbig{a}{b}{c}\qquad\tdowny(a,b,c)\equiv\cfdownbig{a}{b}{c}.
\end{gather*}
\end{defn}
Note that ciliations are normally chosen to be on the local extrema, and degree-3 vertices, when
expanded, also have a number of ciliated vertices. The need to keep track of these ciliations makes
diagram manipulation a more delicate operation.


\subsection{Trivalent Diagram Manipulations}\label{spinprops}
This section describes in detail the relations which may be used to manipulate trivalent spin networks.
For the remainder of this chapter, we assume that all sets of labels incident to a common vertex in a
diagram are admissible. Moreover, whenever we sum over a label in a diagram, the sum is taken over
all possible values of that label which make the requisite triples in the diagram admissible.


Any closed trivalent spin network may be interpreted as a constant. The simplest such diagrams are given by

\begin{prop}
Let $\Theta(a,b,c)=\tbubbled(a,b,c)$ and $\Delta(c)=\tcirc(c)$. Then $\Theta(a,b,c)$ is symmetric in $\{a,b,c\}$ and
explicitly (recall the $\alpha,\beta,\gamma,\delta$ given in Convention \ref{alphabeta}):
\begin{align}
\label{delta} \Delta(c) &=c+1=\dim(V_c);\\
\label{theta1} \Theta(a,b,c) &=
\tfrac{\pfrac{-a+b+c}{2}!\pfrac{a-b+c}{2}!\pfrac{a+b-c}{2}!\pfrac{a+b+c+2}{2}!}{a!b!c!}
                              = \tfrac{\alpha!\beta!\gamma!(\delta+1)!}{a!b!c!}; \\
\label{theta2} \Theta(1,a,a+1) &= \Delta(a+1) = a+2.
\end{align}
\begin{proof}
The first equation \eqref{delta} is a consequence of the \emph{Looping Relation} \eqref{rec4}. That
$\Theta(1,a,a+1)=\Delta(a+1)$ is a consequence of the stacking relation, and demonstrates
$\eqref{theta2}$. We refer the reader to \cite{CFS} for the $\Theta(a,b,c)$ formula.
\end{proof}
\end{prop}

Ratios of $\Delta$ and $\Theta$ show up in the next two propositions, which tell us how to ``pop
bubbles'' and how to ``fuse together'' two thick edges. The first demonstrates the usefulness of
Schur's Lemma \index{Schur's Lemma} (Proposition \ref{schurs}) in diagrammatic techniques.

\begin{prop}[Bubble Identity]\label{bubble}\index{bubble identity}
$\tbubble(a,b,c,d) = \left(\tfrac{\Theta(a,b,c)}{\Delta(c)}\:\tup(c)\right)\delta_{cd}$, where
$\delta_{cd}$ is the Kronecker delta.
\begin{proof}
Schur's Lemma requires $\tbubble(a,b,c,d)=C\:\tup(c)\delta_{cd}$ for some constant $C$, since
$\tbubble(a,b,c,d)$ is a map between irreducible representations.
This equation remains true if we ``close off'' the diagrams, giving:%
$$\tbubbled(a,b,c)=C\tcirc(c) \qimplies C=\frac{\Theta(a,b,c)}{\Delta(c)}.\mbox{\qedhere}$$%
\end{proof}
\end{prop}

\begin{prop}[Fusion Identities]\label{fusion}\index{fusion identity}
\begin{align*}
  \tfancyup(a,b)&=\sum_{\Adm abc}\pfrac{\Delta(c)}{\Theta(a,b,c)}\tfused(a,b,b,a,c)\\
  \tx(a,b)&=\sum_{\Adm abc}(-1)^{\half(a-b+c)}\pfrac{\Delta(c)}{\Theta(a,b,c)}\thright(a,b,a,b,c).
\end{align*}
\begin{proof}
Maps of the form $\tfused(a,b,b,a,c)$ for $\Adm abc$ form a basis for the space of $\sltc$-equivariant maps $V_a\tensor V_b\to V_a\tensor
V_b$ \cite{CFS}. Thus, we may express
the first diagram as a linear combination:%
  $$\tfancyup(a,b) = \sum_{\Adm abc} C(c) \tfused(a,b,b,a,c).$$%
For a fixed $d\in\lceil{a,b}\rfloor$, the constant $C(d)$ is computed by composing this expression with $\tupy(a,b,d)$:
\begin{align*}
	\tupy(a,b,d) &= \sum_{\Adm abc} C(c) \tupy(a,b,c) \circ \tbubble(a,b,c,d)
           = \sum_{\Adm abc} C(c) \pfrac{\Theta(a,b,c)}{\Delta(c)} \tupy(a,b,c) \circ \tup(d)\delta_{cd}\\%
           &= C(d) \pfrac{\Theta(a,b,d)}{\Delta(d)}\tupy(a,b,d)
           \qimplies C(d)=\frac{\Delta(d)}{\Theta(a,b,d)}.
\end{align*}
The second equation follows from the first and from Proposition \ref{funsigns} below:
\begin{align*}
  \tx(a,b)=(-1)^b\txkink(a,b)
    &=\sum_{\Adm abc}(-1)^b\pfrac{\Delta(c)}{\Theta(a,b,c)}\thtwistedx(a,b,a,b,c)\\
    &=\sum_{\Adm abc}(-1)^{\half(a-b-c)}\pfrac{\Delta(c)}{\Theta(a,b,c)}\thtwisted(a,b,a,b,c)\\
    &=\sum_{\Adm abc}(-1)^{\half(a-b+c)}\pfrac{\Delta(c)}{\Theta(a,b,c)}\thright(a,b,a,b,c).\qedhere
\end{align*}

\end{proof}
\end{prop}


The identity $\upcll$ gives rise to the following compendium of sign changes through
diagram manipulations:
\begin{prop}\label{funsigns}
\begin{align}
  \label{sign1} \tkink(n)&=(-1)^n\:\tup(n);\\
  \label{sign2} \tdownyx(a,b,c)&=(-1)^{\half(a+b-c)}\tdowny(a,b,c);\\
  \label{sign3} \tupykink(a,b,c)&=(-1)^{\half(-a+b+c)}\tdowny(a,b,c);\\ %
  \label{sign4} \thtwisted(a,b,c,d,e)&=(-1)^{\half(a+b+c+d-2e)}\thleft(a,b,c,d,e);\\ %
  \label{sign5} (-1)^{\half(a+c)}\thleft(a,b,c,d,e)&=(-1)^{\half(b+d)}\thright(a,b,c,d,e);\\ %
  \label{sign6} \thtwisted(a,b,c,d,e)&=(-1)^{b+d-e}\thright(a,b,c,d,e).%
\end{align}
\begin{proof}
First, \eqref{sign1} is just a restatement of $\kinkn{n}=(-1)^n\slantn{n}$, and \eqref{sign2} follows directly from
the Proposition \ref{reflprop}, since $\tdowny(a,b,c)$ contains $\gamma=\half(a+b-c)$ local extrema and
$\tdownyx(a,b,c)=\tdowny(b,a,c)$.

For \eqref{sign3}, notice that in the simplest case $$\triplel,$$ the negative
sign comes from the strand on top of the diagram. Similarly, the general case for transforming
$\tupykink(a,b,c)$ into $\tdowny(a,b,c)$ has a sign for each strand
between $b$ and $c$, giving $(-1)^\alpha=(-1)^{\half(-a+b+c)}$. This identity is used twice to give
\eqref{sign4}.

Finally, \eqref{sign5} follows from:
$$\thleft(a,b,c,d,e)=(-1)^e\thkink(a,b,c,d,e)
  =(-1)^{e+\half(d+e-a+b+e-c)}\thright(a,b,c,d,e),$$
and \eqref{sign6} is given by combining \eqref{sign4} and \eqref{sign5}.
\end{proof}
\end{prop}

The above relations permit the definition of a ``$\frac\pi4$-reflection'' on certain types of
diagrams, which will be important later:
\begin{prop}\label{funsigns2}\index{spin network!reflection}
If a relation consists entirely of terms of the form $\thright(a,b,c,d,e)$ and $\tfused(a,b,c,d,f)$, then one may
``reflect about the line through $a$ and $c$'' in the following sense:%
\begin{equation*}
  \sum_e\alpha_e\thright(a,b,c,d,e)=\sum_f\beta_f\tfused(a,b,c,d,f)\iff\sum_e\alpha_e\tfused(a,d,c,b,e)=\sum_f\beta_f\thright(a,d,c,b,f).
\end{equation*}
\begin{proof}
By horizontally reflecting the first relation, using Theorem \ref{reflector},
\begin{align*}
  \sum_e\alpha_e\thright(a,b,c,d,e)&=\sum_f\beta_f\tfused(a,b,c,d,f)\\
    \iff\sum_e\alpha_e(-1)^{\half(a+b+c+d-2e)}\thleft(b,a,d,c,e)&=\sum_f\beta_f(-1)^{\half(a+b+c+d-2f)}\tfused(b,a,d,c,f)\\
    \iff\sum_e\alpha_e\thleft(b,a,d,c,e)&=\sum_f\beta_f\tfused(b,a,d,c,f),
\end{align*}
where the signs cancel due to the admissibility conditions.

Now, add strands to both sides, so that the right side $\tfused(b,a,d,c,f)$ becomes
$$\thtwisted(a,d,c,b,f)=(-1)^{b+d-f}\thright(a,d,c,b,f).$$
Likewise, on the left side, $\thleft(b,a,d,c,e)$ becomes $(-1)^{b+d-e}\tfused(a,d,c,b,e)$. Once again,
admissibility implies that $e$ and $f$ must have the same parity, so these signs cancel.
\end{proof}
\end{prop}

Two alternate versions of this proposition follow (see \cite{Pet}).

\begin{cor}\label{funsigns3}\index{spin network!reflection}
\begin{equation*}
  \sum_e\alpha_e\thleft(a,b,c,d,e)=\sum_f\beta_f\tfused(a,b,c,d,f)\iff\sum_e\alpha_e\tfused(a,d,c,b,e)=\sum_f\beta_f\thleft(a,d,c,b,f)%
\end{equation*}
\begin{multline*}
  \quad\sum_e\alpha_e\thright(a,b,c,d,e)=\sum_f\beta_f\tfused(a,b,c,d,f) \\
    \iff\sum_e\alpha_e(-1)^{\half(e-b)}\tuptreer(a,b,c,d,e)=\sum_f\beta_f(-1)^{\half(d-f)}\tuptreel(a,b,c,d,f).
\end{multline*}
\end{cor}

%% file: spinsl2-decomposition.tex

\section{Decomposition of $\C[G]$}\label{decomp}

The following theorem is a consequence of the ``unitary trick''\cite{Dol},\index{unitary trick} the Peter-Weyl Theorem,\index{Peter-Weyl Theorem} and the fact that the set of
matrix coefficients of $G$ is exactly its coordinate ring \cite{CSM}. We offer a self-contained constructive proof in Section
\ref{proofdecomp}, since it gives an explicit correspondence between regular functions and spin networks.

\begin{thm} \label{pwdecomp} \index{coordinate ring!decomposition}
There is a $G$-module isomorphism $$\cb[G]\isom\sum_{n\geq 0}V_n^*\otimes V_n.$$
\end{thm}


\subsection{Central Functions} \label{cfs}
Theorem \ref{pwdecomp} allows $\C[G\times G]^G$ to be described in terms of an additive basis of
class functions that have an elegant realization as spin networks.  Indeed, together with the
Clebsch-Gordan decomposition, it implies
\begin{eqnarray*}
\cb[G\times G] & \isom& \cb[G]\otimes \cb[G]\\
 & \isom& \left(\sum_{a\geq 0} V_a^*\otimes V_a\right) \otimes \left(\sum_{b\geq 0}V_b^*\otimes V_b\right)\\
 & \isom& \sum_{a\geq 0} \sum_{b\geq 0} V_a^*\otimes V_a\otimes V_b^* \otimes V_b\\
 & \isom& \sum_{0\leq a, b < \infty} \! \! \left(V_a^* \otimes V^*_b\right) \otimes \left(V_a \otimes
             V_b\right)\\
 & \isom& \sum_{0\leq a, b < \infty}  \! \! \! \left(\sum_{i=0}^{\mathrm{min}(a,b)}V^*_{a+b-2i}\right)
             \otimes \left(\sum_{j=0}^{\mathrm{min}(a,b)}V_{a+b-2j}\right)\\
 & \isom& \sum_{\substack{0\leq a,b<\infty\\ 0\leq i,j \leq \mathrm{min}(a,b)}} V^*_{a+b-2i} \otimes V_{a+b-2j}.
\end{eqnarray*}
Since the above maps are $G$-equivariant,
\begin{equation}\label{cgg1}
  \cb[G\times G]^G \isom \sum_{\substack{0\leq a,b<\infty\\ 0\leq i,j \leq \mathrm{min}(a,b)}} \left(V^*_{a+b-2i} \otimes V_{a+b-2j}\right)^G.
\end{equation}
By Schur's Lemma (Proposition \ref{schurs}),\index{Schur's Lemma}
$$\mathrm{dim}_{\cb}\left(V^*_{a+b-2i} \otimes V_{a+b-2j}\right)^G= \left\{
  \begin{array}{ll} 1 & \textrm{if $i=j$}\\ 0 & \textrm{if $i\not=j$}\end{array}\right.,$$
so
$$\cb[G\times G]^G \isom\sum_{\substack{0\leq a,b<\infty\\ 0\leq j \leq \mathrm{min}(a,b)}}
              \mathrm{End}(V_{a+b-2j})^G.$$
\begin{defn}
Given the above isomorphism, for each $\Adm abc$ (see Definition \ref{admissible}),  there exists a
class function $\ch abc\in\C[G\times G]^G$ which corresponds to a generating homothety (unique up
to scalar) in $\mathrm{End}(V_c)^G$. We refer to the functions $\ch abc$ as \emph{central
functions}.\index{central function}
\end{defn}

Denote by $\cspan\ch abc\subset \cb[G\times G]^G$ the linear span over $\cb$ of $\ch abc$. Then \eqref{cgg1} may be rewritten as

$$\cb[G\times G]^G \isom \sum_{\substack{0\leq a,b<\infty\\ \Adm abc}}
              \cspan\ch abc.$$
Thus, the central functions $\ch abc$ form an additive basis for the ring of regular functions on $\XX=\mathrm{Spec}_{max}(\cb[\mathcal{R}]^G)=\mathcal{R}/\!/G$. In Section \ref{rank2},
we describe the multiplicative structure of $\cb[G\times G]^G$ in terms of this basis.

The central functions may be described using the Clebsch-Gordan injection $\iota_c^{a,b}:V_c\hookrightarrow V_a\tensor V_b$: 
	$$\ch{a}{b}{c}(\xb_1,\xb_2)=
		\mathrm{tr}\bigg(\iota({\sf c}_i^*)\Big((\xb_1,\xb_2)\cdot\iota({\sf c}_j)\Big)\bigg)_{\!\!ij},$$
where $\{{\sf c}_j\}$ is a basis for $V_c$. We will omit indices on $\iota$ when they are clear from context.

The functions $\ch abc$ take a natural diagrammatic form. If the matrix $\xb$ is represented
diagrammatically by $\ims{\bfx}:V\to V$, then its action on $V_a$ can be represented by
$\tmup(a)(\bfx)\equiv\srinva{a}{\bfx}.$ A closed spin network with $r$ different matrices is
an invariant regular function $G^{\times r}\to\C$. In particular, since $\tupy(,,)$ and $\tdowny(,,)$
are the Clebsch-Gordan injection and projection, respectively,
	$$\ch abc(\xb_1,\xb_2)=\cftwo(\bfx_1,\bfx_2)(a,b,c)=\cftwob(\bfx_1,\bfx_2)(a,b,c).$$
\index{central function!diagrammatic}
As a special case, setting $\xb_1=\xb_2=\id$, where $\id$ is the identity matrix in $G$, gives $\ch
abc (\id,\id)=\Theta(a,b,c)$.



\subsection{Proof of $\C[G]$ Decomposition Theorem} \label{proofdecomp}
Define $$\Upsilon :\sum_{n\ge0} V_n^*\otimes V_n\longrightarrow \cb [G]$$ by linear extension of
the mapping
	$${\sf n}^*_{n-k}\otimes {\sf n}_{n-l}\mapsto {\sf n}^*_{n-k}(\xb\cdot {\sf n}_{n-l}),$$
where $\xb=\imx{x}$ is a matrix variable.

\begin{prop}
$\Upsilon$ is a well-defined $G$-equivariant morphism.
\begin{proof}
The image of $\Upsilon$ consists of regular functions since
\begin{align*}
{\sf n}^*_{n-k}(\xb\cdot {\sf n}_{n-l})
&=  {\sf n}^*_{n-k}\left((x_{11} e_1+x_{21} e_2)^{n-l}(x_{12} e_1 + x_{22} e_2)^l\right)\\
&=\sum_{\substack{i+j=k\\0 \le i \le n-l \\ 0 \le j \le l}}
\tbinom{n}{k}^{-1}\tbinom{n-l}{i}\tbinom{l}{j}x_{11}^{n-l-i}x_{12}^{l-j}x_{21}^ix_{22}^j.
\end{align*}
Equivariance is verified by the calculation:
\begin{align*}
  \Upsilon(g\cdot({\sf n}^*_{n-k}\otimes {\sf n}_{n-l}))&=\Upsilon\left((g\cdot {\sf n}^*_{n-k})\otimes (g\cdot {\sf n}_{n-l})\right)\\ %
    &=(g\cdot {\sf n}^*_{n-k})(\xb\cdot (g\cdot {\sf n}_{n-l})) ={\sf n}^*_{n-k}((g^{-1}\xb g)\cdot{\sf n}_{n-l})\\ %
    &= g\cdot {\sf n}^*_{n-k}(\xb\cdot {\sf n}_{n-l}) = g\cdot \Upsilon({\sf n}^*_{n-k}\otimes {\sf n}_{n-l}).\qedhere%
\end{align*}
\end{proof}
\end{prop}

There is a right action of $G$ on $\cb[G]$ given by $f\cdot g(\xb)=f(\xb g).$ Denote by $\cgright$
the ring $\cb[G]$ with this right action, to distinguish it from the conjugation action already
imposed on $\cb[G]$. Additionally, $G$ acts on the left of $\mathrm{Hom}_G(V_n, \cgright)$ by
	$$(g\cdot \gamma)(v)(\xb)=\gamma_v(g^{-1}\xb),$$
where $\gamma_v=\gamma(v).$ This action is well-defined since
  $$(g\cdot \gamma)(g'\cdot v)(\xb) =\gamma_{g'\cdot v}(g^{-1}\xb) =\gamma_v(g^{-1}\xb g')
    =\big((g\cdot\gamma)(v)\big)\cdot g'(\xb).$$

The next two lemmas, whose proofs are deferred, define two additional maps which will be used to prove the theorem.
\begin{lem} \label{pwl1}
The map%
$$\Phi: \sum_{n\geq 0}\mathrm{Hom}_G(V_n, \cgright)\otimes V_n\longrightarrow \cb [G]$$%
defined by linearly extending the mappings $\gamma \otimes v \mapsto \gamma(v)$ is an isomorphism
of $G$-modules.
\end{lem}

\begin{lem} \label{pwl2}
Define $\Psi_n: V_n^*\to\mathrm{Hom}_G(V_n,\cgright)$
by $w^*\mapsto \mathbf{F}_{w^*}$, where $\mathbf{F}_{w^*}(v)(\xb)=w^*(\xb\cdot v)$.
Then the map%
	$$\Psi:\sum_{n\ge0} V_n^*\tensor V_n\longrightarrow\sum_{n\ge0} \Hom_G(V_n,\cgright)\tensor V_n$$
given by $\Psi=\sum(\Psi_n\otimes \mathrm{id})$ is an isomorphism of $G$-modules.
\end{lem}

Assuming the above lemmas, Theorem \ref{pwdecomp} is equivalent to showing that the following diagram commutes:
$$\xymatrix{
  {\displaystyle\sum_{n\ge0}} V_n^*\tensor V_n \ar[rr]^{\Upsilon}\ar[dr]^{\Psi}&&\C[G]\\
  &\hspace{-.5in}{\displaystyle\sum_{n\ge0}} \Hom_G(V_n,\cgright)\tensor V_n \ar[ur]^{\Phi}\hspace{-.5in}&
}$$
The proof of commutativity follows:
	$$\Phi\circ\Psi(w^*\otimes v)=\Phi(\mathbf{F}_{w^*}\otimes v)=\mathbf{F}_{w^*}(v)
    =w^*(\xb\cdot v) =\Upsilon(w^*\otimes v).\qed$$

It remains to establish Lemmas \ref{pwl1} and \ref{pwl2}. The proof of Lemma \ref{pwl1} requires some preliminary technical results.

\begin{lem} \label{pwl3}
Every regular function is contained in a finite-dimensional sub-representation of $\cb[G]$.
\end{lem}
\begin{proof}[Proof of Lemma \ref{pwl3}]
The following $G\times G$-action encompasses both the right and diagonal $G$-actions defined above. Let $$\alpha:G\times G\times
G\longrightarrow G$$ be defined by $(g_1,g_2,\xb)\mapsto g_1\xb g_2^{-1}$, and further let
\begin{equation}\label{pullback}
  \alpha^*:\cb[G]\longrightarrow \cb[G\times G\times G]\isom\cb[G]^{\otimes 3}
\end{equation}
be defined by $f\mapsto f\circ \alpha$, the {\it pull-back} of regular functions on $G$ to regular functions on $G\times G\times
G$. For $f\in \cb[G]$, \eqref{pullback} implies that there exist $n_f\in\N$ and regular functions $f_i,f_i',f_i''$ for $1\leq
i\leq n_f$ such that
$$\alpha^*(f)=\sum_{i=1}^{n_f}f_i\otimes f_i'\otimes f_i''.$$
Therefore
$$\alpha^*(f)(g_1^{-1},g_2^{-1},\xb)=\sum_{i=1}^{n_f}f_i(g_1^{-1})f_i'(g_2^{-1})f_i''(\xb).$$
On the other hand,
$$\alpha^*(f)(g_1^{-1},g_2^{-1},\xb)=f(\alpha(g_1^{-1},g_2^{-1},\xb))=f(g_1^{-1}\xb g_2)=((g_1,g_2)\cdot
f)(\xb),$$ which implies
\begin{equation}\label{orbitbasis}
  (g_1,g_2)\cdot f =\sum_{i=1}^{n_f}f_i(g_1^{-1})f_i'(g_2^{-1})f_i''.
\end{equation}
Let $(G\times G)f=\{(g_1,g_2)\cdot f: f\in G\}$ be the $G\times G$-orbit of $f$, and let $W_f$ be
the linear subspace spanned over $\cb$ by $(G\times G)f$ in $\cb[G]$. By \eqref{orbitbasis},
$\{f_i''\}$ is a spanning set for $W_f$, and so $W_f$ is finite-dimensional. Clearly $W_f$ is
$G\times G$-invariant, and so invariant with respect to the diagonal and right $G$-actions. Thus,
it is a finite-dimensional sub-representation containing $f$.
\end{proof}

\begin{lem} \label{pwl4}
$\cb[G]$ is completely $G\times G$-reducible.
\end{lem}

\begin{proof}[Proof of Lemma \ref{pwl4}]
Let $\mathcal{I}$ be the set of direct sums of irreducible finite-dimensional sub-representations
of $\cb[G]$. $\mathcal{I}$ is partially ordered by set inclusion and is nonempty. Thus, by Zorn's
lemma there exists a maximal element $M\in\mathcal{I}$. If $M\neq\cb[G]$, then consider any
$f\notin M$.  By Lemma \ref{pwl3}, there exists a finite-dimensional sub-representation $W_f$ that
contains $f$. Let $K=\mathrm{SU}(2)$ be the maximal compact subgroup of $G$. Restrict the action of
$G\times G$ to $K\times K$ to find an invariant orthogonal complement to $W_f$ in $M\cup W_f$.
Denote this complement by $M^\perp$. Then $M^\perp\oplus W_f \in \mathcal{I}$, since $K\times K$
representations extend to $G\times G$ representations. Hence $M$ is not maximal, which is a
contradiction. Therefore $\cb[G]$ is completely reducible with respect to the $G\times G$-action,
and so%
  $$\cb[G]\isom\sum_{j\geq 0}c_j V_j,$$%
where $c_j$ is the (possibly infinite) multiplicity of $V_j$ in $\C[G]$. This decomposition also holds for $\C[G]$ with
both the right and diagonal actions since they are restrictions of the same $G\times G$-action.
\end{proof}

\begin{proof}[Proof of Lemma \ref{pwl1}]
By Lemma \ref{pwl4},
$$\Phi: \sum_{n\geq 0} \big(\mathrm{Hom}_G(V_n, \cgright)\otimes V_n\big) \longrightarrow \cb [G]$$
is an isomorphism if and only if
$$\sum_{n\geq 0}\left(\sum_{j\geq 0} \mathrm{Hom}_G(V_n, c_jV_j)\otimes V_n\right)\longrightarrow\sum_{j\geq 0}c_j V_j$$
is an isomorphism. By Schur's Lemma, this reduces to
$$\sum_{n\ge0}\left(c_n\C\otimes V_n\right)\isom\sum_{n\geq 0}\left(\mathrm{Hom}_G(V_n, c_nV_n)\otimes V_n\right)\longrightarrow\sum_{n\geq 0}c_n V_n.$$
However, this is the map sending $\sum\lambda\otimes v\mapsto \sum \lambda v$  for $\lambda \in
\cb$ and $v\in V_n$, which is canonically an isomorphism.
\end{proof}

The final task is to show that $\Psi$ is an isomorphism:

\begin{proof}[Proof of Lemma \ref{pwl2}]
Recall that $$\Psi_n: V_n^*\longrightarrow \mathrm{Hom}_G(V_n,\cgright)$$
was defined by $w^*\mapsto \mathbf{F}_{w^*}$, where $\mathbf{F}_{w^*}(v)(\xb)=w^*(\xb\cdot v)$.
$\Psi_n$ is well-defined since
  $$\mathbf{F}_{w^*}(g\cdot v)(\xb)=w^*(\xb\cdot(g\cdot v))=w^*((\xb g)\cdot v)=\mathbf{F}_{w^*}(v)(\xb g)=(\mathbf{F}_{w^*}(v))\cdot g(\xb),$$
and is $G$-equivariant because
\begin{align*}
\Psi_n(g\cdot w^*)(v)(\xb)&=\mathbf{F}_{g\cdot w^*}(v)(\xb)= (g\cdot w^*)(\xb\cdot v)=w^*((g^{-1}\xb)\cdot v)\\ %
&=\mathbf{F}_{w^*}(v)(g^{-1}\xb)= (g\cdot \mathbf{F}_{w^*})(v)(\xb)= g\cdot \Psi_n (w^*)(v)(\xb).
\end{align*}

Since $V^*_n$ is irreducible, Schur's Lemma implies $\Psi_n$ is injective. We now show
surjectivity. Consider $\gamma \in \mathrm{Hom}_G(V_n, \cgright).$  For $\mathbb{I}\in G$,
$\gamma(v)(\mathbb{I})$ is a linear functional on $V_n$. Hence there exists $w^* \in V_n^*$ such
that $w^*(v)= \gamma(v)(\mathbb{I})$ for all $v\in V_n$. The following computation establishes that
$\Psi_n(w^*)=\gamma$:
  $$\mathbf{F}_{w^*}(v)(\xb)=w^*(\xb\cdot v)=\gamma(\xb\cdot v)(\mathbb{I}) %
    =(\gamma(v))\cdot \xb(\mathbb{I})=\gamma(v)(\mathbb{I}\xb)=\gamma(v)(\xb).$$
Therefore $\Psi_n$ is an isomorphism and so is $\Psi=\sum(\Psi_n\otimes \mathrm{id})$:
	$$\sum_{n\geq 0} V_n^*\otimes V_n\isom\sum_{n\geq0}\left(\mathrm{Hom}_G(V_n, \cgright)\otimes V_n\right).\qedhere$$
\end{proof}


\subsection{Ring Structure of $\cb[G]^G$}
We have established $$\cb[G]\isom\sum_{n\geq 0} V_n^*\otimes V_n.$$  By Schur's Lemma and the fact
that $V_n^*\tensor V_n\isom\End(V_n)$,
$$\cb[G]^G \isom \sum_{n\geq 0} (V_n^*\otimes V_n)^G \isom \sum_{n\geq 0}
\cspan\chxx^n,$$ where $\chxx^n \in \mathrm{End}(V_n)^G$ is a multiple of the identity.

The isomorphism $\mathrm{End}(V_n)\to V_n^*\otimes V_n$ is given by%
$${\sf n}_{n-l}({\sf n}_{n-k})^T\mapsto\tbinom{n}{k}{\sf n}^*_{n-k}\otimes{\sf n}_{n-l}.$$ Therefore, the
central function $\chxx^n$ \index{central function!rank one} corresponds to an invariant function in $\C[G]^G$ by%
$$\chxx^n=\sum_{i=0}^n{\sf n}_i({\sf n}_i)^T
  \longmapsto\sum_{i=0}^n\tbinom{n}{i}{\sf n}^*_i\otimes{\sf n}_i
  \overset{\Upsilon}{\longmapsto}\sum_{i=0}^n\tbinom{n}{i}{\sf n}_i^*(\xb\cdot{\sf n}_i).$$
We will freely identify $\chxx^n$ with its image in $\C[G]^G$.


For example, the trivial representation $V_0$ gives $\chxx^0=1$. The standard representation $V_1$
has diagonal matrix coefficients $x_{11}$ and $x_{22}$, hence
$$\chxx^1= x_{11}+x_{22}=\tr(\xb).$$

The remaining functions may be computed directly, or by using the following product formula:

\begin{thm}[Product Formula]\index{central function!product}
\begin{eqnarray}\label{rank1prod}
  \chxx^a \chxx^b = \sum_{c\in \lceil a,b\rfloor} \chxx^c
\end{eqnarray}
\end{thm}
\begin{proof}
From the Clebsch-Gordan decomposition,%
  $$(V_a\otimes V_b)^*\otimes (V_a\otimes V_b)\isom\sum_{\Adm ab{c,d}}V^*_c\otimes V_d.$$
Hence $$\End(V_a\tensor V_b)^G\isom\sum_{\Adm abc}\End(V_c)^G$$ and the characters satisfy
  $$\chxx^a\chxx^b=\chxx_{\sss(V_a\otimes V_b)}=\chxx_{\sss\oplus_cV_c}
    =\sum_{\Adm abc}\chxx^c.$$
There is an alternate diagrammatic proof of this statement, which uses the fusion and bubble
identities in Propositions \ref{bubble} and \ref{fusion}. If the matrix $\xb$ is represented by
$\ims{\bfx}$, then:\index{central function!product}
\begin{align*}
  \chxx^a\chxx^b&=\cfrelaa(\bfx,\bfx)(a,b)
                =\sum_{\Adm abc}\pfrac{\Delta(c)}{\Theta(a,b,c)}\cftwoc(\bfx,\bfx)(a,b,c)\\
              &=\sum_{\Adm abc}\pfrac{\Delta(c)}{\Theta(a,b,c)}\cfrelab(\bfx)(a,b,c)
                =\sum_{\Adm abc}\pfrac{\Delta(c)}{\Theta(a,b,c)}\cfrelac(\bfx)(a,b,c)\\
              &=\sum_{\Adm abc}\pfrac{\Delta(c)\Theta(a,b,c)}{\Theta(a,b,c)\Delta(c)}\cfone(\bfx)(c)%
                =\sum_{\Adm abc}\cfone(\bfx)(c)=\sum_{\Adm abc}\chxx^c.\qedhere
\end{align*}
\end{proof}

The product formula \eqref{rank1prod} and the initial calculations of $\chxx^0$ and $\chxx^1$ may
be used to show:

\begin{thm}
$\cb[G]^G\isom\cb[t]$.
\end{thm}
\begin{proof}
Consider the ring homomorphism $\Phi:\cb[t]\to \cb[G]^G$ defined by $f\mapsto f\circ \mathrm{tr}.$
Suppose $f(\tr(g))=0$ for all $g\in G$. If $f\neq0$, then since $f$ has a finite number of zeros,
$\tr(g)$ must have a finite number of values. However, $$\tmx{t&1}{-1&0}\in G$$ for all values of
$t$.  Hence, $f=0$ and $\Phi$ is injective.  It remains to establish surjectivity. We have already
shown $t\mapsto \chxx^1$ and $1\mapsto \chxx^0.$  Suppose $a\ge2$ and $\chxx^b$ is in the image of
$\Phi$ for all $b<a$.  Equation \eqref{rank1prod} implies $\chxx^1\chxx^{a-1}=\chxx^a+\chxx^{a-2}.$
Thus, by induction, $$t \Phi^{-1}(\chxx^{a-1})-\Phi^{-1}(\chxx^{a-2})\mapsto
\chxx^a.\mbox{\qedhere}$$
\end{proof}

The following closed formula for $\chxx^n$ is given in \cite{Pet}:
  $$\chxx^n(t)=\sum_{r=0}^{\lfloor\frac n2\rfloor}(-1)^r\binom{n-r}{r}t^{n-2r}.$$


The characters $\chxx^n$ may also be expressed as functions of eigenvalues, since $\chxx^n$ is
determined by its values on normal forms
$$\tmx{\lambda & *}{0 & \inv\lambda}\in G.$$ Explicitly, $\tmxt{\lambda & *}{0 & \inv\lambda}$ acts on
$V_n$ by the matrix%
{\small
$$\left[\begin{array}{ccccc}
\lambda^n& *&*&\cdots&*\\
0&\lambda^{n-2}&*&\cdots &*\\
\vdots&0&\ddots &*&*\\
0 &\vdots& 0 & \lambda^{2-n}&*\\
0&0&\cdots&0 &\lambda^{-n}
\end{array}\right].$$}%
Hence,
$$\chxx^n=\lambda^n+\lambda^{n-2}+\cdots+\lambda^{2-n}+\lambda^{-n}=\frac{\lambda^{n+1}-\lambda^{-n-1}}{\lambda-\lambda^{-1}}=[n+1]_\lambda,$$
where $[n+1]_\lambda$ is the quantized integer for $q=\lambda$.

%% file: spinsl2-rank2.tex

\section{Structure of $\C[G\times G]^G$} \label{rank2}

\index{central function!rank two}
Recall the decomposition $$\cb[G\times G]^G\isom\sum_{\substack{a,b\in\N\\
\Adm abc}}\cspan\ch{a}{b}{c},$$ where $\ch{a}{b}{c}$ corresponds by $\Upsilon$ to the image of
  $$\sum_{k=0}^c{\sf c}_k({\sf c}_k)^T\mapsto\sum_{k=0}^c\tbinom{c}{k}\bs{\sf c}{k}{k}$$
under the injection $V^*_c \otimes V_c\hookrightarrow V^*_a\otimes V^*_b \otimes V_a\otimes V_b$.
This inclusion is determined by the Clebsch-Gordan injection $\iota:V_c\hookrightarrow V_a\otimes
V_b.$ Hence, an explicit formula for $\iota$ provides a means to compute $\ch abc$ directly. We
freely use $\ch abc$ to denote its image in $\C[G\times G]^G$.

A few simple examples will motivate the construction of $\iota$. For $k=1,2$, let
$\xb_k=[x_{ij}^k]$ be $2\times2$ matrix variables, and let
\begin{align*}
  x&=\tr(\xb_1)=x^1_{11}+x^1_{22},\\
  y&=\tr(\xb_2)=x^2_{11}+x^2_{22},\\
  z&=\tr(\xb_1 \xb_2^{-1})=(x^1_{11}x^2_{22}+x^1_{22}x^2_{11})-(x^1_{12}x^2_{21}+x^1_{21}x^2_{12}).
\end{align*}
Recall that the map $\cupp:V_0\hookrightarrow V_1\otimes V_1$ given by
$${\sf c}_0\mapsto{\sf a}_1\otimes {\sf b}_0-{\sf a}_0\otimes {\sf b}_1$$ is invariant, using the notation defined in section \ref{reptheory}.
More generally, the injection $V_0 \hookrightarrow V_a\otimes V_a$ is given by
\begin{equation}\label{extcup}
  \tcup(a):{\sf c}_0\longmapsto\sum_{m=0}^{a}(-1)^m\tbinom{a}{m}{\sf a}_{a-m}\tensor{\sf b}_m.
\end{equation}
Hence, $\ch000=1$ and $\ch110$ may be computed by:
\begin{align*}
  \ch110 & \mapsto \bs{\sf c}{0}{0}\\
    & \mapsto ({\sf a}^*_1\otimes {\sf b}_0^*-{\sf a}_0^*\otimes {\sf b}^*_1)\otimes({\sf a}_1\otimes {\sf b}_0-{\sf a}_0\otimes {\sf b}_1)\\
    & \mapsto (\bs{\sf a}{1}{1})\otimes(\bs{\sf b}{0}{0})-(\bs{\sf a}{0}{1})\otimes (\bs{\sf b}{1}{0}) \\
        & \hspace{.5in} -(\bs{\sf a}{1}{0})\otimes (\bs{\sf b}{0}{1})+(\bs{\sf a}{0}{0})\otimes (\bs{\sf b}{1}{1})\\
    & \mapsto x^1_{11}\otimes x^2_{22}-x^1_{12}\otimes x^2_{21}-x^1_{21}\otimes x^2_{12}+x^1_{22}\otimes x^2_{11}\\
    & \mapsto (x^1_{11}x^2_{22}+x^1_{22}x^2_{11})-(x^1_{12}x^2_{21}+x^1_{21}x^2_{12})=z.
\end{align*}

The representation $V_c$ may be identified with a subset of $\vprod c$ via the equivariant maps
  $$\xymatrix{V_c\ar@/^1pc/[r]^-{\sss\sf Sym} & \vprod{c}\ar@/^1pc/[l]^-{\sss\sf Proj}}$$
where ${\sf Proj}\circ{\sf Sym}={\rm id}$. Thus, when $c=a+b$, $\iota$ is given by the commutative diagram
  $$\xymatrix{
     \vprod{c}\ar@{=}[r]\ar@{}[dr]|-{\text{\Large $\circlearrowright$}} & \vprod a\tensor\vprod b\ar[d]^{\sss\sf Proj\tensor Proj}\\
     V_c\ar[r]_-\iota\ar[u]^{\sss\sf Sym}              & V_a\tensor V_b.}$$
In particular,
\begin{equation}\label{mixsym}
  \tbinom{c}{k}{\sf c}_k\overset{\iota}{\longmapsto}\sum_{\substack{0\le i \le a\\ 0 \leq j \leq b \\
  i+j=k}}\tbinom{a}{i}{\sf a}_i\otimes \tbinom{b}{j}{\sf b}_j.
\end{equation}
For example, consider $\ch101$. In this case, ${\sf c}_0 \mapsto {\sf a}_0\otimes {\sf b}_0$ and
${\sf c}_1 \mapsto {\sf a}_1\otimes {\sf b}_0$. Hence, {\small
\begin{align*}
  \ch101&\mapsto\bs{\sf c}00+\bs{\sf c}11\mapsto(\bs{\sf a}{0}{0})\otimes(\bs{\sf b}00)+(\bs{\sf a}11)\otimes(\bs{\sf b}00)\\
   &\mapsto x^1_{11}\otimes 1+x^1_{22}\otimes 1\mapsto x^1_{11}+x^1_{22}=x.
\end{align*}}
A similar computation shows that $\ch011\mapsto y$.

The general form of $\iota$ is determined by combining \eqref{extcup} and \eqref{mixsym} in the
following diagram:
  $$\xymatrix{
      V_c\ar[r]^-\iota\ar[d]_-\iota\ar@{}[dr]|-{\text{\Large $\circlearrowright$}}
                                      & V_{\beta}\tensor V_{\alpha}\ar[d]^{{\rm id}\tensor{\tcup(\gamma)}\tensor{\rm id}}\\
      V_a\tensor V_b                  & V_{\beta}\tensor V_\gamma\tensor V_\gamma\tensor V_{\alpha}\ar[l].
  }$$
It follows that the mapping $\iota:V_c\to V_a\tensor V_b$ is explicitly given by:
\begin{align*}
  \tbinom{c}{k}{\sf c}_k
    &\longmapsto\sum_{\substack{0\le i\le\beta\\0\le j\le\alpha\\0\le m\le\gamma\\i+j=k}}%
      \tbinom{\beta}{i}{\sf a}_i
        \tensor\left[(-1)^m\tbinom{\gamma}{m}{\sf a}_{\gamma-m}\tensor{\sf b}_m\right]
        \tensor\tbinom{\alpha}{j}{\sf b}_j\\ %
    &\longmapsto\sum_{\substack{0\le i\le\beta\\0\le j\le\alpha\\0\le m\le\gamma\\i+j=k}}%
      (-1)^m\tbinom{\beta}{i}\tbinom{\alpha}{j}\tbinom{\gamma}{m}{\sf a}_{i+\gamma-m}\tensor{\sf b}_{j+m}.%
\end{align*}


\subsection{Symmetry of Central Functions}

Our first theorem regarding central functions is a symmetry property that is essentially trivial in diagram form, despite being highly nontrivial algebraically. A portion of the Fricke-Klein-Vogt Theorem (\ref{fricke}) is required to state the theorem. We begin with a diagrammatic proof of this classical result, in which the the binor identity plays the role of the characteristic equation in the classical proof.
\index{Fricke-Klein-Vogt Theorem}

\begin{lem}\label{nonconstructive}
	Each central function $\ch abc$ is associated to a unique polynomial $\mathfrak{p}_{\sss a,b,c}$,
	denoted for all pairs $(\bfx_1,\bfx_2)\in G\times G$ by
		$$\ch abc(\bfx_1,\bfx_2)=\mathfrak{p}_{\sss a,b,c}(\tr(\bfx_2),\tr(\bfx_1),\tr(\bfx_1\bfx_2^{-1})).$$
\begin{proof}
Expanding the symmetrizers in $\ch abc$ gives a collection of circles with matrix elements, each of which correspond to a product
of traces of words in $\bfx_1$ and $\bfx_2$, so it suffices to show that every loop can be reduced to a collection of loops containing one of $\bfx_1$, $\bfx_2$, or $\bfx_1\bfx_2^{-1}$.

This reduction depends entirely on the binor identity, which when composed with
$\bfx_1\tensor\bfx_2=\mbinorrela(\bfx_1,\bfx_2)$ gives:
	\begin{equation}\label{binorgh}
		\mbinor(\bfx_1,\bfx_2).
	\end{equation}%
Denote $\bfx_1^{-1}$ by $\bar\bfx_1$. Two special cases of \eqref{binorgh} follow:
	$$
		\mbinor(\bfx_1,\bar\bfx_1)=\mbinorrela(\bfx_1,\bar\bfx_1)-\mbinorrelb(\bfx_1^2)
		\qquad\text{and}\qquad
		\mbinor(\bfx_1,\bfx_1)=\mbinorrela(\bfx_1,\bfx_1)-\mbinorrelb(\mathbb{I}).
	$$%
The first relation allows us to assume no loop has both $\bfx_1$ and $\inv\bfx_1$, while the second
allows us to assume no loop has more than one of any matrix. The remaining cases are the traces
$\tr(\bfx_1)$, $\tr(\bfx_2)$, $\tr(\bfx_1\bfx_2)$, and $\tr(\bfx_1\inv{\bfx_2})$. Finally, closing off (\ref{binorgh})
permits the reduction of $\tr(\bfx_1\bfx_2)$:
$$\tr(\bfx_1\bfx_2)=\tr(\bfx_1)\tr(\bfx_2)-\tr(\bfx_1\inv{\bfx_2}).\qedhere$$
\end{proof}
\end{lem}

We can now prove the symmetry result. In the statement and proof below, $\sigma(\dmd_1,\dmd_2,\dmd_3)$ denotes the ordered triple $(\dmd_{\sigma(1)},\dmd_{\sigma(2)},\dmd_{\sigma(3)})$ obtained by applying a given permutation $\sigma\in\Sigma_3$ to the triple $(\dmd_1,\dmd_2,\dmd_3)$. This result was first outlined in \cite{Res}.

\begin{thm}[Symmetry of Central Functions]\label{symmetry}\index{central function!symmetry}
The family of polynomials $\ch abc(\bfx_1,\bfx_2)=\mathfrak{p}_{\sss a,b,c}(\tr(\bfx_2),\tr(\bfx_1),\tr(\bfx_1\bfx_2^{-1}))$ posseses the following symmetry:
	$$\mathfrak{p}_{\sss \sigma(a,b,c)}(y,x,z)=\mathfrak{p}_{\sss a,b,c}(\inv\sigma(y,x,z)).$$
\begin{proof}
Define the following function $G\cross G\cross G\to\C$:%
	$$\chx\alpha\beta\gamma(\bfx,\bfy,\bfz)=\tcfsyma(\bfx,\bfy,\bfz)(\alpha,\beta,\gamma)$$
where the symmetrizer on the right is assumed to `wrap around' to the one on the left (imagine this
diagram being drawn on a cylinder). By construction this function is symmetric, in the sense that:
	$$\chxx_{\sss\sigma(\alpha,\beta,\gamma)}\left(\sigma\left(\ims{\bfx_1},\ims{\bfx_2},\ims{\bfx_3}\right)\right)
		=\chx\alpha\beta\gamma\left(\ims{\bfx_1},\ims{\bfx_2},\ims{\bfx_3}\right).$$

A central function $\ch abc(\bfx_1,\bfx_2)$ may be drawn as:
	$$\cftwo(\bfx_1,\bfx_2)(a,b,c)=\tcfsymb(\bfx_1,\bfx_2)(\frac{a-b+c}{2},\frac{a+b-c}{2},\frac{-a+b+c}{2})
    =\tcfsyma(\bar\bfx_1,\bfx_1\bar\bfx_2,\bfx_2)(\beta,\gamma,\alpha),$$
with the symmetrizers in the last two diagrams assumed to wrap around as before. Thus, $\mathfrak{p}_{\sss
a,b,c}(y,x,z)=\chx\alpha\beta\gamma(\bfx_2,\inv{\bfx_1},\bfx_1\inv{\bfx_2})$ and so:
\begin{align*}
\mathfrak{p}_{\sss \sigma(a,b,c)}(y,x,z) &=\chxx_{\sss\sigma(\alpha,\beta,\gamma)}(\bfx_2,\inv{\bfx_1},\bfx_1\inv{\bfx_2})\\
                         &=\chx\alpha\beta\gamma(\inv\sigma(\bfx_2,\inv{\bfx_1},\bfx_1\inv{\bfx_2}))\\
                         &=\mathfrak{p}_{\sss a,b,c}(\inv\sigma(y,x,z)).\qedhere
\end{align*}
\end{proof}
\end{thm}
Table \ref{cfsymtable} contains six central functions illustrating this symmetry.
{\renewcommand\arraystretch{1.5}
  \begin{table}[!h]\label{cfsymtable}
    \begin{center}
      \begin{tabular}{|c|c|}
        \hline
        $\ch123=xy^2-\tfrac23(yz+x)$ & $\ch321=xz^2-\tfrac23(yz+x)$\\
        \hline
        $\ch231=yz^2-\tfrac23(xz+y)$ & $\ch132=y^2z-\tfrac23(xy+z)$\\
        \hline
        $\ch312=x^2z-\tfrac23(xy+z)$ & $\ch213=x^2y-\tfrac23(xz+y)$\\
        \hline
      \end{tabular}
    \medskip
    \caption{Rank Two Central Function Symmetry.}
    \end{center}
  \end{table}
}\index{central function!table}


\subsection{A Recurrence Relation for Central Functions}

Define the \emph{degree} of a central function to be:%
	$$\delta=\deg(\ch abc)=\tfrac12(a+b+c).$$
\index{central function!degree}
We will obtain a recurrence relation for an arbitrary central function $\ch abc$ by manipulating diagrams to express the product
$$\tr(\bfx_1)\cdot\ch abc(\bfx_1,\bfx_2)$$ as a sum of central functions. This formula can be rearranged to write $\ch abc$ as a
linear combination of central functions \emph{with lower degree}. There are three main ingredients to the diagram manipulations:
the \emph{bubble identity} and the \emph{fusion identity} from Section \ref{spinprops}, and two \emph{recoupling formulae} which
we prove in the following lemma.

\begin{lem} \label{6js} For $i=\half(a+1-b+c)$ and appropriate triples admissible,
	\begin{alignat}3\index{recoupling formula}
		\label{6j1}
			\thright(1,a,b,c,{\ssz c-1}) &=
			& &\tfused(1,a,b,c,{\ssz a+1})
			&-(-1)^i\ptfrac{a+b-c+1}{2(a+1)} &\tfused(1,a,b,c,{\ssz a-1});\\
		\label{6j2}
			\thright(1,a,b,c,{\ssz c+1})&=
			&(-1)^i\ptfrac{-a+b+c+1}{2(c+1)} &\tfused(1,a,b,c,{\ssz a+1})
			&+\ptfrac{(a+b+c+3)(a-b+c+1)}{4(a+1)(c+1)} &\tfused(1,a,b,c,{\ssz a-1}).%
	\end{alignat}
\begin{proof}
Note that $i$ is just the number of strands connecting $\suppnbig{a+1}$ to $\suppnbig{c}$ in $\tfused(1,a,b,c,{\ssz a+1})
= \tdowny(c,b,{\ssz a+1})$. For \eqref{6j1}, use $n=a+1$ and $i$ in recurrence relation \eqref{rec2} to
get:
$$\srrecb{a+1}{i}.$$
Compose this equation with $\cfdownpartial(c,b)({\ssz i},{\ssz a+1-i})$
to get, via the \emph{stacking relation}:
	$$\tfused(1,a,b,c,{\ssz a+1})=\tdowny(c,b,{\ssz a+1})
		=\thright(1,a,b,c,{\ssz c-1})+(-1)^i\pfrac{a+1-i}{a+1}\tfused(1,a,b,c,{\ssz a-1}),$$
which is the desired result.

To prove \eqref{6j2}, notice that if we switch $a$ and $c$ in the previous relation, and apply a $\frac\pi4$-reflection to the
relation about the $1\leftrightarrow b$ axis as in Proposition \ref{funsigns2}, then $i$ is unchanged and the equation becomes:
	$$\thright(1,a,b,c,{\ssz c+1})=\tfused(1,a,b,c,{\ssz a-1})+(-1)^i\pfrac{c+1-i}{c+1}\thright(1,a,b,c,{\ssz c-1}).$$
Rearrange this equation, and use \eqref{6j1} in its exact form to get:
	\begin{align*}
	\thright(1,a,b,c,{\ssz c+1})&=\tfused(1,a,b,c,{\ssz a-1})+(-1)^i\ptfrac{c+1-i}{c+1}\left(\tfused(1,a,b,c,{\ssz a+1})
		-(-1)^i\ptfrac{a+1-i}{a+1}\tfused(1,a,b,c,{\ssz a-1})\right)\\
    &=(-1)^i\ptfrac{c+1-i}{c+1}\tfused(1,a,b,c,{\ssz a+1})
    	+\left(1-\tfrac{(a+1-i)(c+1-i)}{(a+1)(c+1)}\right)\tfused(1,a,b,c,{\ssz a-1})\\
    &=(-1)^i\ptfrac{-a+b+c+1}{2(c+1)}\tfused(1,a,b,c,{\ssz a+1})+\ptfrac{(a+b+c+3)(a-b+c+1)}{4(a+1)(c+1)}\tfused(1,a,b,c,{\ssz a-1}).
	\end{align*}
To show the last computation, note that $a+1-i=\half(a+b-c+1)$ and $c+1-i=\half(-a+b+c+1)$, so the
numerator of the last term is:%
{\scriptsize\begin{align*}
4((a+1)(c+1)-(a+&1-i)(c+1-i))\\
 &= 4(a+1)(c+1)-((b+1)+(c-a))((b+1)-(c-a))\\
 &= 4(a+1)(c+1)-(b+1)^2+(a-c)^2\\
 &= ((a+1)-(c+1))^2+4(a+1)(c+1)-(b+1)^2\\
 &= ((a+1)+(c+1))^2-(b+1)^2\\
 &= (a+1+c+1+b+1)(a+1+c+1-b-1)\\
 &= (a+b+c+3)(a-b+c+1).\qedhere
\end{align*}}
\end{proof}
\end{lem}

The coefficients we have computed are examples of \emph{6j-symbols}\index{6j-symbols}, most easily defined to be the
coefficients $\sixjt abcdef'$ in the following \emph{change of basis} equation:
	$$\thright(a,b,c,d,{\ssz e})=\sum_{f\in\lceil{a,b}\rfloor\cap\lceil{c,d}\rfloor}\sixjt abcdef'\cdot\tfused(a,b,c,d,{\ssz f}).$$
We use a prime because we will need an alternate version later:
\begin{defn}
The \emph{6j-symbols} $\sixjt abcdef$ are the coefficients given by
	$$\tuptreer(a,b,c,d,{\ssz e})= \sum_{f\in\lceil{a,b}\rfloor\cap\lceil{c,d}\rfloor}\sixjt abcdef\cdot\tuptreel(a,b,c,d,{\ssz f}).$$
\end{defn}

Both versions given here differ from those in the literature \cite{CFS,Kau}. It is not hard to
show, using Corollary \ref{funsigns3}, that
	$$\sixjt abcdef'=(-1)^{\half(b+d-e-f)}\sixjt abcdef.$$
Thus, as a corollary to the above lemma we have the following $6j$-symbols, given by replacing $c$ with $c+1$ or $c-1$, which will be used to prove the next theorem:

\begin{cor} \label{6j3}
	\begin{alignat*}{3}
		\sixjt 1ab{c+1}c{a+1}&=1;\:
			&\sixjt 1ab{c+1}c{a-1}
			&=(-1)^{\half(a-b+c+2)}\tfrac{(a+b-c)}{2(a+1)};\\%
		\sixjt 1ab{c-1}c{a+1}&=(-1)^{\half(a-b+c)}\tfrac{(-a+b+c)}{2c};\:
			&\sixjt1ab{c-1}c{a-1}
			&=\tfrac{(a+b+c+2)(a-b+c)}{4(a+1)c}.%
	\end{alignat*}
\end{cor}

We can now prove the ``multiplication by $x$'' formula.
\begin{thm} \label{xmult} \index{central function!product}
The product $x\cdot\ch abc(x,y,z)$ can be expressed by:
\begin{multline}\label{multeq}
  x\cdot\ch abc = \ch {a+1}b{c+1} + \tfrac{(a+b-c)^2}{4a(a+1)} \ch {a-1}b{c+1}
    +\tfrac{(-a+b+c)^2}{4c(c+1)} \ch {a+1}b{c-1} \\
    + \tfrac{(a+b+c+2)^2(a-b+c)^2}{16a(a+1)c(c+1)}
      \ch {a-1}b{c-1}.
\end{multline}
This equation still holds for $a=0$ or $c=0$, provided we exclude the terms with $a$ or $c$ in the
denominator.
\begin{proof}
Diagrammatically, $x \cdot \ch abc(x,y,z)$ is represented by%
	$$\cfrelba(\bfx_1,\bfx_1,\bfx_2)(a,b,c,1),$$
since $x=\tr(\bfx_1)=\circml{\bfx_1}$ and multiplication is automatic on disjoint diagrams. Now
manipulate the diagram to obtain a sum over $\chi$'s with the following three steps.

First, apply the fusion identity to connect the lone $\ims{\bfx_1}$ strand to the
$\ch abc$:%
	\begin{equation}\label{multa}
  	\cfrelba(\bfx_1,\bfx_1,\bfx_2)(a,b,c,1)=\cfrelbb(\bfx_1,\bfx_1,\bfx_2)(a,b,c,c,c+1,1)
  		+\frac{c}{c+1}\cfrelbb(\bfx_1,\bfx_1,\bfx_2)(a,b,c,c,c-1,1),
	\end{equation}%
where the coefficients are evaluated from %
	$$
		\frac{\Delta(c \pm 1)}{\Theta(1,c,c \pm 1)}=\frac{c\pm1+1}{c+\frac32\pm\half}.
	$$

Second, use the $6j$-symbols computed in Corollary \ref{6j3} above to move the $a$ strand from one
side of the diagram to the other:
	\begin{alignat}{3}\label{multc}
  	\!\cfrelbb(\bfx_1,\bfx_1,\bfx_2)(a,b,c,c,c+1,1)\!\!\!&=
  		&\!\!\!\!\!\!\!\!&\cfrelbd(\bfx_1,\bfx_2)(a+1,b,c+1)\!\!\!
  		&+\tfrac{(a+b-c)^2}{4(a+1)^2}\!\!&\cfrelbc(\bfx_1,\bfx_1,\bfx_2)(a,b,a\!-\!1,a\!-\!1,c+1,1)\\
  	\!\cfrelbb(\bfx_1,\bfx_1,\bfx_2)(a,b,c,c,c-1,1)\!\!\!&=
  		&\tfrac{(-a+b+c)^2}{4c^2}\!\!\!\!\!\!\!\!&\cfrelbd(\bfx_1,\bfx_2)(a+1,b,c-1)\!\!\!
  		 &+\tfrac{(a+b+c+2)^2(a-b+c)^2}{16(a+1)^2c^2}\!\!&\cfrelbc(\bfx_1,\bfx_1,\bfx_2)(a,b,a\!-\!1,a\!-\!1,c-1,1)\!\!.%
	\end{alignat}%
In each case, we are recoupling twice: once for the top piece \tuptreer(,,,,) and once for the
corresponding bottom piece. In doing this, we would actually get four terms, but since the $a\pm1$
labels \emph{must be the same} on both the top and the bottom (a consequence of Schur's Lemma or
the bubble identity), two of the terms vanish.

In the final step, use the bubble identity to collapse the final pieces:
	\begin{align*}
  	\cfrelbc(\bfx_1,\bfx_1,\bfx_2)(a,b,{a\!+\!1},{a\!+\!1},{c\pm1},1)\!\!
  		 &=\ptfrac{\Theta(1,a,a+1)}{\Delta(a+1)}\!\!\!\!\!\!\!\cfrelbd(\bfx_1,\bfx_2)(a+1,b,c\pm1)\!\!=\ch{a+1}{b}{c\pm1};\\ %
  	\cfrelbc(\bfx_1,\bfx_1,\bfx_2)(a,b,{a\!-\!1},{a\!-\!1},{c\pm1},1)\!\!
  		 &=\ptfrac{\Theta(1,a,a-1)}{\Delta(a-1)}\!\!\!\!\!\!\!\cfrelbd(\bfx_1,\bfx_2)(a-1,b,c\pm1)\!\!=\ptfrac{a+1}{a}\ch{a-1}{b}{c\pm1}. %
	\end{align*}
At this point, obtaining \eqref{multeq} is simply a matter of multiplying the coefficients obtained
in the previous formulae.

Now consider the special cases. For $a=0$, since $b=c$ and consequently
$\frac{c}{c+1}=\frac{(-a+b+c)^2}{4c(c+1)}$, the desired formula is exactly \eqref{multa}.
Similarly, for $c=0$, the desired formula is \eqref{multc}.
\end{proof}
\end{thm}

We find it interesting that, for all our discussion of signs introduced by non-topological
invariance, all signs introduced are eventually squared and thus do not show up in this result.

We can rearrange the terms in \eqref{multeq} and re-index to get:
\begin{cor}[Central Function Recurrence]\label{recurrence}\index{central function!recurrence}
Provided $a>1$ and $c>1$, we can write
\begin{multline*}
  \ch abc=x\cdot\ch {a-1}{b}{c-1}-\tfrac{(a+b-c)^2}{4a(a-1)}\ch {a-2}bc%
    -\tfrac{(-a+b+c)^2}{4c(c-1)}\ch ab{c-2}\\
    -\tfrac{(a+b+c)^2(a-b+c-2)^2}{16a(a-1)c(c-1)}\ch{a-2}b{c-2}.
\end{multline*}
The relation still holds for $a=1$ or $c=1$, provided we exclude the terms with $a-1$ or $c-1$ in
the denominator.
\end{cor}
The condition $a>1,\:c>1$ arises because decrementing $a$ and $c$ in \eqref{multeq} means $(a-1,b,c-1)$ must now be admissible.
Also, note that formulae for multiplication by $y$ and $z$ may be obtained by applying the symmetry relation of Theorem
\ref{symmetry}. This fact is indispensable in our proof of Theorem \ref{fricke}.

\subsection{Graded Structure of the Central Function Basis}\index{central function!grading}
The majority of the content in this section was suggested to us by Carlos Florentino \cite{Flo} after he read an early draft of this chapter.

Recall the $\alpha,\beta,\gamma$ notation used earlier, and the notation
  $$\chx\alpha\beta\gamma(\bfx_2,\inv\bfx_1,\bfx_1\inv\bfx_2)=\ch abc(\bfx_1,\bfx_2)$$ %
introduced in the proof of Theorem \ref{symmetry}. The recurrence in Corollary \ref{recurrence} may be rewritten as
\begin{multline*}
\chx\alpha\beta\gamma=\chx010\chx\alpha{\beta-1}\gamma
    -\tfrac{\gamma^2}{a(a-1)}\chx{\alpha+1}{\beta-1}{\gamma-1}
    -\tfrac{\alpha^2}{c(c-1)} \chx{\alpha-1}{\beta-1}{\gamma+1}\\
    -\tfrac{\delta^2(\beta-2)^2}{a(a-1)c(c-1)}\chx\alpha{\beta-2}\gamma.
\end{multline*}
The interchangeability of $(a,\alpha)$ and $(c,\gamma)$ is guaranteed by the symmetry theorem.

\begin{prop}\label{monic}
The polynomial $\ch abc=\chx\alpha\beta\gamma$ is monic, with highest degree monomial $x^\beta
y^\alpha z^\gamma$.
\begin{proof}
Induct on the degree $\delta=\alpha+\beta+\gamma$ of central functions. The statement is clearly true
for the base cases, since $\chx000=1, \chx010=x, \chx100=y$, and $\chx001=z$.  The recurrence
relation implies that the highest order term of $\chx\alpha\beta\gamma$ is $x$ times the highest
order term of $\chx\alpha{\beta-1}\gamma$, hence $x(x^{\beta-1}y^\alpha z^\gamma)=x^\beta y^\alpha
z^\gamma$. This fact, together with the appropriate symmetric facts for $y$ and $z$, completes the
induction.
\end{proof}
\end{prop}

The basis also preserves a certain grading on $\C[x,y,z]$. To define this grading, partition the standard basis
$\mathcal{B}=\{x^ay^bz^c\}$ of this space as follows. Let $\gr:\mathcal{B}\to\Z_2\times\Z_2$ be defined by:
  $$\gr(x^ay^bz^c)=(a+c,b+c)\mo2.$$
If $\mathcal{B}$ is considered as a semigroup under multiplication, then $\gr$ is a homomorphism since
\begin{align*}
  \gr(x^ay^bz^c)+\gr(x^{a'}y^{b'}z^{c'})&=(a+c,b+c)+(a'+c',b'+c')\mo2\\
                                        &=(a+a'+c+c',b+b'+c+c')\mo2\\
                                        &=\gr(x^{a+a'}y^{b+b'}z^{c+c'})\mo2.
\end{align*}%
Therefore, $\gr$ defines a grading on this basis.

\begin{prop}\label{grading}
The basis $\{\ch abc\}$ respects the $\Z_2\times\Z_2$-grading on $\C[x,y,z]$ defined by $\gr$, in the sense that
	$$\ch abc\in\C\left(\inv\gr(a,b)\right).$$
\begin{proof}
This is another proof by induction on the degree $\delta$. Clearly, $\ch 000=1\in\inv \gr(0,0)$, and likewise $\ch 101=x\in\inv
\gr(1,0)$, $\ch 011=y\in\inv \gr(0,1)$, and $\ch 110=z\in\inv \gr(1,1)$. In the induction step, note that
  $$(a,b)=(1,0)+(a-1,b)=(a-2,b)\mo2,$$
so all terms on the righthand side of the recurrence relation in Corollary \ref{recurrence} have the same grading. Thus $\ch
abc\in\inv \gr(a,b)$.
\end{proof}
\end{prop}


\subsection{Multiplication of Central Functions} It is not difficult to write down the formula for
the product of two central functions, although the formula is by no means simple. The proof that follows was motivated by
\cite{Res}. We begin with a lemma which encapsulates the most tedious diagram manipulations:

\begin{lem}
	$$\trelba(a,b,c)(a',b',c')=\sum_{i,j,k,l,m}{\mathrm C}^{\sss abc,a'b'c'}_{\sss j_1k_1l_1,j_2k_2l_2,m}
		\trelbb(a,b)(a',b')(k_1,k_2,l_1,l_2,m),$$
where the coefficients are given by the formula%
\begin{equation*}
{\mathrm C}^{\sss abca'b'c'}_{\sss j_1k_1l_1,j_2k_2l_2,m}=
    \tfrac{\Theta(c,c',m)}{\Delta(m)}\prod_{i=1,2}
    \tfrac{\Delta(j_i)}{\Theta(a',b,j_i)}
    \cdot{\sixjt a{a'}{j_i}cb{k_i}}{\sixjt {b'}b{j_i}{c'}{a'}{l_i}}{\sixjt {k_i}{l_i}{c'}c{j_i}m},
\end{equation*}
and the following 15 triples are assumed to be admissible:%
{\begin{center}
  {$(a,a',k_i)$, $(b,b',l_i)$, $(c,c',m)$, $(a',b,j_i)$, $(c,j_i,k_i)$,
    $(c',j_i,l_i)$, $(b,j_i,l_i)$, $(k_i,l_i,m)$}.%
\end{center}}
\begin{proof}
We will just demonstrate the diagram manipulation for the top half of the diagram, which by
symmetry must be the same as for the bottom half. Combining these two manipulations and applying a
bubble identity will give the desired result. We will save enumeration of admissible triples until
after the manipulation, but keep a close eye on signs in the meantime.

{\small
	\begin{align*}
  	\trela(a,b,c)(a',b',c')
  	=\sum_{j}(-1)^{\half(a'-b+j)}\tfrac{\Delta(j)}{\Theta(a',b,j)}
    	&\trelb(a,b,c)(a',b',c')(j)\\
  	=\sum_{j,k}(-1)^{\half(a'-b+j)+j}\tfrac{\Delta(j)}{\Theta(a',b,j)}{\sixjt a{a'}jcbk}
    	&\trelc(a,b,c)(a',b',c')(j,k)\\
  	=\sum_{j,k,l}(-1)^{\half(a'-b-j)}\tfrac{\Delta(j)}{\Theta(a',b,j)}{\sixjt a{a'}jcbk}{\sixjt {b'}bj{c'}{a'}l}
      &\treld(a,b,c)(a',b',c')(j,k,l)\\
  	=\sum_{j,k,l}(-1)^{\half(a'-b-j)+\half(j+l-c')}\tfrac{\Delta(j)}{\Theta(a',b,j)}{\sixjt a{a'}jcbk}{\sixjt {b'}bj{c'}{a'}l}
      &\trele(a,b,c)(a',b',c')(j,k,l)\\
  	=\sum_{j,k,l,m}{(-1)^{\sss \half(a'-b+c-c'-j-m)+l}}\tfrac{\Delta(j)}{\Theta(a',b,j)}
    	{\sixjt a{a'}jcbk}{\sixjt {b'}bj{c'}{a'}l}{\sixjt kl{c'}cjm}%
      &\trelf(a,b,c)(a',b',c')(k,l,m)\\
	\end{align*}}

The $(-1)$ terms all cancel in the end, a consequence of the fact that the following triples must
be admissible:
\begin{center}
  {\small $(a,a',k)$, $(b,b',l)$, $(c,c',m)$, $(a',b,j)$, $(c,j,k)$,
    $(c',j,l)$, $(b,j,l)$, $(k,l,m)$.}
\end{center}
One computes the 13-parameter coefficients ${\mathrm C}^{\sss abc,a'b'c'}_{\sss
j_1k_1l_1,j_2k_2l_2,m}$ above by reflecting this result vertically, taking two sets of indices for
the variables $j,k,l,m$ on the two halves, and noting that the resulting bubble in the middle
collapses with a factor of $\tfrac{\Theta(c,c',m)}{\Delta(m)}$ for $m=m_1=m_2$.
\end{proof}
\end{lem}

With that out of the way, we can describe the central function multiplication table explicitly.
Note the symmetry with respect to $k,l,m$, which is guaranteed by Theorem \ref{symmetry}.

\begin{thm}[Multiplication of Central Functions]\label{mult}\index{central function!product}
The product of two central functions $\ch abc$ and $\ch {a'}{b'}{c'}$ is given by:
$$\ch abc\ch{a'}{b'}{c'}=\sum_{j_1,j_2,k,l,m}{\mathrm C}_{j_1klm} {\mathrm C}_{j_2klm}
        \tfrac{\Theta(a,a',k)\Theta(b,b',l)\Theta(c,c',m)}{\Delta(k)\Delta(l)\Delta(m)}\ch klm,$$
where the sum is taken over admissible triples
  $$\text{\small $(a,a',k)$, $(b,b',l)$, $(c,c',m)$, $(a',b,j_i)$, $(c,j_i,k)$,
    $(c',j_i,l)$, $(b,j_i,l)$, $(k,l,m)$}$$
and the coefficients are given by:
$${\mathrm C}_{j_iklm}=\tfrac{\Delta(j_i)}{\Theta(a',b,j_i)}{\sixjt a{a'}{j_i}cbk}{\sixjt {b'}b{j_i}{c'}{a'}l}{\sixjt kl{c'}c{j_i}m}.$$%
\begin{proof}
By the previous lemma and the \emph{bubble identity}, we have:
\begin{align*}
  \cfrelca(\bfx_1,\bfx_2)(a,b,c)(a',b',c')
    &=\sum_{j_1,k_1,l_1,j_2,k_2,l_2,m}{\mathrm C}^{abc,a'b'c'}_{j_1k_1l_1,j_2k_2l_2,m}
    	\cfrelcb(\bfx_1,\bfx_2)(a,b)(a',b')(k,l,m)\\ %
    &=\sum_{j_1,j_2,k,l,m}{\mathrm C}^{abc,a'b'c'}_{j_1kl,j_2kl,m}\left(\frac{\Theta(a,a',k)\Theta(b,b',l)}{\Delta(k)\Delta(l)}\right)
    	\cfrelbd(\bfx_1,\bfx_2)(k,l,m)\\
    &=\sum_{i,j,k,l}{\mathrm C}_{j_1klm}{\mathrm C}_{j_2klm}
    	\tfrac{\Theta(a,a',k)\Theta(b,b',l)\Theta(c,c',m)}{\Delta(k)\Delta(l)\Delta(m)}
      \cfrelbd(\bfx_1,\bfx_2)(k,l,m).
    \qedhere
\end{align*}
\end{proof}
\end{thm}


\subsection{Applications}

Spin networks offer a novel approach to a classical theorem of Fricke, Klein, and Vogt \cite{FK,Vogt}.
We give here a new \emph{constructive} proof which depends on the symmetry, recurrence, and multiplication formulae for central functions.

\begin{thm}[Fricke-Klein-Vogt Theorem]\label{fricke}\index{Fricke-Klein-Vogt Theorem}
Let $G=\sltc$ act on $G\times G$ by simultaneous conjugation.  Then $$\mathbb{C}[G\times G]^G \isom \mathbb{C}[t_x,t_y,t_z],$$
the complex polynomial ring in three indeterminates. In particular, every regular function $f:\sltc\times\sltc\to\C$ satisfying%
$$f(\bfx_1,\bfx_2)=f(g\bfx_1g^{-1},g\bfx_2g^{-1})\qquad\text{for all }g\in\sltc,$$%
can be written uniquely as a polynomial in the three trace variables $x=\tr(\bfx_1)$, $y=\tr(\bfx_2)$, and
$z=\tr(\bfx_1\inv{\bfx_2})$.
\begin{proof}
Define the ring homomorphism $$\Gamma:\cb[t_x,t_y,t_z]\to \cb[G\times G]^G$$ by
$f(t_x,t_y,t_z)\mapsto f(\tr(\xb_1),\tr(\xb_2),\tr(\xb_1\xb_2^{-1})).$

We first show that $\Gamma$ is injective. Suppose $f(\tr(\xb_1),\tr(\xb_2),\tr(\xb_1\xb_2^{-1}))
=0$ for all pairs $(\xb_1,\xb_2)\in G\times G$.  Let $(\tau_x,\tau_y,\tau_z)\in \cb^3$,
$\epsilon_x=\tmx{\tau_x&1}{-1&0}$, and $\eta_{y,z}=\tmx{\tau_y&\frac{1}{\zeta}}{-\zeta&0},$ where
$\zeta+\zeta^{-1}=\tau_z$.  Then
$$(\tau_x,\tau_y,\tau_z)=(\tr(\epsilon_x),\tr(\eta_{y,z}), \tr(\epsilon_x \eta_{y,z}^{-1})).$$  Hence $f=0$
on $\cb^3$, $\mathrm{Ker}(\Gamma)=\{0\}$, and $\Gamma$ is injective. This is the ``Fricke slice'' given by Goldman in
\cite{Gol3}.

It remains to show that $\Gamma$ is surjective. Theorem \ref{pwdecomp} implies that the central
functions form a basis for $\C[G\times G]^G$. Since $t_x\mapsto x,t_y\mapsto y$, and $t_z\mapsto
z$, it suffices to show that every $\ch abc$ may be written as a polynomial in $x,y$, and $z$. This was already done via Lemma \ref{nonconstructive}, but we provide here a \emph{constructive} proof.

Proceed by induction on the degree $\delta=\half(a+b+c)$ of a central function $\ch abc$. For the base cases $\delta=0,1$ recall
our earlier computations demonstrating
	$$\ch000=1,\;\ch101=x,\;\ch011=y,\;\ch110=z.$$

For $\delta>0$, we may inductively assume that all central functions with degree less than $\delta$ are in $\C[x,y,z]$. The
admissibility conditions imply that at least two out of the triple $(a,b,c)$ are positive. Without loss of generality, using
Theorem \ref{symmetry}, we may assume that $a$ and $c$ are positive. In this case, the recurrence given by Corollary
\ref{recurrence},
\begin{multline*}
\ch abc=x\cdot\ch{a-1}b{c-1}-\tfrac{(a+b-c)^2}{4a(a-1)}\ch{a-2}bc\\
    -\tfrac{(-a+b+c)^2}{4c(c-1)}\ch ab{c-2}
    -\tfrac{(a+b+c)^2(a-b+c-2)^2}{16a(a-1)c(c-1)}\ch{a-2}b{c-2},
\end{multline*}
allows us to write $\ch abc$ in terms of central functions of lower degree, which by induction must
be in $\C[x,y,z]$. Thus, $\ch abc\in\C[x,y,z]$, and we have established surjectivity.
\end{proof}
\end{thm}

The recursion relations provide an algorithm for writing any $\ch abc$ as a polynomial in $\{x,y,z\}$. Conversely, in
\cite{Pet} the following formula is established, which may be used to express any polynomial in $\C[x,y,z]$ in terms of central
functions:
  \begin{multline*}
    x^ay^bz^c=\hspace{-11pt}\sum_{\substack{r,s,t=0\\ k,l,m}}^{\lfloor\frac a2\rfloor,\lfloor\frac b2\rfloor,\lfloor\frac c2\rfloor}
      \hspace{-11pt}\Big(\tbinom ar-\tbinom a{r-1}\Big)\Big(\tbinom bs-\tbinom b{s-1}\Big)\Big(\tbinom ct-\tbinom c{t-1}\Big)\cdot\\
      \Big(\tfrac{\Delta(l)\Delta(m)\Theta(a-2r,c-2t,k)}{\Delta(k)\Theta(a-2r,b-2s,m)\Theta(b-2s,c-2t,l)}\Big)
      \sixjt{\!a-2r}{c-2t}{l}{m}{b-2s\!}{k}^2
      \ch klm.
  \end{multline*}

Table \ref{cftwotable} lists several central functions that were computed with Mathematica using Corollary \ref{recurrence}. Only
one function per triple of indices is listed; the others follow directly from Theorem \ref{symmetry}.

{\renewcommand\arraystretch{1.5}
  \begin{table}[ht]\label{cftwotable}
    \begin{center}
      \begin{tabular}{|c|c|c||l|}
        \hline
        {\bf $\delta$} & $\ch abc$ & $\chx\alpha\beta\gamma$ & $\mathfrak{p}_{\sss a,b,c}(y,x,z)$\\
        \hline\hline
        $\bf 0$ & $\ch000$ & $\chx000$ & $1$\\
        \hline
        $\bf 1$ & $\ch101$ & $\chx010$ & $x$\\
        \hline
        $\bf 2$ & $\ch202$ & $\chx020$ & $x^2-1$\\
                & $\ch112$ & $\chx110$ & $xy-\half z$\\
        \hline
        $\bf 3$ & $\ch303$ & $\chx030$ & $x^3-2x$\\
                & $\ch213$ & $\chx120$ & $x^2y-\tfrac23(xz+y)$\\
                & $\ch222$ & $\chx111$ & $xyz-\half(x^2+y^2+z^2)+1$\\
        \hline
        $\bf 4$ & $\ch404$ & $\chx040$ & $x^4-3x^2+1$\\
                & $\ch314$ & $\chx130$ & $x^3y-\tfrac34x^2z-\half(3xy-z)$\\
                & $\ch224$ & $\chx220$ & $x^2y^2-xyz+\tfrac16z^2-\half(x^2+y^2)+\tfrac13$\\
                & $\ch323$ & $\chx121$ & $x^2yz-\tfrac23(xz^2+xy^2)-\half x^3-\tfrac19(2yz-13x)$\\
        \hline
      \end{tabular}
    \medskip
    \caption{$\sltc$-Central Functions.}
    \end{center}
  \end{table}
}\index{central function!table}